\input amstex
\documentstyle{amams} 
\document
\input amssym.def
\input amssym.tex

\input boxedeps.tex 
\SetepsfEPSFSpecial 
\HideDisplacementBoxes

\annalsline{153}{2001}
\received{February 22, 1999}
\startingpage{191}
\font\tenrm=cmr10

\catcode`\@=11
\font\twelvemsb=msbm10 scaled 1100

\font\ninemsb=msbm10 scaled 800
\newfam\msbfam
\textfont\msbfam=\twelvemsb  \scriptfont\msbfam=\ninemsb
  \scriptscriptfont\msbfam=\ninemsb
\def\msb@{\hexnumber@\msbfam}
\def\Bbb{\relax\ifmmode\let\next\Bbb@\else
 \def\next{\errmessage{Use \string\Bbb\space only in math
mode}}\fi\next}
\def\Bbb@#1{{\Bbb@@{#1}}}
\def\Bbb@@#1{\fam\msbfam#1}
\catcode`\@=12

 \catcode`\@=11
\font\twelveeuf=eufm10 scaled 1100
\font\teneuf=eufm10
\font\nineeuf=eufm7 scaled 1100
\newfam\euffam
\textfont\euffam=\twelveeuf  \scriptfont\euffam=\teneuf
  \scriptscriptfont\euffam=\nineeuf
\def\euf@{\hexnumber@\euffam}
\def\frak{\relax\ifmmode\let\next\frak@\else
 \def\next{\errmessage{Use \string\frak\space only in math
mode}}\fi\next}
\def\frak@#1{{\frak@@{#1}}}
\def\frak@@#1{\fam\euffam#1}
\catcode`\@=12

\def\sni#1{\smallbreak\noindent{#1}}
\def\ssni#1{\vglue1pt\noindent\hskip12pt {#1}}


 \def\ve{\varepsilon} \def\lr{\leftrightarrow}  
 
\font\sans=cmss10 \font\script=eusm10 \def\V{\text{\script V}_s}
 \def\R{\Bbb R}  \def\C{\Bbb C}  \def\Cc{\C\,'} \def\Z{\Bbb Z}   \def\N{\Bbb N}  \def\bs{\gtrless}
\def\H{\Cal H} \def\G{\Gamma} \def\D{\Delta}\def\ve{\varepsilon} \def\h{\tfrac12}  \def\sm{\smallsetminus}  
\def\e{\doteq} \def\M{\text{\sans Maass}}  \def\FE{\text{\sans FE}_s}  
\def\lr{\leftrightarrow} \def\inv{^{-1}}   
\def\PL{Phragm\'en-Lindel\"of } \def\wf{\widetilde f}  \def\wU{\widetilde U}   \def\s{\sigma} \def\a{\alpha}
\def\p{\partial}
\def\bz{\bar z}
\def\z{\zeta}
\def\abcd{\bigl(\smallmatrix a&b\\c&d\endsmallmatrix\bigr)}
\def\Cs{\Cal C_s}
\def\F{\digamma}
\def\vp{\varphi}
\title{Period functions for Maass wave forms. I} 
   \twoauthors{J\. Lewis}{D\. Zagier}
   \institutions{Massachusetts Institute of Technology, Cambridge, MA\\
 {\eightpoint {\it E-mail address\/}: jlewis@math.mit.edu}\\
\vglue6pt
Max-Planck-Institut f\"ur Mathematik, Bonn, Germany\\
 {\eightpoint {\it E-mail address\/}: zagier@mpim-bonn.mpg.de}}
  
\bigbreak\centerline{\bf Contents}
\medbreak
{\ninepoint

\sni{Introduction}
\sni{Chapter I.  The period correspondence via $L$-series}

\ssni{1. The correspondences $u\lr L_\ve \lr f\lr\psi$} 
\ssni{2. Periodicity, $L$-series, and the three-term functional equation}

\ssni{3. Even and odd}

\ssni{4. Relations between Mellin transforms; proof of Theorem 1}

\sni{Chapter II. The period correspondence via integral transforms} 

\ssni{1. The integral representation of $\psi$ in terms of $u$}

\ssni{2. The period function as the integral of a closed 1-form}

\ssni{3. The incomplete gamma function expansion of $\psi$}

\ssni{4. Other integral transforms and intermediate functions}

\ssni{5. Boundary values of Maass wave forms}

\sni{Chapter III. Periodlike functions} 

\ssni{1. Examples}

\ssni{2. Fundamental domains for periodlike functions}

\ssni{3. Asymptotic behavior of smooth periodlike functions}

\ssni{4. ``Bootstrapping"}

\sni{Chapter IV. Complements} 

 \ssni{1. The period theory in the noncuspidal case}

\ssni{2. Integral values of $s$ and connections with holomorphic modular forms}

\ssni{3. Relation to the Selberg zeta function and Mayer's theorem}

\sni{References} 

}
 \bigbreak
\intro

Recall that a Maass wave form 
\footnote{$^1$}{We use the traditional term, but one should really specify ``cusp form." Also, the word ``form" should more
properly be  ``function," since $u$ is simply invariant under $\G$, with no automorphy factor. We often abbreviate ``Maass wave
form" to ``Maass form."}   on the full modular group $\G={\rm PSL}(2,\Z)$ is a smooth $\G$-invariant function
$u$ from the upper half-plane $\H=\{x+iy\mid y>0\}$ to $\C$ which is small as $y\to\infty$ and satisfies $\D u=\lambda\,u$ 
for some $\lambda\in\C$, where $\D=-y^2\bigl(\frac{\partial^2}{\partial x^2}+\frac{\partial^2}{\partial y^2}\bigr)$ is
the hyperbolic Laplacian. These functions give a basis for $L_2$ on the modular surface $\Gamma \backslash \H$,
in analogy with the usual trigonometric waveforms on the torus
$ \R^2 / \Z^2 $, which are also (for this surface) both the Fourier building 
blocks for $L_2$ and eigenfunctions of the Laplacian. Although therefore
very basic objects, Maass forms nevertheless still remain mysteriously 
elusive fifty years after their discovery; in particular, no explicit 
construction exists for any of these functions for the full modular group. 
The basic information about them (e.g. their existence and the density of 
the eigenvalues) comes mostly from the Selberg trace formula; the rest is
conjectural with support from extensive numerical computations. 

  Maass forms arise naturally in such diverse fields as number theory, 
dynamical systems and quantum chaos; hence concrete analytic information 
about them would be of interest and have applications in a number of areas of 
mathematics, and for this reason they are still under active investigation.

In [10], it was shown that there exists a one-to-one correspondence between the space of even Maass wave forms
(those with $u(-x+iy)=u(x+iy)$) with eigenvalue $\lambda=s(1-s)$ and the space of holomorphic functions 
on the cut plane $\Cc=\C\sm(-\infty,0]$  satisfying the functional equation
  $$\psi(z)=\psi(z+1)+z^{-2s}\,\psi\bigl(\frac {z+1}z\bigr)\qquad(z\in\Cc)\tag0.1$$
together with a suitable growth condition.  However, the passage from $u$ to $\psi$ was given in [10] by an integral transform (eq.~(2.2)
below) from which the functional equation~(0.1) and other properties of $\psi$ are not at all evident.  In the present paper we will:
\medbreak
 \item{(i)} find a different and simpler description of the function $\psi(z)$ and give a more conceptual proof
of the $u\lr\psi$ correspondence, for both even and odd wave forms, in terms of $L$-series; 
\item{(ii)} give a natural interpretation of the integral representation in [10] and of some of the related functions
 introduced there; 
\item{(iii)} study the general properties of the solutions of the functional equation (0.1), and determine sufficient conditions 
for such a solution to correspond to a Maass wave form;
\item{(iv)} show that the function $\psi(z)$ associated to a Maass wave form is the analogue   
 of the classical Eichler-Shimura-Manin period polynomial of a holomorphic cusp form, and describe the relationship between the two theories; and
\item{(v)} relate the theory to D.~Mayer's formula [13] expressing the Selberg zeta
 function of $\G$ as the Fredholm determinant of a certain operator on a space of holomorphic functions. \medbreak\noindent
Because of the connection (iv), we call $\psi$ the {\it period function} associated to the wave form $u$ and the mapping
$u\lr\psi$ the {\it period correspondence}.

In this introduction, we state the correspondence in its simplest form (e.g\. by assuming that the spectral parameter $s$
has real part $\h$, which is true for Maass forms anyway) and discuss some of its salient aspects.  In the statement below,
the functional equation (0.1) has been modified in two ways. First, the last term in the equation has been replaced by 
$(z+1)^{-2s}\,\psi\bigl(\dfrac z{z+1}\bigr)$; this turns out to give the functional equation which corresponds to arbitrary,
rather than just even, Maass forms. Secondly, we consider solutions of the functional equation on just the positive real axis, since
it will turn out that under suitable analytic conditions any such solution will automatically extend to $\Cc$.   

\phantom{hard}

\nonumproclaim{Theorem} Let $s$ be a complex number with $\Re(s)=\h${\rm .}
 Then there is an isomorphism between the space of
Maass cusp forms with eigenvalue $s(1-s)$ on $\G$ and the space of real\/{\rm -}\/analytic solutions of 
the {\rm three-term functional equation}
$$\psi(x)=\psi(x+1)+(x+1)^{-2s}\,\psi\bigl(\frac x{x+1}\bigr)\tag 0.2$$    
on $\R_+$ which satisfy the growth condition
   $$\psi(x)=\text{\rm o}(1/x)\qquad(x\to0),\qquad\quad \psi(x)=\text{\rm o}(1)\qquad(x\to\infty)\,.\tag0.3$$  \endproclaim
 
\vglue9pt 

In particular, the dimension of the space  of solutions of (0.2) satisfying the growth condition (0.3) is finite
for any $s$ and is zero except for a discrete set of values of~$s$. This is especially striking since we will show 
in Chapter III that if we relax the growth condition (0.3) minimally to
   $$ \psi(x)=\text O(1/x)\qquad(x\to0),\qquad\quad \psi(x)=\text O(1)\qquad(x\to\infty)\,,\tag0.4$$
then the space of all real-analytic solutions of (0.2) on $\R_+$ is infinite-dimen\-sional for any $s$. In the opposite direction,
however, if we weaken the growth condition further or even drop it entirely, then the space 
of solutions does not get any bigger, since we will show later that {\it any} real-analytic solution of (0.2)
on the positive real axis satisfies (0.4). Conversely, if $\psi$ does correspond to a Maass form, then we will see that (0.3)
can be strengthened to $\;\psi(x)=\text O(1)$ as $x\to0\,$ and $\;\psi(x)=\text O(1/x)$ as $x\to\infty$; moreover, $\psi$ in this case
extends holomorphically to $\Cc$ and these stronger asymptotic estimates hold uniformly in wedges $|\arg(x)|<\pi-\ve$.
%

The theorem formulated above is a combination of two results, Theorems~1~and~2, from the main body of 
the paper.  The first of these, whose statement and proof occupy Chapter I, gives a very
simple way to go between a Maass wave form $u$ and a solution of (0.2) holomorphic in the whole cut plane $\Cc$.  
We first associate to $u$ a {\it periodic and holomorphic} function
$f$ on $\C\sm\R$ whose Fourier expansion is the same as that of the Maass wave form $u$, but with the Bessel
functions occurring there replaced by exponential functions. Then we define $\psi$ in $\C\sm\R$ by the equation
$$c(s)\,\psi(z)\,=\,f(z)-z^{-2s}\,f\bigl(\frac{-1}z\bigr)\,,\tag0.5$$
where $c(s)$ is a nonzero normalizing constant. The three-term functional equation for $\psi$ is then a purely
algebraic consequence of the periodicity of $f$, while the fact that $f$ came from a $\G$-invariant function $u$ 
is reflected in the analytic continuability of $\psi$ across the positive real axis.  In the  converse direction, if $\psi(z)$ 
is any holomorphic function in $\Cc$, then the function $f$ defined by
$$c^\star(s)\,f(z)\,=\,\psi(z)+z^{-2s}\,\psi\bigl(\frac{-1}z\bigr)\tag0.6$$
(where $c^\star(s)$ is another normalizing constant)
is automatically periodic if $\psi$ satisfies the three-term functional equation, and corresponds to a Maass wave
form if $\psi$ also satisfies certain growth conditions at infinity and near the cut.  The proof of the correspondence in
both directions makes essential use of the Hecke $L$-series of~$u$ and of Mellin transforms, the estimates on $f$ and $\psi$
permitting us to prove the required identities by rotating the line of integration.  The key fact is that the {\it same} Hecke
$L$-series can be represented as the Mellin transform of {\it either} $u$ on the positive imaginary axis {\it or} $\psi$ 
on the positive real axis, but with different gamma-factors in each case.  The functional equation of the $L$-series is 
then reflected both in the $\G$-invariance of $u$ and in the functional equation (0.2) of $\psi$ on the positive reals.

In Chapter II we study the alternative construction of $\psi$ as an integral transform of $u$ (or, in the odd case, of its 
normal derivative) on the imaginary axis.  As already mentioned, this was the construction originally given in [10], but 
here we present several different points of view which make its properties more evident: we show that the function defined by the
integral transform has the same Mellin transform as the function constructed in Chapter I and hence agrees with it; 
re-interpret the integral as the integral of a certain closed form in the upper half-plane; give an expansion of the
integral in terms of the Fourier coefficients of $u$; construct an auxiliary entire function $g(w)$ whose special values at integer
arguments are the Fourier coefficients of $u$ and whose Taylor coefficients at $w=0$ are proportional to those of $\psi(z)$ at $z=1$;
and finally give a very intuitive description in terms of formal ``automorphic boundary distributions" defined by the limiting
behavior of $u(x+iy)$ as $y\to0$.

In Chapter III we study the general properties of solutions of the three-term functional equations (0.1) and (0.2).  We
start by giving a number of examples of solutions of these equations which do not necessarily satisfy the growth conditions and hence
need not come from Maass wave forms.  These are constructed both by explicit formulas and by a process analogous to the use of
fundamental domains in constructing functions invariant with respect to a group action.  Next, we study the properties
of smooth solutions of (0.2) on the real line and prove that they automatically satisfy the estimate (0.4).  The most amusing aspect
(Theorem 2, stated and proved in \S4) is a surprising ``bootstrapping" phenomenon which
says that any analytic solution of (0.2) on $\R_+$ satisfying the weak growth condition (0.3) (or a slight strengthening
of it if we do not assume that $\Re(s)=\frac12$) {\it automatically} extends to all of $\Cc$ as a holomorphic function 
satisfying the growth conditions required to apply Theorem~1. This provides the second key ingredient of the theorem formulated above.

In Chapter IV we return to the modular world. We treat three topics. The first is the extension of the theory to the noncuspidal case.
We find that the correspondence $u\lr\psi$ remains true if $u$ is allowed to be a noncuspidal Maass wave form and the growth condition on
$\psi$ at infinity is replaced by a weaker asymptotic formula.  The second topic is the relation to the classical holomorphic theory.
One of the consequences of the formulas from Chapter II is that the Taylor coefficients of $\psi(z)$ around $z=0$ and $z=\infty$ are 
proportional to the values at integral arguments of the Hecke-Maass $L$-function associated to $u$.  This is just like the correspondence 
between cusp forms and their period polynomials in the holomorphic case, where the coefficients of the period polynomial are multiples
of the values at integral arguments in the critical strip of the associated Hecke $L$-functions.  
In fact, it turns out that the theory developed in this paper and the classical theory of period polynomials are not only analogous, 
but are closely related when $s$ takes on an integral value.
 
Finally, the whole theory of period functions of Maass wave forms has a completely
different motivation and explanation coming from the work of Mayer [13] expressing the Selberg zeta function of $\G$
as the Fredholm determinant of a certain operator on a space of holomorphic functions: the numbers $s$ for which there
is a Maass form with eigenvalue $s(1-s)$ are zeros of the Selberg zeta function, and the holomorphic functions which are
fixed points of Mayer's operator satisfy the three-term functional equation (with shifted argument).  This point of view is described 
in the last section of Chapter IV and was also discussed in detail in the expository paper [11].

This concludes the summary of the contents of the present paper. Several of its main results were announced in [11],
which is briefer and more expository than the present paper, so that the reader may wish to consult it for an overview.
In the planned second part of the present paper we will discuss further aspects of the theory. In particular, we will
describe various ways to realize the period functions as cocycles, generalizing the classical interpretation of period
polynomials in terms of Eichler cohomology.  (These results developed partly from discussions with
Joseph Bernstein.) We will also treat a number of supplementary topics such as Hecke operators, Petersson scalar product,
extension to congruence subgroups of ${\rm SL}(2,\Z)$, and numerical aspects, and will discuss some arithmetic and nonarithmetic examples.


\demo{Acknowledgment} The first author would like to thank the Max-Planck-Institut f\"ur Mathematik, Bonn,
for its continued support and hospitality while this paper was being written.

\bigbreak\centerline{
\bf Chapter I. The period correspondence via $L$-series}
\bigbreak

In this chapter we state and prove the correspondence between Maass forms and solutions of the three-term functional equation in
the cut plane $\Cc$. The easier and more formal parts of the proof will be given in Sections 2 and 3.  The essential analytic step
of the proof, which involves relating both the Maass wave form $u$ and its period function $\psi$ to the $L$-series
of $u$ via Mellin transforms and then moving the path of integration, will be described in Section 4.

\section{The correspondences $u\lr L_\ve \lr f\lr\psi$} 

The following theorem gives four equivalent descriptions of Maass forms, the first equivalence $u\lr(L_0,L_1)$ being due to Maass.

\nonumproclaim{Theorem 1} Let $s$ be a complex number with $\s:=\Re(s)>0${\rm .}
 Then there are canonical correspondences between
objects of the following four types\/{\rm : }  
\medbreak
 {\rm (a)} a Maass cusp form $u$ with eigenvalue $s(1-s)${\rm ;}
\smallbreak
 {\rm (b)} a pair of Dirichlet $L$-series $L_\ve(\rho)\quad(\ve=0,\,1)${\rm ,} convergent in some
right half\/{\rm -}\/plane{\rm ,} such that the functions $L_\ve^*(\rho)=\gamma_s(\rho+\ve)\,L_\ve(\rho)${\rm ,} where
$$ \gamma_s(\rho)=\frac1{4\pi^\rho}\,\G\bigl(\frac{\rho-s+\h}2\bigr)\,\G\bigl(\frac{\rho+s-\h}2\bigr)\,,$$  
are entire functions of finite order and satisfy
    $$L^*_\ve(1-\rho)=(-1)^\ve L^*_\ve(\rho)\,;\tag1.1$$

 {\rm (c)} a holomorphic function $f(z)$ on $\C\sm\R${\rm ,}
 invariant under $z\mapsto z+1$\break and bounded by $|\Im(z)|^{-A}$ 
 for some $A>0$, such that the function\break $f(z)-z^{-2s}\,f(-1/z)$ extends holomorphically across the positive real axis and
is bounded by a multiple of $\,\min\{1,|z|^{-2\s}\}\,$ in the right half-plane{\rm ;}
\smallbreak
 {\rm (d)} a holomorphic function $\psi:\Cc\to\C$ satisfying the functional equation 
$$\psi(z)=\psi(z+1)+(z+1)^{-2s}\,\psi\bigl(\frac z{z+1}\bigr)\qquad(z\in\Cc)\tag1.2$$ 
and the estimates
$$\psi(z)\ll\cases |\Im(z)|^{-A}\,(1+|z|^{2A-2\s})&\text{if $\Re(z)\le0\,,$}\\ 1&\text{if $\Re(z)\ge0,\;|z|\le1\,,$}\\ 
  |z|^{-2\s}&\text{if $\Re(z)\ge0,\;|z|\ge1$} \endcases\tag1.3 $$
for some $A>0${\rm .}
\medbreak\noindent The correspondences $u\mapsto L_\ve$ and $f\mapsto L_\ve$ are given by the Mellin transforms
$$L_\ve^*(\rho)=\int_0^\infty u_\ve(y)\,y^{\rho-1}\,dy  \tag1.4$$   and
 $$(2\pi)^{-\rho}\,\Gamma(\rho)\,L_\ve(\rho-s+\h)=\int_0^\infty\bigl(f(iy)-(-1)^\ve f(-iy)\bigr)\,y^{\rho-1}\,dy\,, \tag1.5$$
where $$u_0(y)=\frac1{\sqrt y}\,u(iy)\,, \qquad 
  u_1(y)=\frac{\sqrt y}{2\pi i}\,u_x(iy)\qquad\bigl(u_x=\frac{\partial u}{\partial x}\bigr)\,.\tag1.6$$ 
The correspondence $f\mapsto\psi$ is given by formula  {\rm (0.5),} for some fixed $c(s)\ne0${\rm .}  \endproclaim

\demo{{R}emark} In the last two lines of (1.3), we wrote the estimates on $\psi$ which we will actually obtain for the period
functions $\psi$ attached to Maass wave forms (for which in fact $\s=1/2$), but for the proof of the implication
 (d)$\,\Rightarrow\,$(a)  we will in fact only need to assume a weaker estimate in the right half-plane; namely 
$$\psi(z)=\cases\text O\bigl(|z|^{-\s+\delta}\bigr)&(\Re(z)\ge0,\;|z|\le1),\\ 
\text O\bigl(|z|^{-\s-\delta}\bigr)&(\Re(z)\ge0,\;|z|\ge1)\endcases\tag1.7$$
for some $\delta>0$. Combining the two implications, we see that a Maass period function $\psi$ satisfying (1.7) for any 
positive $\delta$ automatically satisfies it with $\delta=\s$. This is another instance of the ``bootstrapping" phenomenon mentioned
in the introduction and studied in Chapter III. In fact, we will see in Chapter~III that such a $\psi$ actually has asymptotic expansions
in the right half-plane (or even in any wedge $|\arg z|<\pi-\ve$) of the form
$$\aligned \psi(z)&\sim C_0^*\,+\,C_1^*z\,+\,C_2^*z^2\,+\,\cdots\qquad\qquad\qquad\text{as $|z|\to0$}\,,\\
 \psi(z)&\sim \,-C_0^*z^{-2s}+C_1^*z^{-2s-1}-C_2^*z^{-2s-2}+\cdots\qquad\text{as $|z|\to\infty$}\,.\endaligned\tag1.8$$  

Notice, too, that the ``{(a)$\,\Rightarrow\,$(d)}" direction of Theorem 1 gives one part of the Theorem stated
in the introduction, since if we start with a Maass wave form then we get a holomorphic function $\psi$ satisfying (1.2) and (1.3),
and its restriction to $\R_+$ is a real-analytic function satisfying (0.2) and (0.3) (or even a strengthening of (0.3),
namely $\psi(x)=\text O(1)$ as $x\to0$, $\psi(x)=\text O(1/x)$ as $x\to\infty$).

\section{Periodicity, $L$-series, and the three-term functional equation} 

The modular group $\G$ is generated by the two transformations $z\mapsto z+1$ and $z\mapsto-1/z$.  The content of Theorem~1
can then be broken down correspondingly into two main parts.  The first part will be treated now: we will show that
the periodicity of $u$ is equivalent to the periodicity of $f$ (here and in future, ``periodic" means 
``1-periodic," i.e.\ invariant under $z\mapsto z+1$), to the property that $L_0$ and $L_1$ are (ordinary) Dirichlet
series, and to the fact that $\psi(z)$ satisfies the three-term functional equation.  The second (and harder) part, which will be 
done in Section~4, says that under suitable growth assumptions the following three conditions are equivalent: the
invariance of $u$ under $z\mapsto-1/z$, the functional equations of $L_0$ and $L_1$, and the continuability of $\psi$ 
from $\C\sm\R$ to $\Cc$.

\nonumproclaim{Proposition 1} Let $s$ be a complex number{\rm ,} $\s=\Re(s)${\rm .}  Then equations  {\rm (1.4)--(1.6)}
  give one\/{\rm -}\/to\/{\rm -}\/one
correspondences  between the following three classes of functions\/{\rm :}
\medbreak
\item{\rm (a)} a periodic solution $u$ of $\D u=s(1-s)u$ in $\H$ satisfying the
growth condition $u(x+iy)={\rm O}(y^A)$ for some $A<\min\{\s,1-\s\}${\rm ;}\smallbreak
\item{\rm (b)} a pair of Dirichlet $L$-series $L_\ve(\rho)\quad(\ve=0,\,1)${\rm ,} convergent in some half\/{\rm -}\/plane{\rm ;}
\smallbreak
\item{\rm (c)} a periodic holomorphic function $f(z)$ on $\C\sm\R$ satisfying 
 $$f(z)={\rm O}\bigl(|\Im(z)|^{-A}\bigr)$$ for some $A>0${\rm .}

 \endproclaim   
 
\demo{Proof} As is well-known, the equation $\D u=s(1-s)u$ together with the periodicity of $u$ and the growth estimate
given in {(a)} are jointly equivalent to the representability of $u$ by a Fourier series of the form
$$ u(z)=\sqrt y\,\sum_{n\ne0}A_n\,K_{s-\h}(2\pi|n|y)\,e^{2\pi inx}\qquad(z=x+iy,\quad y>0)\tag 1.9$$
with coefficients $A_n\in\C$ of polynomial growth. (We need the growth condition to eliminate exponentially large terms
$\sqrt y\,I_{s-\h}(2\pi|n|y)\,e^{2\pi inx}$ and ``constant" terms $\alpha y^s+\beta y^{1-s}$ in the Fourier expansion of $u$.)
 To such a $u$ we associate the two Dirichlet series $L_0$ and $L_1$ defined by
$$L_\ve(\rho)=\sum_{n=1}^\infty\frac{A_{n,\ve}}{n^\rho}\qquad(\ve=\text{0 or 1},\quad A_{n,\ve}=A_n+(-1)^\ve A_{-n})\tag1.10$$
and the periodic holomorphic function $f$ in $\C\sm\R$ defined by 
$$ f(z)=\cases \phantom-\sum_{n>0} n^{s-\h}\,A_n\,e^{2\pi inz}\qquad\;\,\text{if $\,\Im(z)>0$},\\
  -\sum_{n<0} |n|^{s-\h}\,A_n\,e^{2\pi inz}\qquad\text{if $\,\Im(z)<0$}.\endcases\tag1.11$$
(The minus sign in front of the second sum will be important later.)
The polynomial growth of the $A_n$ implies that $L_0$ and $L_1$ converge in a half-plane and that $f(x+iy)$ is bounded by
a power of $|y|$ as $|y|\to0$.  Conversely, if we start either with two Dirichlet series $L_0$ and $L_1$ which are convergent 
somewhere, or with a periodic and holomorphic function $f(z)$ in $\C\sm\R$ which is $\,\text O\bigl(|\Im(z)|^{-A}\bigr)\,$ 
for some $A>0$, then the expansion (1.10) or (1.11) defines coefficients $\{A_n\}_{n\ne0}$ which are of polynomial growth in $n$ 
(evidently so in the former case, and by the standard Hecke argument in the latter case).  Then if we define $u$ 
by (1.9), we find that the functions $u_\ve$ defined by (1.6) have the Bessel expansions
$$u_\ve(y)=\sum_{n=1}^\infty\,(ny)^\ve A_{n,\ve}\,K_{s-\h}(2\pi ny) \qquad(y>0,\quad\ve=0\text{ or 1})\,;$$ 
this in conjunction with the standard Mellin transform formulas
$$\int_0^\infty e^{-2\pi y}\,y^{\rho-1}\,dy=(2\pi)^{-\rho}\,\Gamma(\rho)\,,\qquad
 \int_0^\infty K_{s-\h}(2\pi y)\,y^{\rho-1}\,dy = \gamma_s(\rho) $$
shows that the functions $L_\ve$, $f$ and $u$ are indeed related to each other by formulas (1.4) and (1.5).  \enddemo 
 
We have now described the passage from a periodic solution of $\D u=s(1-s)u$ to a holomorphic periodic function $f$.
The passage from $f$ to $\psi$ is given by the following purely algebraic result.

\nonumproclaim{Proposition 2} Let $\psi(z)$ be a function in the complex upper half\/{\rm -}\/plane and $s$ a
 complex number{\rm ,}
$s\notin\Z${\rm .} Then $\psi(z)$ satisfies the functional equation $(1.2)$ 
if and only if the function $\psi(z)+z^{-2s}\psi(-1/z)$ is
periodic{\rm .}

More precisely{\rm ,}  formulas $(0.5)$ and $(0.6)${\rm ,} for any two constants $c(s)$ and $c^\star(s)$ 
satisfying $c(s)\,c^\star(s)=1-e^{-2\pi is}${\rm ,}
define a one\/{\rm -}\/to\/{\rm }\/-one correspondence between solutions $\psi$ of $(1.2)$
 and periodic functions $f$ in $\H${\rm .} The same holds true in the
lower half\/{\rm -}\/plane{\rm ,} but with the condition on $c(s)$ and $c^\star(s)$ now being
 $c(s)\,c^\star(s)=1-e^{2\pi is}${\rm .}\endproclaim

\demo{Proof} If $\psi$ satisfies (1.2), then
$$\align 0&=\biggl[\psi(z+1)-\psi(z)+\bigl(\frac1{z+1}\bigr)^{2s}\,\psi\bigl(\frac z{z+1}\bigr)\biggr]\\
  &\qquad - (z+1)^{-2s}\biggl[\psi\bigl(\frac z{z+1}\bigr)-\psi\bigl(\frac{-1}{z+1}\bigr)
  +\bigl(\frac{z+1}z\bigr)^{2s}\,\psi\bigl(\frac{-1}z\bigr)\biggr]\\
&=\biggl[\psi(z+1)+(z+1)^{-2s}\psi\bigl(\frac{-1}{z+1}\bigr)\biggr]
  -\biggl[\psi(z)+z^{-2s}\psi\bigl(\frac{-1}z\bigr)\biggr]\,,\endalign$$
so $\psi(z)+z^{-2s}\psi(-1/z)$, and hence also the function $f(z)$ defined by (0.6), is periodic. Conversely, if $f$ is periodic 
and we substitute for $\psi$ from (0.5), then we find that the difference of the left- and right-hand sides of (1.2) is a linear combination 
of three expressions of the form $f(w+1)-f(w)$ and hence vanishes. It is easy to check that (0.5) and (0.6) are inverse to each other
if the product of $c(s)$ and $c^\star(s)$ has the value given in the proposition. \enddemo

\demo{Remarks}  1.
The choice of the normalizing constants $c(s)$ is not important, but to have a well-defined correspondence we must make it explicit. We choose   
$$c(s)= \frac{i\,\pi^{-s}}{\G(1-s)}\,,\qquad c^\star(s)=\pm\,\frac{2\pi^{s+1}}{\G(s)}\,e^{\mp i\pi s}\,.\tag1.12$$ 
In the second formula, in which $c^\star(s)$ should more properly be denoted $c^\star_\pm(s)$,
the upper sign is to be chosen in the upper half-plane and the lower one in the lower half-plane. We have chosen to 
split the product $1-e^{\mp2\pi is}$ into two (reciprocals of) gamma functions because when we discuss the degeneration of
our story at integral values of $s$ (which we will do in Chapter IV) it will be convenient to have $c(s)$ nonzero for negative
integral $s$ and $c^\star(s)$ nonzero for positive integral $s$. 

\medbreak 2.   We would have liked to add a part ``{(d)}" to Proposition 1 giving a fourth  class of functions $\psi(z)$
equivalent to the other three.  Unfortunately, in the nonMaass case the growth conditions cannot be made to match\break up exactly 
as they do in Theorem~1.  The condition $f(z)\ll|\Im(z)|^{-A}$\break implies that $z^{-2s}_{\phantom{|}}f(-1/z)$, and hence $\psi(z)$, is
bounded by a multiple of\break $|\Im(z)|^{-A}(1+|z|^{2A-2\s})$ in $\C\sm\R$, but in the converse direction, imposing this growth condition
on $\psi(z)$ only permits us to deduce that the function $f(z)$ defined by (0.6) satisfies the same estimate.  Since $f$ is periodic
this is equivalent to saying that $f(z)$ is $\text O(|\Im(z)|^{-A})$ as $|\Im(z)|\to0$ and has at most polynomial growth 
as $|\Im(z)|\to\infty$, and since $f$ is also holomorphic and hence has an expansion in powers of $e^{2\pi iz}$, this implies that $f(z)$
in each half-plane is the sum of a constant $a_\pm$ and an exponentially small remainder term $\,\text O(e^{-2\pi|\Im(z)|})\,$ 
as $\Im(z)\to\pm\,\infty$. Because of these two constants this class of functions $f(z)$ is not exactly the one occurring in Proposition 1(c), 
but contains it as a subspace of codimension 2. (These extra two constants correspond to the possible terms $y^s$ and $y^{1-s}$ in
a solution of $\D u=s(1-s)\,u$ in the upper half-plane.) This point is related to the existence of Eisenstein series and the corresponding
modifications to the theory when noncuspidal wave forms are allowed, as discussed in Section~1 of Chapter~IV.

\section{Even and odd}  
As  mentioned in the introduction, a Maass wave form $u(z)$
is called {\it even} or {\it odd} if $u(-\bar z)=\pm u(z)$.  Since $u(-\bar z)$ for any Maass wave form $u$
is another wave form, it is clear that we can decompose the space $\M_s$ of Maass wave forms with eigenvalue
$s(1-s)$ into the direct sum of the spaces $\M_s^\pm$ of even and odd forms. (In all known cases, and conjecturally
for all $s$, $\dim(\M_s)=0$~or~1, so one of $\M_s^\pm$ is always 0.)  If we restrict to even
or odd forms, then the description of the correspondences $u\lr\{A_n\}\lr L\lr f$ becomes somewhat simpler: we need
only the coefficients $A_n$ for $n\ge1$ (since $A_{-n}=\pm A_n$), we have only one $L$-series 
$L(\rho)=\sum_{n=1}^\infty A_n n^{-\rho}$ (since either $u_0(y)$ or $u_1(y)$ is identically zero), and the function $f(z)$ 
need only be specified in the upper half-plane (since $f$ in the lower half-plane is then determined by $f(-z)=\mp f(z)$).  

On the period side we have a similar decomposition.  Let $\FE$ denote the space of solutions of the
three-term functional equation (1.2) in $\C\sm\R$ (holomorphic or continuous, with or without growth conditions, 
and defined in $\C\sm\R$, $\H$, $\H^-$, $\Cc$ or $\R_+$), and denote by $\FE^\pm$ the space of functions
of the same type satisfying the functional equation
  $$\psi(z)=\psi(z+1)\pm z^{-2s}\,\psi\bigl(\frac {z+1}z\bigr)\,.\tag 1.13$$   
\nonumproclaim{Proposition} $\;\FE=\FE^+\oplus\FE^-\,.$\endproclaim

\demo{Proof} If $\psi(z)$ satisfies (1.2), then one checks directly that the function $$\psi^\tau(z)\;:=\;z^{-2s}\psi(1/z)\tag1.14$$
also does. The involution $\tau$ therefore splits $\FE$ into a $(+1)$- and a $(-1)$-eigen\-space.  We claim that these are just
the spaces $\FE^+$ and $\FE^-$.  Indeed, if $\psi$ satisfies (1.13) then $\psi$ is $(\pm1)$-invariant
under $\tau$ (since the right-hand side of (1.13) is) and then the last term of (1.13) can be replaced by
$+z^{-2s}\psi^\tau(z^{-1}+1)=(z+1)^{-2s}\psi(z/(z+1))$, so $\psi\in\FE$; and conversely, if $\psi\in\FE$
is  $(\pm1)$-invariant under $\tau$ then we can replace the last term in (1.2) by $\pm (z+1)^{-2s}\psi^\tau(z/(z+1))$
to get (1.13). \enddemo

There are then even and odd versions of Theorems 1, giving isomorphisms between $\M_s^\pm$ and
the subspaces of $\FE^\pm$ consisting of solutions of (1.13) in $\Cc$ or $\R_+$ satisfying
appropriate growth conditions.  The even version reads:

\nonumproclaim{Theorem} Let $\{A_n\}_{n\ge1}$ be a sequence of complex numbers of polynomial growth{\rm .}
 Then the following
are equivalent\/{\rm :} \bigbreak
 \item{\rm (a)} the function $\sqrt y\sum\limits_{n=1}^\infty A_n\,K_{s-\h}(2\pi ny)\,\cos 2\pi nx\,$ 
is invariant under $z\mapsto\break
-1/z$ {\rm (}\/and 
hence is an even Maass wave form\/{\rm );}\/\medbreak
 \item{\rm (b)} the function $\gamma_s(\rho)\,\sum\limits_{n=1}^\infty A_nn^{-\rho}$ is entire of finite order 
and is invariant under $\rho\mapsto1-\rho${\rm ;}\medbreak
 \item{\rm (c)} the function defined by $$\,\pm\sum\limits_{n=1}^\infty n^{s-1/2}A_n(e^{\pm2\pi inz}-z^{-2s}e^{\mp2\pi in/z})\,$$ for
$\Im(z)\!\bs\!0$  extends holomorphically to $\Cc$ and is bounded in the right half\/{\rm -}\/plane.

\endproclaim

The odd version is similar except that we must replace ``cos" by ``sin" (and ``even" by ``odd") in part {(a)}, replace
$\gamma_s(\rho)$ by $\gamma_s(\rho+1)$ and ``invariant" by ``anti-invariant" in part {(b)}, and omit the $\pm$ sign before
the summation in part {(c)}.   The direct proof of either the odd or even version is slightly simpler than the proof of
Theorem 1 because there is only one nonzero function $u_\ve$ and only one Dirichlet series to deal with; but on the other hand
there are two cases to be considered  rather than one, so that we have preferred to give a uniform treatment.

\section{Relations between Mellin transforms; proof of Theorem 1}

In Section 2 we saw how the invariance of $u$ under $z\mapsto z+1$ corresponds to the existence of the two Dirichlet
series $L_0$ and $L_1$ and to the three-term functional equation of $\psi(z)$ in $\C\sm\R$, and also how the 
invariance under $z\mapsto-1/z$ translates into the functional equations of $L_0$ and $L_1$.  In this section we 
give the essential part of the proof of Theorem 1 by showing how the functional equations of the $L$-series $L_\ve$ 
both implies and follows from the extendability of $\psi$ to all of $\Cc$ (assuming appropriate growth conditions).

The main tool will be Mellin transforms and their inverse transforms, which are integrals along vertical lines,
so we will often need growth estimates on such lines.  We introduce the terminology ``$\a$-exponential decay"
to denote a function which grows at most like $\,\text O\bigl(|\rho|^A\,e^{-\a|\rho|}\bigr)\,$ for some $A$ as
$|\rho|\to\infty$ along a vertical line or in a vertical strip.  For instance, the gamma function is of
$(\pi/2)$-exponential decay in every vertical strip.

\demo{The implications {\rm (a)}$\,\Leftrightarrow\,${\rm (b)}} This part of Theorem 1, which is due to Maass (see [12]), 
follows easily from the discussion in Section~2.  Indeed, from $\D u=s(1-s)u$ and the fact that $\D$ commutes with the action
of
$\G$,  it follows that the function $u(z)-u(-1/z)$ satisfies the same  differential equation and therefore vanishes identically 
if it vanishes to second order on the positive imaginary axis, i.e\., if the two functions $u_0$ and $u_1$ satisfy
$u_\ve(1/y)=(-1)^\ve y\,u_\ve(y)$ ($\ve=0,\,1$), which by virtue of (1.4) translates immediately into the functional 
equation (1.1). The only point which has to be made is that for the converse
direction, which depends on writing $u_\ve(y)$ as an inverse Mellin transform
$$ u_\ve(y)=\frac1{2\pi i}\int_{\Re(\rho)=C}L^*_\ve(\rho)\,y^{-\rho}\,d\rho\qquad(C\gg0)\,,$$
we need an estimate on the growth of $L^*_\ve$ in vertical strips in order to ensure the convergence of the
integral.  Such an estimate is provided by the \PL theorem: the Dirichlet series $L_\ve(\rho)$
is absolutely convergent and hence uniformly bounded in some right half-plane, so by the functional equation
it also grows at most polynomially on vertical lines $\Re(\rho)=C$ with $C\ll0$; and since by assumption $L_\ve(\rho)$ 
is entire of finite order, the \PL theorem then implies that it grows at most polynomially in $|\rho|$ on any vertical 
line. It then follows from the definition of the functions $L^*_\ve(\rho)$ that they are of $(\pi/2)$-exponential
decay, since the gamma factor is.  This growth estimate ensures the rapid convergence of the inverse Mellin transform 
integral above and will be needed several times below.
\enddemo

\demo{The implication {\rm (b)}$\,\Rightarrow\,${\rm (c)}}  
Let $f$ be the periodic holomorphic function associated to $L_0$ and $L_1$ by (1.10) and   (1.11). Since the $A_n$ have polynomial growth,
it is clear that $f(z)$ is exponentially small as $|\Im(z)|\to\infty$ and of at most polynomial growth as  $|\Im(z)|\to0$.
We have to show that the function $f(z)-z^{-2s}f(-1/z)$ continues analytically from $\,\C\sm\R\,$ to $\,\Cc\,$ and satisfies the
given growth estimate in the right-half plane. (Here the specific normalization in (1.11), with the 
extra minus sign in the lower half-plane, will be essential.)  To do this, we proceed as follows. Denote by
  $$\wf_\pm(\rho)=\int_0^\infty f(\pm iy)\,y^{\rho-1}\,dy\qquad(\Re(\rho)\gg0)\tag1.15$$
the Mellin transform of the restriction of $f$ to the positive or negative imaginary axis. (The integral converges 
by the growth estimates just given.) Then for $\Re(\rho)$ large enough we have
$$\align
\noalign{\vskip8pt}
&\tag1.16a\\
\noalign{\vskip-26pt} \hskip.6in \wf_\pm(\rho)&= \int_0^\infty\biggl(\pm\sum_{n=1}^\infty n^{s-\h}\,A_{\pm n}\,e^{-2\pi
ny}\biggr)\,y^{\rho-1}\,dy
\\
&=\pm\frac{\G(\rho)}{2(2\pi)^\rho}\,\biggl(L_0(\rho-s+\h)\pm L_1(\rho-s+\h)\biggr) \\
&=\pm\frac1{\pi^{s+1}}\,\G\bigl(\frac{\rho+1}2\bigr)\,\G\bigl(\frac{2s-\rho+1}2\bigr)\,
\cos\pi\bigl(s-\frac\rho2\bigr)\,L^*_0(\rho-s+\h) \tag 1.16b\\
&\quad+\frac1{\pi^s}\,\G\bigl(\frac\rho2\bigr)\,\G\bigl(\frac{2s-\rho}2\bigr)\,
\sin\pi\bigl(s-\frac\rho2\bigr)\,L^*_1(\rho-s+\h)\,,\endalign$$
where to get the last equality we have used the standard identities
$$\G(x)=2^{x-1}\,\pi^{-1/2}\,\G\bigl(\frac x2\bigr)\G\bigl(\frac{x+1}2\bigr)\,,\quad\G(x)\G(1-x)=\frac\pi{\sin\pi x}\;.$$
From (1.16a) it follows that $\wf_\pm(\rho)$ is holomorphic except for simple poles at $\rho=0,\,-1,\,-2$,\dots and
is of $(\pi/2)$-exponential decay on any vertical line.  Therefore the inverse Mellin transform
$$f(\pm iy)=\frac1{2\pi i}\int_{\Re(\rho)=C}\wf_\pm(\rho)\,y^{-\rho}\,d\rho\qquad(y>0)$$
converges for any $C>0$ and extends analytically to all $y$ with $|\arg(y)|<\pi/2$, giving the integral representation
$$f(z)=\frac1{2\pi i}\int_{\Re(\rho)=C}\wf_\pm(\rho)\,e^{\pm i\pi\rho/2}\,z^{-\rho}\,d\rho\quad(z\in\C,\;\Im(z)\bs0)$$
and therefore 
$$\align f(z)-z^{-2s}\,f(&-1/z) =\frac1{2\pi i}\int_{\Re(\rho)=C}
\wf_\pm(\rho)\,\bigl[e^{\pm i\pi\rho/2}\,z^{-\rho}-e^{\mp i\pi\rho/2}\,z^{-2s+\rho}\bigr]\,d\rho\\
   &=\frac1{2\pi i}\int_{\Re(\rho)=C} \bigl[e^{\pm i\pi\rho/2}\,\wf_\pm(\rho)
   -e^{\mp i\pi(2s-\rho)/2}\,\wf_\pm(2s-\rho)\bigr]\,z^{-\rho}\,d\rho \endalign$$
for $0<C<2\Re(s)$ and $\Im(z)\bs0$.
But formula (1.16b) together with the functional equations of $L^*_0$ and $L^*_1$ and the elementary trigonometry identities
$$ \align \pm &e^{\pm i\pi\rho/2}\,\cos\pi(s-\rho/2)\mp e^{\mp i\pi(2s-\rho)/2}\,\cos\frac{\pi\rho}2=\,i\,\sin\pi s \\
  &e^{\pm i\pi\rho/2}\,\sin\pi(s-\rho/2)+ e^{\mp i\pi(2s-\rho)/2}\,\sin\frac{\pi\rho}2\,=\;\sin\pi s\endalign $$ give
\smallbreak\noindent (1.17)
$$\align&\frac{\pi^{s+1}}{i\,\sin\pi s}\,\bigl[\,e^{\pm i\pi\rho/2}\,\wf_\pm(\rho)\,-\,e^{\mp i\pi(2s-\rho)/2}\,\wf_\pm(2s-\rho)\,\bigr]
 \\  &\;\,=\,\G\bigl(\frac{\rho+1}2\bigr)\G\bigl(\frac{2s-\rho+1}2\bigr)\,L^*_0(\rho-s+\h)
   -i\pi\,\G\bigl(\frac{\rho}2\bigr)\G\bigl(\frac{2s-\rho}2\bigr)\,L^*_1(\rho-s+\h)\,,\endalign$$
so miraculously the {\it two} integral representations of $f(z)-z^{-2s}f(-1/z)$ in the upper and in the lower half-plane coalesce 
into a {\it single} integral representation
$$\align
& \hskip-12pt  f(z)-z^{-2s}f(-1/z)\\ \noalign{\vskip4pt}
&\qquad =\frac{\sin\pi s}{2\pi^{s+2}}\;\int_{\Re(\rho)=C} 
\biggl[\G\bigl(\frac{\rho+1}2\bigr)\G\bigl(\frac{2s-\rho+1}2\bigr)\,L^*_0(\rho-s+\h)\\
&\hskip1.6in
-i\pi\,\G\bigl(\frac{\rho}2\bigr)\G\bigl(\frac{2s-\rho}2\bigr)\,L^*_1(\rho-s+\h)\biggr]\,z^{-\rho}\,d\rho\,,\endalign$$ and
this now converges for all $z$ with $|\arg(z)|<\pi$ (i.e\. for all $z\in\Cc$) because the expression in square brackets is of
$\pi$-exponential decay.    

Finally, the estimates for $f(z)-z^{-2s}f(-1/z)$ in the right half-plane follow easily from the
integral representation just given.  Indeed, the integral immediately gives a (uniform)
bound $\,\text O\bigl(|z|^{-C}\bigr)\,$ in this half-plane for any $C$ between 0 and $2\s$, but since the integrand
is meromorphic with simple poles at $\rho=0$ and $\rho=2s$ we can even move the path of integration to a vertical line 
$\Re(\rho)=C$ with $C$ slightly to the left of~0 or to the right of $2\s$, picking up a residue proportional to 1 or
to $z^{-2s}$, respectively.  This gives an even stronger asymptotic estimate than the one  in {(c)}, and by moving the path of 
integration even further we could get the complete asymptotic expansions (1.8), with the coefficients $C_n^*$
being multiples of the values of the $L$-series $L_0(\rho)$ and $L_1(\rho)$ at shifted integer arguments. We omit the
details since we will also obtain these asymptotic expansions in Chapter III by a completely different method. 
\enddemo

\demo{The implication {\rm (c)}$\,\Rightarrow\,${\rm (d)}}  This is essentially just the algebraic identity proved in Section~2
(Proposition 2), since the function $\psi$ defined by (0.5) automatically satisfies the three-term functional equation (1.2).
The first estimate in (1.3) is trivial since it is satisfied separately by $f(z)$ and $z^{-2s}f(-1/z)$, and the other two estimates
were given explicitly in {(c)} as conditions on the function $f$.
\enddemo

\demo{The implication {\rm (d)}$\,\Rightarrow\,${\rm (b)}}  Now suppose that we have a function $\psi(z)$ satisfying the
conditions in part {(d)} of Theorem 1.  We already saw in Section2  (Proposition~2 and Remark~2) that the functional
equation (1.2)  and the estimates (1.3) imply the existence of Dirichlet series $L_\ve$ and of a periodic holomorphic function
$f(z)$ related to each other and to $\psi$ by (1.10), (1.11), (0.5) and (0.6).  The last two growth conditions in (1.3) imply that
the Mellin transform integral 
$$\widetilde\psi(\rho)=\int_0^\infty\psi(x)\,x^{\rho-1}\,dx\tag1.18$$ 
converges for any $\rho$ with $0<\Re(\rho)<2\s$ and moreover that we can rotate the path of integration from the
positive real axis to the positive or negative imaginary axis to get
$$i^\rho\,\int_0^\infty\psi(iy)\,y^{\rho-1}\,dy=\widetilde\psi(\rho)
   =i^{-\rho}\,\int_0^\infty\psi(-iy)\,y^{\rho-1}\,dy\,.\tag 1.19$$
The same estimates also imply that $f(iy)=\text O(1)$ as $|y|\to0$, and of course $f(iy)$ is exponentially small for 
$|y|\to\infty$ by (1.11), so the Mellin transform integral (1.15) converges for all $\rho$ with $\Re(\rho)>0$ and we have  
$$ \align c(s)\,\int_0^\infty\psi(\pm iy)\,y^{\rho-1}\,dy
  & =\int_0^\infty\bigl[f(\pm iy)-e^{\mp i\pi s}\,y^{-2s}\,f(\pm i/y)\bigr]\,y^{\rho-1}\,dy \tag1.20\\
  & = \wf_\pm(\rho)-e^{\mp i\pi s}\,\wf_\pm(2s-\rho)\endalign$$
for $\rho$ in the strip $0<\Re(\rho)<2\s$.  Therefore (1.19) gives
$$ e^{i\pi\rho/2}\,\bigl[\wf_+(\rho)-e^{-i\pi s}\,\wf_+(2s-\rho)\bigr]\,=
\,e^{-i\pi\rho/2}\,\bigl[\wf_-(\rho)-e^{i\pi s}\,\wf_-(2s-\rho)\bigr]\,.$$
Substituting for $\wf_\pm(\rho)$ in terms of $L^*_\ve(\rho)$ from (1.16b), which is valid in the strip  $0<\Re(\rho)<2\s$ 
by the same argument as before, and moving the appropriate terms on each side of the equation to the other side, we obtain
$$\align &\frac{i\pi^{1-s}}{\G\bigl(\frac{1-\rho}2\bigr)\,\G\bigl(\frac{1+\rho-2s}2\bigr)}
    \,\bigl[L^*_0(\rho-s+\h)-L^*_0(s-\rho+\h)\bigr]\\
  &\hskip.5in =\,\frac{\pi^{-s}}{\G\bigl(\frac{2-\rho}2\bigr)\,\G\bigl(\frac{2+\rho-2s}2\bigr)}
    \,\bigl[L^*_1(\rho-s+\h)+L^*_1(s-\rho+\h)\bigr]\,.\endalign $$
But the left-hand side of this equation changes sign and the right-hand side is invariant under 
$\rho\mapsto 2s-\rho$, so both sides must vanish. This gives the desired functional equations of $L_0$ and $L_1\,$.

We still have to check that $L^*_\ve$ is entire and of finite order.  We already know that  $\wf_\pm(\rho)$ is holomorphic
for $\Re(\rho)>0$, so formula (1.5) implies that $L_\ve(\rho)$ is also holomorphic in this half-plane.  If $0<\s<1$, then by 
looking at the poles of the gamma-factor $\gamma_s(\rho+\ve)$ we deduce that $L^*_\ve(\rho)$ has no poles in the smaller 
right half-plane $\Re(\rho)>|\h-\s|$.  The functional equation then implies that $L^*_\ve(\rho)$ also has no poles in the
left half-plane $\Re(\rho)<1-|\h-\s|$, and since these two half-planes intersect, $L^*_\ve(\rho)$ is
in fact an entire function of~$\rho$.  Furthermore, any of the integral representations which show that $L_\ve^*(\rho)$ is 
holomorphic also shows that it is of at most exponential growth in any vertical strip, which together with the
functional equation and the boundedness of $L_\ve$ in a right half-plane implies that $L_\ve^*$ is of finite order. 

This completes the proof of the theorem if $0<\s<1$. To extend the result to all $s$ with $\s>0$, we use the fact (which will be 
proved in \S3 of  Chapter III) that any function $\psi$ satisfying the assumptions of {(d)} has
asymptotic representations of the form (1.8) near 0 and $\infty$.  These asymptotic expansions give the locations of the
poles of the Mellin transform $\widetilde\psi(\rho)$ of $\psi$; namely, it has simple poles at $\rho=-m$ and $\rho=2s+m$ with
residues $C_m^*$ and $(-1)^{m+1}C_m^*$, respectively.  But equations (1.19) and (1.20) together with equation (1.16b),
the functional equation (1.1) (now established), and the trigonometric identities preceding equation~(1.17) combine to give 
$$ \align C(s)\,\widetilde\psi(\rho)\;&=\;\G\bigl(\frac{\rho+1}2\bigr)\,\G\bigl(\frac{2s-\rho+1}2\bigr)\,L_0^*(\rho-s+\h)\\
&\qquad-i\,\pi\;\G\bigl(\frac\rho2\bigr)\,\G\bigl(\frac{2s-\rho}2\bigr)\,L_1^*(\rho-s+\h) \tag 1.21 \endalign$$
for some nonzero constant $C(s)$. Replacing $\rho$ by $2s-\rho$ just changes the sign of the second term on the 
right, so the first and second terms on the right are proportional to $\widetilde\psi(\rho)+\widetilde\psi(2s-\rho)$ and
$\widetilde\psi(\rho)-\widetilde\psi(2s-\rho)$, respectively, and this implies that both $L_0^*$ and $L_1^*$ are entire,
since the poles of $\widetilde\psi(\rho)\pm\widetilde\psi(2s-\rho)$ are cancelled by those of the gamma factors.

Finally, we observe that if we had used the weaker condition (1.7) instead of the last two conditions of (1.3) 
then the same proof would have gone through, except that we would have had to work with $\rho$ in the smaller strip 
$\s-\delta<\Re(\rho)<\s+\delta$ instead of $0<\Re(\rho)<2\s$.

\vglue24pt \centerline{\bf Chapter  II. The period correspondence via integral transforms}
\vglue16pt

Let $u(z)$ be a Maass wave form with spectral parameter $s$. In Chapter I we defined the associated period function $\psi$ 
in the upper and lower half-planes by the formula
  $$\psi(z) \e \pm\sum_{n=1}^\infty n^{s-1/2}\,A_{\pm n}\,\bigl(e^{\pm 2\pi inz} - z^{-2s}e^{\mp 2\pi in/z}\bigr)
   \qquad(\Im(z)\bs0)\,.\tag2.1$$
(Here and throughout this chapter, the symbol $\,\e\,$ denotes equality up to a factor depending only on $s$.)
On the other hand, as was mentioned in the introduction to the paper, the original definition of the period function
as given (in the even case) in [10] was by an integral transform; namely
 $$\psi_1(z) \e \int_0^\infty zt^s\,\bigl(z^2+t^2\bigr)^{-s-1}\,u(it)\,dt\qquad\bigl(\Re(z)>0\bigr)\,,\tag2.2$$
where we have written ``$\psi_1$" instead of ``$\psi$" to avoid ambiguity until we have proved the proportionality of the
two functions. This definition is more direct, but does not make apparent either of the two main properties of the period function,
viz., that it extends holomorphically to the cut plane $\Cc$ and that it satisfies the three-term functional equation (0.1). 
In this chapter we will study the function defined by the integral (2.2) from several different points of view, each of
which leads to a proof of these two properties and each of which brings out different aspects of the theory.  More specifically:

In Section~1 we extend the definition (2.2) to include odd as well as even wave forms and show, using the representation of $\psi(z)$
as an inverse Mellin transform of the $L$-series of $u$ which was the central result of Chapter I, that the functions $\psi$ and
$\psi_1$ are proportional in their common region of definition.  This shows that $\psi_1$ extends to the left half-plane 
and satisfies the functional equation (since these properties are obvious from the representation (2.1)), and at the same time 
that $\psi$ extends holomorphically across the positive real axis (since this property is clear from (2.2)).

In Section~2 we give a direct proof of the two desired properties. It turns out that the integrand in (2.2) can be written in a 
canonical way as the restriction to the imaginary axis of a closed 1-form defined in the whole upper half-plane. This permits 
us to deform the path of integration, and from this the analytic extendability, the functional equation of $\psi_1(z)$, and
the proportionality of $\psi_1$ and $\psi$ follow in a very natural way.

In Section~3 we study the properties of the function $\psi_1$ when $u(z)$ is assumed to be an eigenfunction of the Laplace operator but
no longer to be $\G$-invariant.  Specifically, we show that the invariance of $u(z)$ under the transformation $$S:\,z\mapsto-1/z$$ 
is equivalent to the identity $\psi_1(1/z)=z^{2s}\psi_1(z)$, while the invariance of $u(z)$ under the transformation $$T:\,z\mapsto z+1$$
is reflected in the fact that the function $\psi_2$ defined by
\vglue-6pt
  $$ z^{-2s}\,\psi_2\bigl(1+\frac1z\bigr) = \psi_1(z)-\psi_1(z+1)\qquad(\Re(z)>0)\tag 2.3$$
\vglue6pt\noindent 
extends holomorphically to $\Cc$ and satisfies the same identity $\psi_2(1/z)=\break z^{2s}\psi_2(z)$. It is then easy to deduce
that if $u$ is invariant under both $S$ and $T$ then $\psi_1$ equals $\psi_2$ (i.e., $\psi_1$ satisfies the three-term functional
equation) and is proportional to $\psi$. We also show how to interpret these relationships in terms of a summation formula due to Ferrar.

 In Section~4 we explain how to pass from the Maass form $u(z)$ to the function $\psi(z)$ via an intermediate function $g(w)$ which is
related to both of them by integral transforms. This is the approach used in [10], but the derivation given here is simpler.  The 
function $g$ links $u$ and $\psi$ in a very pretty way: it is an entire function whose values at integral multiples 
of $2\pi i$ are the Fourier coefficients of $u$ and whose Taylor coefficients at $w=0$ are proportional to the Taylor coefficients 
of $\psi$ at $z=1$.  We also show that the Taylor coefficients of $\psi$ at $z=0$ are proportional to special values of the $L$-series 
of $u$, in exact analogy with the corresponding fact for the coefficients of the period
polynomials of holomorphic modular forms which is discussed in Chapter~IV.

Finally, in Section~5 we give an expression for $\psi(z)$ as a formal integral transform of an automorphic distribution on $\R$
which is obtained from the function $u(z)$ by a limiting process as $z$ approaches the boundary.  This representation makes the
properties of $\psi$ intuitively clear and ties together several of the other approaches used in the earlier sections of the
chapter. 

\advance\sectioncount by -4
 
\section{The integral representation of $\psi$ in terms of $u$}  
In this section we will use the $L$-series proof given in Section~4 of Chapter~I to prove that the Mellin transforms
of the restrictions of $\psi_1$ and $\psi$ to $\R_+$ are proportional and hence that $\psi_1$ is a multiple of $\psi$.
This helps to understand the properties of the period function, since each representation puts different aspects into
evidence: formula (2.1) and the elementary algebraic lemma (Proposition~2) of~Section~2 of Chapter~I make it clear that
$\psi$ satisfies the  three-term functional equation (1.2), but not at all obvious that it extends holomorphically from
$\C\sm\R$ to
$\Cc$,  while from (2.2) (or its odd analogue) it is obvious that $\psi_1$ extends across the positive real axis but not that
it extends beyond the imaginary axis or that it satisfies the three-term equation.

We will state the result in a uniform version which includes both even and odd Maass forms. The definition (2.2)
of the function $\psi_1(z)$ must then be replaced by
$$\psi_1(z) = 2sz \int_0^\infty \frac{t^{s+1/2}\,u_0(t)}{(z^2+t^2)^{s+1}}\,dt
  -2\pi i\int_0^\infty \frac{t^{s-1/2}\,u_1(t)}{(z^2+t^2)^s}\,dt\qquad\bigl(\Re(z)>0\bigr)\,, \tag 2.4$$
where $u_0$ and $u_1$, as in Chapter I (eq.~(1.6)), denote the renormalized value and normal derivative of $u$ 
restricted to the imaginary axis. Of course only the first term is present is $u$ is even and only the second if $u$ is odd.
\nonumproclaim{Proposition} For $u$ a Maass wave form the function $\psi_1(z)$ defined by $(2.4)$ is proportional to the
period function $\psi(z)$ defined in Theorem {\rm 1} of Chapter {\rm I.}\endproclaim 

\demo{Proof} We will prove this by comparing the Mellin transforms of $\psi_1$ and~$\psi$.  Since $u_0$ and $u_1$ are of
rapid decay at both 0 and infinity, the function defined by (2.4) is holomorphic in the right half-plane and its restriction
to the positive real axis is $\,{\rm O}(1)$ near 0 and $\,{\rm O}(z^{-2\Re(s)})$ at infinity.  Hence the Mellin transform 
$\int_0^\infty\psi_1(x)\,x^{\rho-1}\,dx$ exists in the strip $0<\Re(\rho)<2\Re(s)$.  We can compute it easily by
interchanging the order of integration and recognizing the inner integral as a beta integral. The result of the computation is
$$ \align \int_0^\infty \psi_1(x)\,x^{\rho-1}\,dx
 &\doteq \G\bigl(\frac{\rho+1}2\bigr)\,\G\bigl(\frac{2s-\rho+1}2\bigr)\,\int_0^\infty u_0(t)\,t^{\rho-s-1/2}\,dt \\
 &\qquad-i\pi\,\G\bigl(\frac\rho2\bigr)\,\G\bigl(\frac{2s-\rho}2\bigr)\,\int_0^\infty u_1(t)\,t^{\rho-s-1/2}\,dt\,,\endalign $$
which agrees (up to a constant factor) with the expression on the right-hand side of (1.21).  Applying the inverse Mellin 
transform, we deduce that the functions $\psi$ and of $\psi_1$ are proportional when restricted to the positive real axis and
hence, by analytic continuation, also in the right half-plane. \enddemo

\section{The period function as the integral of a closed $1$-form} 
In the last section we proved that $\psi_1(z)$ satisfies the conditions of a Maass period function 
(extendability to $\Cc$, functional equation, and growth) by using the results of Chapter I. Now we give a direct proof
by interpreting (2.4) as the integral of a closed 1-form along a path and then deforming the path.

We start with a more general construction. If $u$ and $v$ are two differentiable functions 
of a complex variable $z=x+iy$, let $\{u,v\}=\{u,v\}(z)$ be the differential 1-form (Green's form) defined by the formula
 $$  \{u,v\} =\bigl(v\,\frac{\p u}{\p y}-u\,\frac{\p v}{\p y}\bigr)\,dx+\bigl(u\,\frac{\p v}{\p x}-v\,\frac{\p u}{\p x}\bigr)\,dy 
   = \left|\matrix u&\p u/\p x&\p u/\p y\\v&\p v/\p x&\p v/\p y\\0&dx&dy\endmatrix\right|\,. \tag2.5 $$
We also consider the complex version of this, defined by 
  $$ [u,v]=[u,v](z) =v\,\frac{\p u}{\p z}\,dz\,+\,u\,\frac{\p v}{\p\bz}\,d\bz\,,$$
where $dz=dx+i\,dy$, $\,d\bz=dx-i\,dy$, $\,\frac\p{\p z}=\h\frac\p{\p x}-\frac i2\frac{\p}{\p y}$, 
$\,\frac\p{\p\bz}=\h\frac\p{\p x}+\frac i2\frac\p{\p y}$. 

\nonumproclaim{Lemma} The forms $\{u,v\}$ and $[u,v]$ have the following properties\/{\rm :}\/ \bigbreak
\item{\rm (i)} $[u,v]+[v,u]=d(uv),\qquad[u,v]-[v,u]=-i\,\{u,v\}${\rm .} \medbreak 
\item{\rm (ii)} If $u$ and $v$ are eigenfunctions of the Laplacian with the same eigenvalue{\rm ,} then $\{u,v\}$ and $[u,v]$ are
closed forms{\rm .}\medbreak 
\item{\rm (iii)} If $z\mapsto g(z)$ is any holomorphic change of variables{\rm ,} then 
$\{u\circ g,v\circ g\} =\{u,v\}\circ g$ and $[u\circ g,v\circ g]=[u,v]\circ g${\rm .}

\endproclaim

\demo{Proof} (i) Direct calculation.

(ii) The statement holds for both the euclidean Laplacian $\D_0=\dfrac{\partial^2}{\partial x^2}+\dfrac{\partial^2}{\partial y^2}$ and the 
hyperbolic Laplacian $\D=-y^2\D_0$, since (for $\{u,v\}$, which suffices by~(i))
$$\align d\{u,v\}  &= \biggl[-\frac\partial{\partial y}\biggl(v\,\frac{\p u}{\p y}-u\,\frac{\p v}{\p y}\biggr)
    +\frac\partial{\partial x}\biggl(u\,\frac{\p v}{\p x}-v\,\frac{\p u}{\p x}\biggr)\biggr]\,dx\wedge dy \\
   &= \bigl(u\,\D_0v-v\,\D_0u\bigr)\,dx\wedge dy =  \bigl(v\,\D u-u\,\D v\bigr)\,\frac{dx\wedge dy}{y^2}\,.\endalign$$

(iii) Again it suffices by (i) to prove this for one of the two forms; this time $[u,v]$ is easier. Replacing $u$ and $v$ by $u\circ g$ 
and $v\circ g$ multiplies $\p u/\p z$ by $g'(z)$ and $\p v/\p\bz$ by $\overline{g'(z)}$, while replacing $z$ by $g(z)$ 
in $[u,v]$ leaves the coefficients $uv_z$ and $vu_{\bz}$ unchanged but multiplies $dz$ by $g'(z)$ and $d\bz$ by $\overline{g'(z)}$.\enddemo

We will apply this construction when $v$ is the $s^{\rm th}$ power of the function 
$$ R_\z(z)=\frac y{(x-\z)^2+y^2}=\frac i2\,\biggl(\frac1{z-\z}-\frac1{\bz-\z}\biggr)\qquad\quad(\z\in\C,\;z=x+iy\in\H)\,. $$
The main properties of this function are the transformation equation  
$$ R_{g\z}(gz)=(c\z+d)^2\,R_\z(z)\tag2.6$$
for $g=\bigr(\smallmatrix a&b\\c&d\endsmallmatrix\bigr)\in {\rm SL}(2,\R)$ and the differential equation 
$$\D\bigl(R_\z^s\bigr)=s(1-s)\,R_\z^s\qquad(s\in\C)\,.\tag 2.7$$
(Here the $s^{\rm th}$ power is well-defined if $\z\in\R$, since then $R_\z(z)>0$; in the general case we must restrict $z$ to the
complement in $\H$ of some path joining $\z$ and $\bar\z$ and choose the evident branch of $R_\z^s$.)  Both properties
can be proved easily either by direct calculation or (for $\z$ real, which suffices) by observing that $R_\z(z)=c^2\,\Im(gz)$ 
for any $g=\bigr(\smallmatrix a&b\\c&d\endsmallmatrix\bigr)\in {\rm SL}(2,\R)$ with $g^{-1}(\infty)=\z$.

From (2.7) and (ii) of the lemma it follows that if $u:\H\to\C$ is an eigenfunction of $\D$ with
eigenvalue $s(1-s)$, then the differential form $ \{u,\, R_\z^s\}$ is closed.  Explicitly, we have
$$ \align\{u,\, R_\z^s\}(z)&=\biggl(\frac{sy^{s-1}(y^2-(x-\z)^2)}{\bigl((x-\z)^2+y^2\bigr)^{s+1}}\,u(z)
    \,+\,\frac{y^s}{\bigl((x-\z)^2+y^2\bigr)^s}\,u_y(z)\biggr)\,dx \\
  &\quad+\biggl(\frac{-2s(x-\z)y^s}{\bigl((x-\z)^2+y^2\bigr)^{s+1}}\,u(z)
    \,-\,\frac{y^s}{\bigl((x-\z)^2+y^2\bigr)^s}\,u_x(z)\biggr)\,dy\,. \endalign$$
The $dy$ part of this (with $x\mapsto0$, $y\mapsto t$, $\z\mapsto z$) is the integrand in (2.4), so (2.4) can be rewritten 
$$ \psi_1(\z) \,=\, \int_0^{i\infty}\{u,\,R_\z^s\}(z) \,,\tag 2.8$$
where $\Re(\z)>0$ and the integral is taken along the imaginary axis. But now, because the form $\{u,\, R_\z^s\}$ is closed, we can replace
the imaginary axis by any other path from 0 to $i\infty$ which passes to the left of both $\z$ and $\bar\z$. We can then move $\z$ to
any new point with $\z$ and $\bar\z$ to the right of the new path.  This gives the analytic continuation of $\psi_1$ to $\Cc$ for
any $u$ which is sufficiently small at 0 and~$\infty$ to ensure convergence, and in particular for $u$ a Maass form. The three-term
functional equation in the Maass case also follows easily. Indeed, (2.6) and part (iii) of the lemma imply 
$$  (c\z+d)^{-2s}\,\psi_1(g\z)=\int_{g^{-1}(0)}^{g^{-1}(\infty)}\{u,\,R_\z^s\}(z)$$
for any $g=\abcd\in {\rm SL}(2,\R)$ for which $u\circ g=u$, and hence for any $g\in\G$ in the Maass case; thus   we have
$$ \psi_1(\z)-\psi_1(\z+1)-(\z+1)^{-2s}\psi_1\bigl(\frac\z{\z+1}\bigr)
  =\biggl(\int_0^\infty-\int_{-1}^\infty-\int_0^{-1}\biggr) \{u,\,R_\z^s\}(z) =0\,.$$

A modification of the same idea can be used to give a direct proof of the proportionality of the functions $\psi_1$ and $\psi$ in 
the Maass case. We observe first that by part (i) of the lemma we can replace $\{u,R_\z^s\}$ by $[u,R_\z^s]$ in (2.8), since
the integral of $d(uR_\z^s)$ from 0 to $\infty$ vanishes. The invariance properties of $[\cdot,\cdot]$ and $R_\z^s$ then give
$$ \psi_1(\z)\e\int_0^{i\infty}[u,\,R_\z^s](z)= f_1(\z)-\z^{-2s}\,f_1(-1/\z)\tag2.9$$
for $\z\in\H$, where
$$ f_1(\z)\e\int_\z^{i\infty}[u,\,R_\z^s](z)\qquad\qquad(\z\in\H)\,.$$
(The point of replacing $\{u,R_\z^s\}$ by $[u,R_\z^s]$ in (2.8) is that the former has a singularity like $|z-\z|^{-s-1}$ as $z\to\z$
and hence cannot be integrated from $\z$ to $\infty$,  whereas the latter has only a $|z-\z|^{-s}$ singularity at $\z$, which is integrable
for $0<\Re(s)<1$.) The function $f_1$ is holomorphic, since
 $$ \frac{\p f_1(\z)}{\p\bar\z} = -u(z)\,\frac{\p}{\p\bz}\bigl(R_\z(z)^s\bigr)\biggr|_{z=\z}
= -\frac{is}2\,\frac{u(z)\,y^{s-1}}{(z-\z)^{s-1}(z-\bar\z)^{s+1}}\biggr|_{z=\z}=0\,,$$
and is obviously periodic (replace $z$ by $z+1$ in the integral),  so (2.9) gives a second proof of the three-term functional
equation by virtue of Proposition~2 of Chapter~I, Section~2.  Moreover, term-by-term integration using formula 10.2 (13),
  [6, Vol.~II, p.~129]
shows that $f_1$ is proportional to the function $f$ defined by (1.11) in the upper half-plane.  A similar argument works in the
lower  half-plane, but now with $f_1(\z)$ defined as $-\int_{\bar\z}^{i\infty}[R_\z^s,u](z)$. This gives a new proof that the
function $\psi$  defined in $\C\sm\R$ by (2.1) extends to $\Cc$ and that the functions $\psi_1$ and $\psi$ are proportional. 
We have also obtained  an explicit representation of the holomorphic function $f:\H\to\C$ associated to a Maass wave form as
an Abel transform:
\nonumproclaim{Proposition} Let $u$ be a Maass wave form with Fourier expansion $(1.9)${\rm .}
 Then the function $f$ defined by $(1.11)$ has 
the integral representation
$$ f(\z)\e\int_\z^{i\infty} \,\biggl(\frac{\p u(z)}{\p z}\,\frac{y^s}{(z-\z)^s(\bz-\z)^s}\,dz 
  \,+\,\frac{is}2\,u(z)\,\,\frac{y^{s-1}}{(z-\z)^{s-1}(\bz-\z)^{s+1}}\,d\bz\biggr) $$
in the upper half\/{\rm -}\/plane{\rm ,} where the integral can be taken along any path{\rm .}\endproclaim

\section{The incomplete gamma function expansion of $\psi$}

We now consider the behavior of the correspondence $u\mapsto\psi_1$ defined by the integral (2.2) or (2.4) for 
functions $u$ in the upper half-plane which are not necessarily invariant under all of $\G$
but only under $S$ or $T$ separately.  This will lead to a number of alternate descriptions of the
transformation and to a better understanding of its properties, as well as supplying yet another proof
that $\psi_1$ satisfies the three-term functional equation when $u$ is in fact a Maass wave form.
For simplicity we restrict our attention to the even case $u(z)=u(-\bz)$, so that $\psi_1$ is defined by (2.2).
The odd case will be treated briefly at the end of the section.

\demo{$S$\/{\rm -}\/invariance} This is very easy. Substituting $t\to1/t$ in (2.2) gives
  $$ u(-1/z)= u(z)\quad\Rightarrow\quad \psi_1(1/z)=z^{2s}\psi_1(z)\,;\tag2.10$$
i.e., the $S$-invariance of $u$ is reflected in the $\tau$-invariance of $\psi_1$, where $\tau$ is the involution defined in (1.14). 
Conversely, if we know that $\psi_1=\psi_1^\tau$, then we can deduce that $u$ is $S$-invariant if it is assumed to be an eigenfunction
of~$\D$. (Without this assumption it follows only that the restriction of $u$ to the positive imaginary axis is $S$-invariant.)
\enddemo

\demo{$T$\/{\rm -}\/invariance} 
Now suppose that $u$ is a $T$-invariant even eigenfunction of $\D$ and is bounded in the upper half-plane i.e., that $u$ 
satisfies all the properties of an even Maass form {\it except} the invariance under $S$.  Then $u$ has a Fourier expansion
(1.9) with $A_n= {\rm O}\bigl(n^{1/2}\bigr)$ by the usual Hecke argument, and by term-by-term integration we find that 
the function $\psi_1$ defined by (2.2) has the expansion
  $$ \psi_1(z) \e \sum_{n=1}^\infty\,n^{s-1/2}\,A_n\,\Cs(2\pi nz)\,, \tag 2.11 $$
where  
  $$ \Cs(z)\;\e\; \int_0^\infty \frac{z\,t^{s+1/2}}{(z^2+t^2)^{s+1}}\,K_{s-1/2}(t)\,dt\qquad\bigl(\Re(z)>0\bigr)\,. \tag 2.12 $$
The function $\Cs(z)$ (a special case of the ``Lommel function") has very nice properties.  The ones we will use are given in the 
following proposition and corollary.
\enddemo

\nonumproclaim{Proposition 1} Let $s\in\C$, $\Re(s)>0${\rm ,} and $\Re(z)>0${\rm .}  Then
   $$\Cs(z)\,\e\,e^{i\pi(s-1/2)+iz}\G(1-2s,iz)+e^{-i\pi(s-1/2)-iz}\G(1-2s,-iz)\tag 2.13$$
where $\G(a,x)=\int_x^\infty e^{-t}\,t^{a-1}\,dt\,$ is the incomplete gamma function.{\rm }  The function $\Cs(z)$
is also given{\rm ,} up to factors depending only on $s${\rm ,} by each of the following formulas\/{\rm :}
$$\align  \Cs(z)\;&\e\; \int_0^\infty \frac{\cos t}{(z+t)^{2s}}\,dt\,, \tag 2.14 \\
\Cs(z)\;&\e\; \int_0^\infty \frac{w^{2s}}{w^2+1}\,e^{-wz}\,dw\,,\tag 2.15 \\
\Cs(z\;)&\e\; \sin(\pi s+z)-\sum_{n=0}^\infty\frac{(-1)^n\,z^{2n+1-2s}}{\G(2n+2-2s)}\,. \tag 2.16\endalign$$\endproclaim

\demo{{R}emark} The  normalization of $\Cs(z)$ is not particularly important for our purposes, so we did not include it in
the statement.  Occasionally we will want to have fixed a choice. We then take the right-hand side of equation (2.15) as the definition 
of $\Cs(z)$, which determines the implied constants of proportionality in equations (2.12), 
(2.13), (2.14) and (2.16) as $\pi^{-1/2}2^{s+1/2}\G(1+s)$, $\frac12\G(2s)$, $\G(2s)$, and $\pi/\sin2\pi s$, respectively.
\enddemo

\demo{Proof} The equality (up to constants) of the various functions in (2.12)--(2.15) can be verified by standard manipulations
or by looking them up in tables of integrals. A more enlightening proof is to observe that each of these functions has 
polynomial growth at infinity and satisfies the differential equation
 $$  \Cs''(z)+\Cs(z)\,\e\,z^{-2s-1}$$
(where the implied constant, with the normalization just given, is $\G(2s+1)$), and that these properties characterize $\Cs$ uniquely. 
(To get the differential equation for (2.12) and (2.14), integrate by parts and use the second-order differential equations satisfied 
by $\,\cos t\,$ and $K_{s-1/2}(t)$.)  The power series expansion (2.16) of $\Cs(z)$ at 0 is obtained from (2.13) using 
   $e^x\bigl(\G(a)-\G(a,x)\bigr) =\G(a)\,\sum_{n=0}^\infty(-1)^n x^{a+n}/\G(a+n+1)$. \enddemo

\nonumproclaim{{C}orollary} The function $\Cs(z)$ extends holomorphically to the cut plane $\Cc$
 and has the asymptotic expansion{\rm ,}
uniform in any wedge $|\arg z|\le\pi-\ve${\rm ,}  
  $$ \Cs(z)\;\sim\;\sum_{n=0}^\infty\frac{(-1)^n\,\G(2s+2n+1)}{z^{2s+2n+1}}\tag2.17 $$
as $|z|\to\infty$. Moreover{\rm ,} it satisfies
  $$  e^{\pm\pi is}\,\Cs(z) + e^{\mp\pi is}\,\Cs(-z) \e e^{\pm iz}\qquad(z\in\C,\;\Im(z)\bs0)\,.\tag2.18$$\endproclaim

\demo{Proof} The analytic continuation and the symmetry property (2.18) follow easily from either (2.14) or (2.16).  The asymptotic expansion
can be obtained from any of the formulas (2.12)--(2.15) if $\Re(z)>0$, and from formula (2.14) for all $z\in\Cc\,$. \enddemo

\nonumproclaim{Proposition 2} Let $u$ be a bounded{\rm ,} even{\rm ,} and $T$\/{\rm -}\/invariant eigenfunction
 of $\D$ in the upper half\/{\rm -}\/plane{\rm ,} and
define $\psi_1$ and $\psi_2$ by equations $(2.2)$ and $(2.3)${\rm .} Then\/{\rm :}
\smallbreak
{\rm i)} Both $\psi_1(z)$  and $\psi_2(z)$ extend holomorphically to $\Cc${\rm ,} and $\psi_2$ is $\tau$\/{\rm -}\/in\-variant{\rm .}
\smallbreak 

{\rm ii)} The function in $\C\smallsetminus\R$ defined by $\psi_1(z)+e^{\mp2\pi is}\psi_1(-z)$ for 
$\Im(z)\bs0$ is periodic{\rm .} \endproclaim
 \demo{Proof} Each term of the series (2.11) extends holomorphically to $\Cc$ by the corollary, and the series  converges
absolutely and locally uniformly because of the estimates $A_n= {\rm O}(\sqrt n)$ and $\Cs(z)= {\rm O}\bigl(|z|^{-2\Re(s)-1}\bigr)$. 
This gives the extension of $\psi_1$ to $\Cc$. Now using (2.14) we get
$$ \psi_1(z)-\psi_1(z+1) \e \sum_{n=1}^\infty n^{1/2-s}\,A_n\,\int_0^1\,\frac{\cos2\pi nt}{(z+t)^{2s}}\,dt$$
for $z\in\Cc$, the calculation being justified by the absolute convergence. Replacing $z$ by $1/(z-1)$ gives
  $$\psi_2(z) \e \sum_{n=1}^\infty n^{1/2-s}\,A_n\,\int_0^1\frac{\cos2\pi nt}{(1-t+tz)^{2s}}\,dt \tag2.19$$
for $z\in\C\setminus(-\infty,1]$. The right-hand side obviously defines a holomorphic function in $\Cc$, and the symmetry 
under $z\mapsto1/z$ follows by substituting $t\mapsto1-t$. This proves (i), and part (ii) is easily seen to 
be an equivalent statement to part (i).\enddemo

 Equations (2.11) and (2.18) give the following more explicit version of part (ii): 
  $$ \psi_1(z)+e^{\mp2\pi is}\psi_1(-z)\e c^\star(s)\,f(z) \qquad(\Im(z)\bs0)\,,\tag2.20$$
where $f(z)$ is the periodic function in $\C\sm\R$ defined by equation (1.11) and $c^\star(s)$ the factor (which depends on 
the sign of $\,\Im(z)\,$ and hence cannot be absorbed by the $\e$) defined in (1.12).

\demo{The Maass case} Combining the results of the last two subsections immediately yields a proof of both
the three-term functional equation of $\psi_1$ and the proportionality of $\psi_1$ and $\psi$ when $u$ is a Maass form. Indeed, 
$\psi_1$ is then\break $\tau$-invariant by (2.10), so $e^{\mp2\pi is}\psi_1(-z)=z^{-2s}\psi_1(-1/z)$. Part (ii) of Proposition 2 then says that
$\psi_1(z)+z^{-2s}\psi_1(-1/z)$ is periodic, which is equivalent to the  three-term functional equation by Proposition 2 of Chapter I,
Section~2. The same argument applied to (2.20) implies the proportionality of $\psi_1(z)$ and the function $\psi(z)$ defined by
(2.1). (Compare equations (0.6) and (0.5).)

The results of our analysis can be summarized in the following proposition.
\enddemo

\nonumproclaim{Proposition 3} Let $A_n$ $(n\ge1)$ be complex numbers satifying $A_n={\rm O}(n^{1/2})$ and $s\,$ 
a complex number with $\Re(s)>0${\rm .} 
Then the following conditions are equivalent\/{\rm :}\/
\smallbreak
{\rm i)} The numbers $\,A_n\,$ are the Fourier coefficients of an even Maass form with eigenvalue $s(1-s)${\rm ;}
\smallbreak
{\rm ii)} For $z\in\C$ with $\Im(z)\bs0${\rm ,}
  $$\quad\sum_{n=1}^\infty A_n\,n^{s-1/2}\,\Cs(2\pi nz) =\pm c_s\sum_{n=1}^\infty\,A_n\,n^{s-1/2}\,\bigl(e^{\pm2\pi inz}
     -z^{-2s}e^{\mp 2\pi in/z}\bigr)\tag2.21$$  
  with $c_s=\dfrac\pi{2i}\,\csc(\pi s)$ and $\Cs(z)$ as in Proposition $1$ and the following remark\/{\rm ;}
\smallbreak
{\rm iii)} The function  on the left\/{\rm -}\/hand side of $(2.21)$ is $\tau$\/{\rm -}\/invariant{\rm .}  \endproclaim    

The equivalence of condition (iii) with the apparently stronger condition (ii) says that  once we assume that $\psi(z)$ has
an expansion in terms of Lommel functions, then merely requiring its $\tau$-invariance implies the full period property. Thus
the three-term functional equation in this situation contains a lot of ``overkill."  We also remark that in the Maass case 
we have $\psi_1=\psi_2$ (this is just the three-term functional equation), so that (2.19) gives yet another
representation, distinct from (2.1) and (2.11), of the period function of an (even) Maass wave form in terms of its Fourier coefficients.

\demo{Alternate approach\/{\rm :} Ferrar summation}  The relationship stated in Proposition 3 can be seen in another way by
making use  of the Ferrar summation formula~[7]. This summation formula characterizes sequences $\{A_n\}$ whose associated
Dirichlet series 
$L(\rho)=\sum A_nn^{-\rho}$ satisfy a functional equation $L(1-\rho)=\Phi(\rho)L(\rho)$ by the property that 
$$\sum_{n=1}^\infty A_n h(n)=\sum_{n=1}^\infty A_nh^\F(n) \quad\; (\,+\text{ possible residue terms}\,)\tag2.22$$  
for arbitrary (nice) test functions $h$, where  $h\to h^\F$ is the integral transform whose kernel function is the inverse Mellin
transform of $\Phi(\rho)$. For example, the integral transform corresponding to the functional equation satisfied by the Riemann zeta function
is the Fourier cosine transform, and the Ferrar summation formula (2.22) for the sequence $A_n=1$ is just the Poisson summation formula.

When $\Phi(\rho)$ is the function corresponding to the functional equation of an even Maass form with spectral parameter $s$, then we 
find from standard tables of Mellin transforms that the integral transform specified by Ferrar's formula is   
$$h^\F(y)=\int_0^\infty F_s(yx)\,h(x)\,dx $$   with $F_s(\xi)=F_{1-s}(\xi)$ given by 
$$F_s(\xi) = \frac\pi{\cos\pi s}\,\bigl(J_{2s-1}(4\pi\sqrt\xi)-J_{1-2s}(4\pi \sqrt\xi)\bigr) + 4\,\sin\pi s\,K_{2s-1}(4\pi \sqrt\xi)$$  
(and no contribution from the residue terms in this case).  Applying formula (2.22) to the test function 
$h(x)=K_{s-1/2}(ax)$ gives the characteristic inversion formula
  $$\sum_{n=1}^\infty A_n K_{s-1/2}(2\pi ny)=y^{-1}\,\sum_{n=1}^\infty A_n K_{s-1/2}(2\pi n/y)$$
of an (even) Maass form, while applying it to $h(x)=x^{s-1/2}e^{iax}$ and to $h(x)=x^{s-1/2}\Cs(ax)$
gives (after some manipulations) the identities (2.21) and
  $$\sum_{n=1}^\infty A_n n^{s-1/2}\Cs(2\pi nz)=z^{-2s}\,\sum_{n=1}^\infty A_n n^{s-1/2}\,\Cs(2\pi n/z)\,,$$
respectively, which by Proposition 3 above also characterize the Fourier coefficients of Maass wave forms.
\enddemo

\demo{The odd case}  Finally, we should say briefly what happens in the case of functions $u$ which are anti-invariant
under $z\mapsto-\bz$.  In this case $\psi_1$ is defined by the second term in (2.4). The analogue of (2.10) is again easy:
if $u$ is $S$-invariant then $u_1(1/t)=-tu_1(t)$ and $\psi_1$ is anti-invariant under $\tau$.  If $u$ is $T$-invariant, then
computing $u_1(t)$ from the Fourier expansion of $u$ and integrating term by term in (2.4) gives an expansion like (2.11),
but with the cosine-like function $\Cs(z)$ replaced by its ``sine" analogue
  $$ \Cal S_s(z)\;\e\; \int_0^\infty \frac{t^{s+1/2}}{(z^2+t^2)^s}\,K_{s-1/2}(t)\,dt\qquad\bigl(\Re(z)>0\bigr)\,.  $$
This function (suitably normalized) is related to $\Cs$ by
$$  \Cal S_s=\Cal C_{s-1/2}\,,\qquad \Cs=-\Cal S^0_{s-1/2}\,,\qquad \Cal S_s'=-\Cs\,,\qquad \Cs'=\Cal S_s^0\,,$$
where $\Cal S_s^0(z)=\Cal S_s(z)-\G(2s)z^{-2s}= {\rm O}\bigl(z^{-2s-2}\bigr)$. The story now works much as before, the 
$T$-invariance of $u$ being reflected finally in the anti-$\tau$-invariance of $\psi_2$ or equivalently in the periodicity of 
$\psi_1(z)-e^{\mp2\pi is}\psi_1(-z)$.  The only new point is that the series $\sum n^{s-1/2}A_n\,\Cal S_s(2\pi nz)$ does not converge 
absolutely and must be replaced~by $\sum n^{s-1/2}A_n\,\Cal S_s^0(2\pi nz)+\G(2s)L_1(s+\h)(2\pi z)^{-2s}$, where $L_1(\rho)$ is
(the analytic continuation of) the $L$-series associated to $u$.  We omit the details.
\enddemo

\section{Other integral transforms and intermediate functions} 

We again consider a periodic, bounded, and even solution of $\D u=\break s(1-s)u$, where $\Re(s)>0$. In the last section we developed 
the properties of the associated holomorphic functions $\psi_1$ and $\psi_2$ defined by (2.2) and (2.3), and in particular showed that $\psi_2$
continues to $\Cc$ and is $\tau$-invariant.  The former property implies that the power series expansion of $\psi_2(1+z)$ has radius 
of convergence 1, so the function $g$ defined by 
$$g(w)=\sum_{m=0}^\infty\frac{C_m}{\G(m+2s)}\,w^m\,,\qquad \text{where}\qquad\psi_2(1+z)=\sum_{m=0}^\infty C_m\,z^m\,,\tag2.23$$
is entire and of exponential type (i.e. $g(w)\ll e^{(1+\epsilon)|w|}\,$), while the $\tau$-invariance of $\psi_2$ 
is equivalent to the functional equation
  $$ g(-w)=e^w\,g(w) \tag 2.24$$
by a calculation whose proof we leave to the reader as an exercise.  ({\it Hint}: Compute the function $\psi_2$ corresponding 
via (2.23) to $g(w)=w^{2n}e^{-w/2}$.) 

\nonumproclaim{Proposition}  The Fourier coefficients of $u$ are related to the function $g$ by 
  $$  g(\pm2\pi in)\,\e\,n^{-s+1/2}\,A_n\qquad\qquad(n=1,\,2,\,\ldots).\tag 2.25$$  \endproclaim

In other words, the entire function $g$ simultaneously interpolates the\break Fourier coefficients of $u$ and is related via its Taylor
expansion  to the function $\psi_2$. Since $\psi_2$ is proportional to $\psi$ in the case when $u$ is a Maass wave form, 
the proposition gives us an explicit way to recover a Maass form from its associated period function.

\demo{Proof} Using the integral formula (2.19) and expanding $\,(1+tz)^{-2s}$ by the binomial theorem, 
we find that the coefficients $C_m$ defined by (2.23) are given by
 $$C_m \,\e\, (-1)^m\binom{m+2s-1}m\,\sum_{n=1}^\infty n^{1/2-s}\,A_n\,\int_0^1 t^m\,\cos2\pi nt\,dt\, ,$$  
so
 $$\align g(w) &\e \sum_{n=1}^\infty n^{1/2-s}\,A_n\,\int_0^1\,e^{-wt}\,\cos2\pi nt\,dt\tag 2.26\\
 &=\bigl(1-e^{-w}\bigr)\,\sum_{n=1}^\infty n^{1/2-s}\,A_n\,\int_0^\infty\,e^{-wt}\,\cos2\pi nt\,dt\\
 &=\bigl(1-e^{-w}\bigr)\,\sum_{n=1}^\infty n^{1/2-s}\,A_n\,\frac w{w^2+4\pi^2n^2}\,.\endalign$$
The proposition follows immediately.  \enddemo

 The identity (2.26) was proved in [10, Lemma 5.1] in a different way, which we sketch briefly.  Let $u$ be a periodic 
eigenfunction (not necessarily Maass) with Fourier expansion (1.9), say with $A_n= {\rm O}(n^{1/2})$.
Define a function $\phi$ on $\{w\,:\,|\Im(w)|<2\pi\}\,$ as the Hankel transform 
  $$\phi(w)=w^{1-s}\,\int_0^\infty \sqrt{wt}\,J_{s-\h}(wt)\,u(it)\,dt\,.\tag2.27$$
 Integrating term-by-term, using formula 8.13 (2),  [6, Vol.~II, p.~63] gives
  $$\phi(w)=\sum_{n=1}^\infty n^{\h-s}\,A_n\,\frac w{w^2+(2\pi n)^2}\,.\tag$2.28$ $$
Hence $\phi(w)$ continues meromorphically to an odd function in the whole\break complex plane with simple poles of residue $(2\pi
n)^{-s+1/2}A_n/2$ at $w=\pm2\pi in$ $(n=1,\,2,\ldots)$ and no other poles, and is of polynomial growth away from the imaginary axis. 
In other words, the translation invariance of $u$ is reflected in $\phi(w)$ in the properties of being odd and 
having (simple) poles only in $2\pi i\Z$.  Now {\it define} the function $g(w)$ by
   $$g(w)=\bigl(1-e^{-w}\bigr)\,\phi(w)\,;\tag 2.29$$
then these properties translate into the properties that $g$ is entire of exponential type and satisfies the functional
equation (2.24) and the interpolation property (2.25). On the other hand, $\psi_1$ can be expressed in terms of $\phi$ as a Laplace transform
$$\psi_1(z) \e \int_0^\infty e^{-zw}\,w^{2s-1}\,\phi(w)\,dw\qquad(\Re(z)>0)\,.\tag2.30$$
({\it Proof.} Substitute (2.27) into (2.30) and apply the integral 4.14(8) of [6, Vol.~I] to recover (2.2)
up to a constant depending only on $s$.) This and (2.29) give
$$\psi_1(z)-\psi_1(z+1)  \e \int_0^\infty e^{-zw}\,w^{2s-1}\,g(w)\,dw\qquad(\Re(z)>0)\, , $$ 
or, replacing $z$ by $1/z$ and then $w$ by $wz$,
$$\psi_2(1+z)  \e \int_0^\infty e^{-w}\,w^{2s-1}\,g(wz)\,dw\qquad(\Re(z)>0)\,.\tag2.31$$ 
The integral on the right converges for $\Re(z)>-1$ because $g$ is of exponential type, giving the holomorphic continuation of $\psi_2$ 
to the right half-plane, and differentiating (2.31) $n$ times and setting $z=0$, we recover the relationship (2.23) between the 
Taylor expansions of $g$ at $w=0$ and of $\psi_2$ at $z=1$. 

The two approaches can be related by observing that substituting (2.28) into (2.30) and using formula (2.15) for $\Cs(z)$ gives
 the Lommel function expansion (2.11).  

Finally, we observe that applying $\tau$ to equation (2.30) gives
$$\psi_1^\tau(z) \e \int_0^\infty e^{-w}\,w^{2s-1}\,\phi(wz)\,dw\qquad(\Re(z)>0)\,.$$
The same argument which was used to get from (2.31) to (2.23) now shows that the function $\psi_1^\tau(z)$ is $C^\infty$ at $z=0$ 
(more precisely, its derivatives from the right to any order exist) and has an asymptotic expansion given by
$$ \psi_1^\tau(z)\,\sim\,\sum\Sb m\ge0\\\text{$m$ odd}\endSb\G(m+2s)\,\frac{\phi^{(m)}(0)}{m!}\,z^m\qquad(z\to0,\quad\Re(z)>0)\,.\tag2.32$$
Note that the series has radius of convergence 0 because the radius of convergence of the power series of $\phi(w)$ at $w=0$ is finite.
On the other hand, expanding each term in (2.28) in a geometric series in $w^2$ we find that
$$ \frac{\phi^{(m)}(0)}{m!} = \frac{(-1)^{\frac{m-1}2}}{(2\pi)^{m+1}}\,L_0(m+s+\tfrac12)\qquad(\text{$m$ odd})\,,$$
where $L_0(\rho)=\sum A_n n^{-\rho}$ is the $L$-series of $u$.  In particular, in the Maass case, when $\psi_1^\tau=\psi_1=\psi_2\e\psi$, we obtain

\nonumproclaim{Proposition} Let $u$ be an even Maass wave form with eigenvalue 
$s(1-s)$, and $\psi(z)$ the associated period
function{\rm .}  Then
$\psi(z)$ is infinitely often differentiable from the right at $z=0$ with $\psi^{(m)}(0)=0$ for $m$ even and
  $$ \psi^{(m)}(0) \e \frac{m!}{(2\pi i)^m}\,\G(m+2s)\,L_0(m+s+\tfrac12)\qquad(\text{\rm $m$ odd})\,.\tag 2.33$$ 
 \endproclaim

This result, which can also be proved by using the invariance property $\psi=\psi^\tau$, the Lommel function expansion (2.11), and the asymptotic 
formula (2.17) for $\Cs(z)$, will be discussed again (in both the even and odd cases) in Chapter IV in connection with the holomorphic period theory.
In fact, (2.33) is the exact analogue of the fact that the period function associated to a holomorphic modular form is a polynomial 
whose coefficients are simple multiples of special values of its $L$-series.

\section{Boundary values of Maass wave forms}

A combination of two ideas led to the equivalence between Maass forms and their associated period functions developed in [10]: 
the Ferrar summation formula discussed in Section~3 above and the (formally) automorphic boundary form of the 
Maass form. This boundary form is a $\G$-invariant distribution on the boundary of the symmetric space $\H$ whose
existence follows from general results on boundary forms of eigenfunctions of invariant differential operators on symmetric spaces. 
However, the general theory is set up in a way which makes it hard to apply directly to our situation and which obscures the action of the translation and inversion generators $T$ and $S$ of $\G$.  
We will therefore present the ideas first from a naive point of view by showing how to express a Maass wave form $u$ and its associated
functions $f$, $\psi$, $g$, $\phi$ and $L_\ve$ formally as simple integral transforms of various types (Poisson, Stieltjes, Laplace and 
Mellin) of a single ``function" $U(t)$.  This makes transparent the relationships of these various functions to one another
and also lets one directly translate their special properties in the Maass case into a certain formal automorphy property of $U$.
The formal argument can then be made rigorous by interpreting $U$ as a $\G$-invariant object in a suitable space of distributions.

Suppose that we have a function or distribution $U(t)$ on the real line and associate to it three functions $u$, $f$ and $\psi$ as follows:
$$ \align u(z)\;&=\,y^s\int_{-\infty}^\infty |z-t|^{-2s}\,U(t)\,dt \qquad(z\in\H)\,,\tag2.34\\
f(z)\;&=\phantom{\,y^s}\int_{-\infty}^\infty (z-t)^{-2s}\,U(t)\,dt\qquad(z\in\C\sm\R)\,,\tag2.35\\
\psi(z)\;&=\phantom{\,y^s}\int_{-\infty}^0\!(z-t)^{-2s}\,U(t)\,dt\qquad(z\in\Cc)\,.\tag 2.36\endalign$$
Then (assuming that the integrals converge well enough) the function $u$ is an eigenfunction of the Laplace operator with eigenvalue $s(1-s)$
and the functions $f$  and $\psi$ are holomorphic in the domains given. If $U$ is also automorphic in the sense that
$$ U(t)\,=\, |ct+d|^{2s-2}\,U\bigl(\frac{at+b}{ct+d}\bigr) \qquad\text{for all $\abcd\in\G$}\,,\tag 2.37$$
then $u$ is $\G$-invariant, $f$ is periodic and is related to $\psi$ by (0.5) and (0.6), and $\psi$ satisfies the 
three-term functional equation (1.2) by virtue of Proposition 2, Section~2, of Chapter I or alternatively by the following formal
calculation:
$$\align 
&\hskip-.25in \psi(z)-\psi(z+1)\\ &\quad =\;\int_{-1}^0\frac{U(t)}{(z-t)^{2s}}\,dt
   \qquad(\text{by (2.37) with $\bigl(\smallmatrix a&b\\c&d\endsmallmatrix\bigr)=\bigl(\smallmatrix 1&1\\0&1\endsmallmatrix\bigr)
\bigr)$}\\ \noalign{\vskip4pt}
&\quad=\;\int_{-1}^0\frac{(1+t)^{2s-2}}{(z-t)^{2s}}\,U\bigl(\frac t{1+t}\bigr)\,dt
   \qquad(\text{by (2.37) with $\bigl(\smallmatrix a&b\\c&d\endsmallmatrix\bigr)=
\bigl(\smallmatrix 1&0\\1&1\endsmallmatrix\bigr)\bigr)$}\\ \noalign{\vskip4pt}
  &\quad=\int_{-\infty}^0\frac{U(t)\,dt}{(z-zt-t)^{2s}}\;=\;(z+1)^{-2s}\,\psi\bigl(\frac z{z+1}\bigr)\,.\endalign$$

 Now in fact there cannot be a reasonable function satisfying (2.37) (for instance, the value of $U(t)$ for $t$ rational would have to be 
proportional to $\,|\text{denom}(t)|^{2-2s}$), and the existence of even a distributional solution in the usual sense is not at all clear.
 We will come back to this issue a little later. First, we look at the properties of the integral transforms (2.34)--(2.36) for functions $U(t)$ 
for which the integrals do make sense and see how they are related; this will give new insight into the $u\lr\psi$ relationship 
which is the fundamental subject of this paper and will at the same time tell us what the object $U(t)$ should be when $u$ is a Maass wave form.

For convenience, we consider the even case when $U(t)=U(-t)$ and $u(-\bz)=u(z)$ (the odd case would be similar), but make no further
automorphy assumptions on~$U$.  We also change the name of the function defined by (2.36) to $\psi_1$, since it is related to the function (2.34)
by equation (2.2). This can be seen by the calculation
  $$\align
\noalign{\vskip12pt}
&\tag2.38\\
\noalign{\vskip-38pt}\qquad\quad  2\,&\int_0^\infty\frac{z\,\tau^{2s}\,d\tau}{(z^2+\tau^2)^{s+1}|t-i\tau|^{2s}}\,
     =\,\int_0^\infty \frac{(z+t)\,(\tau^2+zt)\,\tau^{2s}\,d\tau}{(\tau^2+z^2)^{s+1}(\tau^2+t^2)^{s+1}}\\
\noalign{\vskip4pt}
     &\;=\,\int_{-\infty}^\infty \frac{(z+t)\,dv}{\bigl(v^2+(z+t)^2\bigr)^{s+1}}\;\e\,(z+t)^{-2s}\qquad(t\in\R_+,\;\Re(z)>0)\,,\endalign$$
where the first equality is obtained (initially for $z\in\R_+$) by symmetrizing with respect to $\tau\mapsto zt/\tau$, the second by substituting 
$v=\tau-zt/\tau$, and the third by the homogeneity property of the integral.  (This also works without the assumption of evenness if we replace 
(2.2) by (2.4).)  A considerably easier calculation shows that the functions (2.35) and (2.36) are related by (2.20). 

Now consider the case when $U$ is periodic.  (As in Chapter I, we always mean by this ``1-periodic," i.e. $U(t+1)=U(t)$.)  Then $u$ and $f$ 
are also clearly periodic, while the $T$-invariance of $U$ is reflected in $\psi_1$ by the property $\psi_2^\tau=\psi_2$ with $\psi_2$
defined by (2.3). (This is the equivalence of parts (i) and (ii) of Proposition 2, Section~3, which was noted at the time.) But we
can be much more explicit, and tie the new approach in with the results of Section~3, by using the Fourier expansion of $U$.  If
we write this expansion (still in the even case) as
  $$ U(t)=\sum_{n=1}^\infty n^{\h-s}\,A_n\,\cos(2\pi nt)\,,\tag2.39$$
then standard integrals show that the periodic functions (2.34) and (2.35) have the Fourier expansions (1.9) and (1.11), respectively, while the 
representation (2.14) of the Lommel function shows that the function (2.36) is given by the expansion (2.11).  We can also connect with 
the results of Section~4 by defining
$$ \align g(w)&=\int_0^1 e^{-wt}\,U(t)\,dt\qquad(w\in\C)\,, \tag2.40\\ 
  \phi(w)&=\int_0^\infty e^{-wt}\,U(t)\,dt\qquad(\Re(w)>0)\tag2.41\endalign$$
and observing that these functions are related to $u$ and $\psi_1$ and each other by equations (2.23), (2.27) and (2.29)
and have the expansions given in (2.26) and (2.28), respectively. Finally, the $L$-series defined by (1.10) is expressed by  
   $$\frac12\,(2\pi)^{-\rho}\,\G(\rho)\,\cos\bigl(\frac{\pi\rho}2\bigr)\,L_0(\rho+s-\h)\,=\,\int_0^\infty U(t)\,t^{\rho-1}\,dt\,,\tag2.42$$
and again the relationship of this function to the others (namely, that it is proportional to the Mellin transforms of $u(iy)$,
$f(iy)$, $\psi(x)$, or $\phi(w)$) follows easily by comparing the various integral representations in terms of $U$. 

If the function $U$ is smooth as well as periodic, then the coefficients $A_n$ in (2.39) are of rapid decay and all the expansions just given
converge nicely. If instead we start with a sequence of coefficients $A_n$ of polynomial growth, then the series (2.39) no longer converges, 
but still defines a $T$-invariant distribution on the real line.  In particular, this is true when the $A_n$ are taken to be the Fourier 
coefficients of a Maass wave form, and in that case the estimate $A_n= {\rm O}(\sqrt n)$ is sufficient to make the various expansions converge, as 
discussed in the previous sections of this chapter. However, it is not immediately clear why the distribution $U(t)$ defined by (2.39) should
have the automorphy property (2.37) in the Maass case, or, for that matter, even what this automorphy property means. We now describe several 
different ways, both formal and rigorous, to see in what sense the series (2.39) can be considered to be an automorphic object when the
$A_n$ are the Fourier coefficients of a Maass form.

\medbreak 1.  The first approach is based on the asymptotic expansion near 0 of the $K$-Bessel functions occurring in the 
Fourier development of $u$; namely:
$$K_{s-1/2}(2\pi t)=\alpha_st^{1/2-s}\,\,+\alpha_{1-s}\,t^{s-1/2}\,+\,  {\rm O}(t^2)\qquad\qquad\text{as $t\to0$}\,.$$
(Here we are using $\Re(s)=\h$; if $\Re(s)$ had a different value the analysis would actually be easier because $K_{s-1/2}(t)$ would behave
asymptotically like a single power of $t$.)  Substituting this into (1.9) gives the formal asymptotic formula  
$$u(z)\sim \alpha_s\,U(x)\, y^{1-s}\,+\,\alpha_{1-s}\,\wU(x)\,y^s\qquad(z=x+iy,\quad y\to0)\,,\tag2.43$$
where $$U(t)=\sum_{n\ne0}|n|^{\h-s}\,A_n\,e^{2\pi int}\,,\qquad \wU(t)=\sum_{n\ne0}|n|^{s-\h}\,A_n\,e^{2\pi int}\,.\tag2.44$$
(The first of these expansions coincides with (2.39) in the even case.) Combining (2.43) with the $\G$-invariance of $u$
we obtain formally equation (2.37) and also the corresponding automorphy property of $\wU$ with $2-2s$ replaced by $2s$.
The rigorous version of this approach is the theory of boundary forms, discussed below. 

\medbreak
{2.} A second way to see formally why the Fourier series (2.39) should be automorphic when $\{A_n\}$ are the
coefficients of a Maass form is based on the properties of the associated period function. We know that the function $\psi$
defined on
$\C\sm\R$ by (2.1) extends analytically  to the positive real axis, so computing $\psi(x)$ for $x>0$ formally as the limit of
$\psi(z)$ from above and below we have the ``equality"
$$\align
& \sum_{n>0} n^{s-1/2}\,A_n\,\bigl(e^{2\pi inx} - x^{-2s}e^{-2\pi in/x}\bigr)\\
&\hskip.5in  =-\sum_{n<0}^\infty|n|^{s-1/2}\,A_n\,\bigl(e^{-2\pi inx} - x^{-2s}e^{2\pi in/x}\bigr),\endalign
$$
where of course none of the four series (taken individually) are convergent, though both sides of the equation are supposed to 
represent the same perfectly good function $\psi|\R_+$. Now moving two of the four terms to the other side of the equation gives
exactly the automorphy of $\wU(x)$ under $S$, and since the invariance under $T$ is obvious this ``proves" the automorphy in general.
Note the formal similarity between this argument and the ``criss-cross" argument used in Chapter I (equation (1.20) and
the following calculations) to prove the extendability to $\Cc$ of the function (2.1).

This approach, too, can be made rigorous, this time by using the theory of hyperfunctions, which is a alternative way to define
functionals on a space of test functions on $\R$ as the differences of integrals against holomorphic functions (here $f(z)$) in the 
lower and upper half-planes.  The theory in the Maass context is developed in [1], where the goal is to give a cohomological 
interpretation of theory of period functions. We refer the interested reader to [1] and also to Part II of the present paper, where
various related approaches will be discussed.

\medbreak
{3.}  A third approach is based on the $L$-series of $u$.  The functional equation (1.1) of the $L$-series says that (the
analytic continuation of) the left-hand side of (2.42) is invariant under $\rho\mapsto2s-\rho$, while the corresponding invariance of the 
right-hand side of (2.42) is formally equivalent to the automorphy under $\tau$ (and hence, since we are in the even case, under $S$) of $U$.

\medbreak {4.} The automorphy of $U$ under $S$ can also be obtained as a formal consequence of the Ferrar summation formula
discussed at the end of \S3. The kernel function $F_s(\xi)$ for the Ferrar transform, which was given there as a complicated
explicit linear combination of Bessel functions, has the simple integral representation
    $$F_s(\xi) = 4\,\xi^{-s+1/2}\,\int_0^\infty x^{2s-2}\,\cos(2\pi x)\,\cos(2 \pi \xi/x)\,dx\,.\tag2.45$$
Therefore Fourier inversion implies that the Ferrar transform of the function $h_t(x)=x^{-s+1/2}\,\cos(2\pi xt)$ is 
$h_t^\F(y)=y^{-s+1/2}\,t^{2s-2}\cos(2\pi y/t)$ for any $t>0$, and the Ferrar summation formula (2.22) applied to $h_t$ reduces
formally to the desired automorphy property $U(t)=|t|^{2s-2}U(1/t)$. To make sense of this latter identity we simply dualize by integrating
against an arbitrary (sufficiently nice) test function $\vp$. The formulas just given for $h_t$ and $h_t^\F$ show that the Ferrar transform of
$\widehat\vp$ (the Fourier cosine transform  of $\vp$) is $\widehat{\vp^\tau}$ (the Fourier cosine transform  of $\vp^\tau$), 
so the Ferrar summation formula becomes 
   $$\sum_{n=1}^\infty n^{\h-s}\,A_n \widehat\vp(n) =\sum_{n=1}^\infty n^{\h-s}\,A_n \widehat{\vp^\tau}(n)\,,\tag2.46$$ 
and this is precisely the desired automorphy property of $U$, if $U$ given by (2.39) is now thought of as a distribution. 

\medbreak {5.} We now make this distributional point of view rigorous by defining a precise space of test functions on which
$\G$ acts and a corresponding space of distributions to which $U$ belongs.  Let $\V$ be the space of $C^\infty$ functions $\vp$
on
$\R$ such that $\vp^\tau$ is also $C^\infty$; i.e.\ $\vp(t)$ has an asymptotic expansion $\vp(t)\sim
|t|^{-2s}\sum_{n\ge0}c_nt^{-n}$ as $|t|\to\infty$.  This space has an action of the group $G={\rm PSL}(2,\R)$ given by
$\bigl(\vp|g)(x):=|cx+d|^{-2s}\,\vp\bigl(\frac{ax+b}{cx+d}\bigr)$ for $g=\abcd\in G$.  This can be checked either directly or,
more naturally, by noting that $\V$ can be identified via $\vp\lr\Phi(x,y)=y^{-2s}\vp(x/y)$ with the space of $C^\infty$
functions $\Phi:\R^2\sm\{(0,0)\}\to\C$ with the homogeneity property $\Phi(tx,ty)=|t|^{-2s}\Phi(x,y)$, with the action of $G$
given simply by $\Phi\mapsto\Phi\circ g$. (A third model consists of $C^\infty$ functions $f$ on the circle, with
$\Phi(re^{i\theta})=r^{-2s}f(\theta)$ or $f(2\theta)=|\cos\theta|^{-2s}\vp(\tan\theta)$ and the corresponding $G$-action.) If
$\Re(s)=\h$, then $\V$ is nothing other than the space of smooth vectors in the standard Kunze-Stein model for the unitary
principal series.

We can now think of the ``function" $U$ defined by (2.39) as a linear map from $\V$ to $\C$, given by
    $$\vp\,\mapsto\,U[\vp] :=\sum_{n=1}^\infty n^{1/2-s}A_n\,\widehat\vp(n)\,.\tag 2.47$$
The series converges rapidly because $\widehat\vp(n)$ decays as $n\to\infty$ faster than any power of $n$, and hence gives an interpretation
of the formal integral $U[\vp]=\int_{-\infty}^\infty U(t)\,\vp(t)\,dt\,$. The meaning of the automorphy equation (2.37) is now simply the following:

\nonumproclaim{Proposition} 
Let $s\in\C$ with $\Re(s)>0${\rm ,} and $\{A_n\}_{n\ge1}$ a sequence of complex numbers of polynomial growth{\rm .}
Then the $A_n$ are the Fourier coefficients of an even Maass wave form $u$ 
with eigenvalue $s(1-s)$ if and only if the linear map
$\V\to\C$ defined by $(2.47)$ is invariant under the action of $\G$ on $\V\,${\rm .}
\endproclaim

\demo{Proof}  The ``if" direction is obtained by applying the linear functional $U$ to the test function $\vp_z\in \V$ defined by 
   $$\vp_z(t):=\frac{y^s}{|z-t|^{2s}}\qquad\qquad(z\in\H)\,,$$   
since $\vp_z$ transforms under ${\rm SL}(2,\R)$ by $\vp_z|g=\vp_{g(z)}$ and $U[\vp_z]=u(z)$ by (2.34).
For the other direction, it suffices to check the invariance of $U[\vp]$ under the generators $T$ and $S$. The former is obvious since replacing
$\vp$ by $\vp|T$ does not change $\widehat \vp(n)$, $n\in\Z$, and the statement for $S$ is simply (2.46).  One must also verify that the 
conditions of  Ferrar's theorem are satisfied for the pair of functions $h=\widehat\vp$, 
$h^\F=\widehat{\vp^\tau}$ for any $\vp\in \V$. We omit this.  \enddemo

We chose to prove this proposition ``by hand" by using integral transforms and the explicit generators $S$ and~$T$ of $\G={\rm PSL}(2,\Z)$,
 in accordance with the themes of this chapter.  Actually, however, the proposition has nothing to do with this particular subgroup,
but is true for any subgroup of $G$. This follows from the fact ([9] and [8], Theorem 4.29) that the ``Poisson map" 
$U\mapsto u(z)\break :=U[\vp_z]$ gives a $G$-equivariant bijection between the continuous dual of $\V$\break (= the space of distributions on
$\Bbb P^1(\R)$) and the space of eigenfunctions of $\D$ on $\H$ with eigenvalue $s(1-s)$   which have at most 
polynomial growth at the boundary.


\def\R{\Bbb R}  \def\C{\Bbb C}  \def\Cc{\C\,'} \def\Z{\Bbb Z}   \def\N{\Bbb N}  \def\z{\zeta} \def\g{\gamma}
\def\H{\Cal H} \def\G{\Gamma}  \def\D{\Delta} \def\ve{\varepsilon} \def\h{\tfrac12}   \def\sm{\smallsetminus}
\def\M{\text{\sans Maass}}    \def\p{Q} \def\FE{\text{\sans FE}_s}  \def\CR{\C\!\sm\!\R} \def\ps{\psi_s}
\def\lr{\leftrightarrow} \def\inv{^{-1}}   
\def\PL{Phragm\'en-Lindel\"of } \def\wf{\widetilde f}  \def\wU{\widetilde U}   \def\s{\sigma} \def\a{\alpha}
\def\sums{\sideset\and^*\to\sum}    \def\O{\text{\rm O}} \def\o{\text{\rm o}}  
\def\ma#1#2#3#4{\bigl(\smallmatrix#1&#2\\#3&#4\endsmallmatrix\bigr)}
\vglue12pt \centerline{\bf Chapter III. Periodlike functions}
\vglue16pt

In the previous two chapters we showed how to associate to any Maass wave form a holomorphic solution of the three-term       
functional equation (1.2) in the cut plane $\Cc$, and conversely showed that any such holomorphic solution satisfying 
suitable growth conditions is the period function associated to a Maass form.  In this chapter we investigate more fully 
the properties of general solutions of the three-term functional equation, which we call {\it periodlike} functions.  
We will be interested both in describing the totality of periodlike functions and in determining sufficient conditions for 
such a function to be the period function of a Maass form.  

As in Section~3, Chapter I, we denote the space of all (resp.~all even or odd) periodlike functions by $\FE$ (resp.~$\FE^+$ or
$\FE^-$). We write $\FE(X)$ or $\FE^\pm(X)$ for the corresponding spaces of solutions in one of the domains $X=\H$, $\,\H^-$
(lower half-plane), 
$\,\CR$, $\,\Cc$ or $\R_+$, and add the subscript ``$\omega$" to denote the subspace of analytic solutions.  We have 
\medbreak
\centerline{\BoxedEPSF{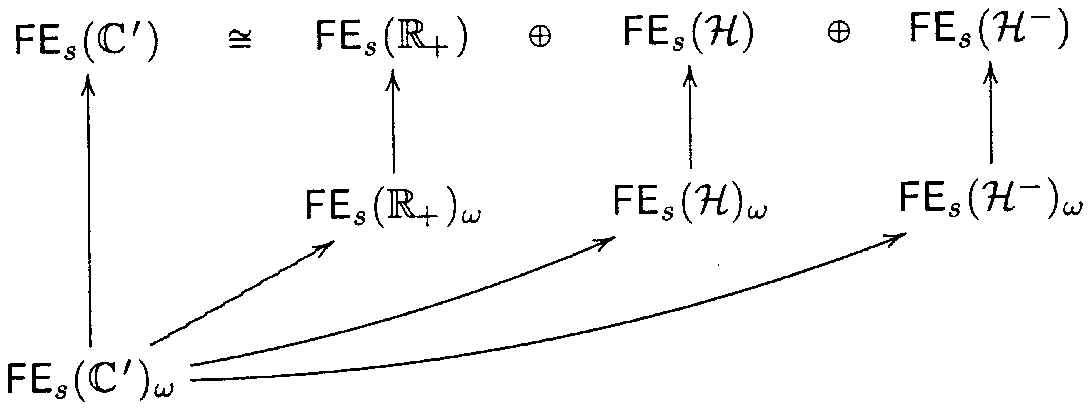 scaled 950}}
\medbreak
\noindent where the direct sum decomposition in the first line comes from the fact that if one
of the arguments
$x$,
$x+1$ or
$x/(x+1)$ of the  three-term functional equation belongs to one of the three sets $\R$, $\H$ or $\H^-$ then all three do, the
vertical arrows are inclusions,  and the remaining maps are injective because an analytic function in a connected domain is
determined by its Taylor expansion at a single point.
  
 In Section~1 we construct explicit families of periodlike functions, the constructions being general enough to show that even the 
smallest space $\FE(\Cc)_\omega$ in the above diagram is infinite-dimensional for every complex number $s$. In Section~2 we give 
a complete description of all (resp\. all continuous or all smooth) periodlike functions on $\R_+$ by a method analogous to 
that of fundamental domains in the theory of automorphic functions, i.e\. we describe various subsets $D\subset\R_+$ with the 
property that every function on $D$ is the restriction of a unique function in $\FE(\R_+)$.  

The other two sections describe the analytic properties of periodlike functions.  In Section~3 we show that every smooth 
periodlike functions on $\R_+$ has an asymptotic expansion of a very precise kind at both  0 and $\infty$.  From this description 
it follows in particular that every such function $\psi(x)$ is bounded by a specific power of $x$ (depending on $s$) in both 
directions; for instance, if $\Re(s)=\h$, the case of most interest to us, then $\psi(x)$ satisfies the estimate (0.4).
In Section~4 we perform the ``bootstrapping" described in the introduction to the paper and show that an element of $\FE(\R_+)_\omega$
satisfying an apparently only slightly stronger growth condition (e.g.~(0.3) in the case $\Re(s)=\h$) is the restriction to $\R_+$ of
an element of $\FE(\Cc)_\omega$ and comes from a Maass wave form.

\advance\sectioncount by -5
 
\section{Examples}

In this section we construct a variety of functions in $\Cc$ or $\R_+$ satisfying the three-term functional equation.

\demo{Example {\rm 1}} The first example is extremely simple: For any complex number $s$ the function
$$ \ps^-(z)=1-z^{-2s} \qquad(z\in\Cc)\tag3.1$$
is holomorphic for $z\in\Cc$ and satisfies equation (1.13) with the minus sign, so 
  $$\ps^-\;\in\;\FE^-(\Cc)_\omega\qquad(\forall s\in\C)\,.$$
The periodic function $f$ associated via the correspondences (0.5) and (0.6) to $\ps^-$ (for $s\notin\Z$) is constant.  
The function $\ps^-(z)$ vanishes identically for $s=0$ but
its derivative with respect to $s$ is the function $2\log z$, which is a  nonzero element of $\text{\sans FE}_0^-$, so we
get a holomorphic and nowhere vanishing section $s\mapsto s^{-1}\psi_s^-(z)$ of the vector bundle $\bigcup_{s\in\C}\FE^-$.
\enddemo

\demo{Example {\rm 2}}  The previous example gave an odd solution of the three-term functional equation for all $s\in\C$. 
Our second example, which is more complicated, will give an even solution for all $s$. Suppose first that $\Re(s)>1$ and define 
$$ \ps^+(z)=\sums_{m,\,n\ge0}\frac1{(mz+n)^{2s}} \qquad(z\in\Cc,\quad\Re(s)>1),\tag3.2$$
where the asterisk on the summation sign means that the ``corner" term (the one with $m=n=0$) is to be omitted
and the ``edge" terms (those with either $m$ or $n$ equal to 0) are to be counted with multiplicity 1/2.  The sum is
(locally uniformly) absolutely convergent and hence defines a holomorphic function of $z$.  Moreover,
$$ \ps^+(z+1)=\sums_{m,\,n\ge0}\frac1{(mz+m+n)^{2s}}=\sums_{n\ge m\ge0}\frac1{(mz+n)^{2s}}\, ;$$
thus
$$ \ps^+(z+1)+z^{-2s}\,\ps^+\bigl(\frac {z+1}z\bigr)
  =\biggl(\;\,\sums_{n\ge m\ge0}+\sums_{m\ge n\ge0}\,\biggr)\frac1{(mz+n)^{2s}}=\ps^+(z)\,,$$
i.e.\ $\ps^+\in\FE^+$.  We can analytically continue $\ps^+(z)$ in a standard way: write
$$ \align \G(2s)\,\ps^+(z) &=\int_0^\infty t^{2s-1}\,\sums_{m,\,n\ge0}e^{-(mz+n)t}\,dt \\ &=\int_0^\infty 
 t^{2s-1}\,\biggl(\frac14\,\frac{1+e^{-zt}}{1-e^{-zt}}\,\frac{1+e^{-t}}{1-e^{-t}}-\frac14\biggr)\,dt \endalign$$
and split up the integral into $\int_0^A+\int_A^\infty$, where $0<A<2\pi\,\min(1,|z|^{-1})$.
The second integral converges rapidly for all $s\in\C$, since the integrand is exponentially small at infinity, and hence defines a 
holomorphic function of $s$, while the first integral can be expanded as
$$ \align \int_0^A t&^{2s-1}\,\biggl(\sum\Sb m,\,n\ge0\\m+n\text{ even}\endSb
  \frac{B_m}{m!}\,\frac{B_n}{n!}\,(zt)^{m-1}\,t^{n-1}\biggr)\,dt \\  
&=\sum\Sb m,\,n\ge0\\ m+n\text{ even}\endSb\frac{B_m}{m!}\frac{B_n}{n!}\,
   z^{m-1}\,\frac{A^{m+n+2s-2}}{m+n+2s-2}\,.\endalign$$
where $B_n$ is the $n^{\rm th}$ Bernoulli number. This expression is meromorphic in $s$ with at most simple poles at $s=1$, 0, $-1$, \dots.  
All of these poles disappear if we divide by $\G(s-1)$, and the three-term functional equation is preserved by the analytic continuation, 
so we get a holomorphic section $s\mapsto (\G(2s)/\G(s-1))\,\ps^+(z)$ of the vector bundle 
$\bigcup_{s\in\C}\FE^+$.  Moreover, the special values of this section for integers $s\le1$ can be evaluated by comparing the 
residues of $\G(2s)\psi_s^+(z)$ and $\G(s-1)$ and are elementary functions of $z$, the values for $s=1$, 0 and $-1$ being (up to
constants) the functions $z^{-1}$, $\,z^{-1}-3+z$, and $z^{-1}-5z+z^3$, whose periodlike property one can verify by hand.
The function $\ps^+(z)$ for $s\notin\Z$ is the period function of the nonholomorphic Eisenstein series $E_s(z)$, while 
for integer values $s=k$ or $s=1-k$ ($k=1$, 2, \dots) it is related to the holomorphic Eisenstein series $G_{2k}(z)$.
These connections will be discussed in Sections~1 and~2 of Chapter IV, respectively. 
\enddemo

\demo{Example {\rm 3}} Generalizing Example 1, we take any periodic function $Q(z)$ and set
$$ \psi(z)=Q(z)-z^{-2s}\,Q(-1/z)\,.\tag3.3$$
Then $\psi$ belongs to $\FE$ (and in fact to $\FE^\mp$ if $Q(-z)=\pm Q(z)$).  This is formally the same construction 
as in (0.5), with $Q$ instead of $f$, but there we required $f$ to be defined and holomorphic on $\CR$, whereas now the 
interesting examples are obtained when the function $Q(z)$ is defined on $\R$ or on some neighborhood of $\R$ in $\C$.  
For instance, taking $Q(z)$ to be $\cot\pi(z-\a)$ with $\a\in\CR$ gives examples of periodlike functions on $\R_+$
which extend meromorphically to $\Cc$ but which can have poles arbitrarily close to the real axis, while taking $Q$ to 
be entire but of large growth gives solutions of (0.2) which extend analytically to all of $\Cc$ but have arbitrarily 
bad growth as $x\to0$ or $x\to\infty$. These examples show that the conclusions of the theorem in the introduction 
to this paper fail radically if no growth condition is imposed on the periodlike function $\psi|\R_+$. The fact that
$\FE(\Cc)_\omega$ is infinite-dimensional (in fact, uncountable-dimensional) for all $s$ also follows from this construction,
for example by considering the linearly independent periodic functions $Q(z)=\exp\bigl(\a\,e^{2\pi iz}\bigr)$, $\;\a\in\C$.
\enddemo

\demo{Example {\rm 4}} We can similarly try to generalize Example 2 by replacing the constant function 1 
with an arbitrary periodic function $Q(z)$.  Define 
$$ \psi(z)=\frac12\,Q(z)+\frac{z^{-2s}}2\,Q\bigl(-\frac1z\bigr)
   +\sum\Sb c,\,d>0\\(c,d)=1\endSb (cz+d)^{-2s}\,Q\bigl(\frac{az+b}{cz+d}\bigr),\tag3.4$$
where $a$ and $b$ are integers chosen so that $ad-bc=1$.
This converges if $\Re(s)>1$ (assuming that $Q$ is continuous) because the arguments of all but finitely many terms lie in a 
thin strip around the real axis and $Q$ is bounded in such a strip by virtue of its periodicity. For $Q\!\equiv\!1$ the 
function $\psi$ is just $\z(2s)^{-1}$ times the function $\ps^+$ of Example 2, and essentially the same calculation as given there
shows that also in the general case it belongs to $\FE$ (and in fact to $\FE^\pm$ if $Q(-z)=\pm Q(z)$), the domain of definition and analyticity 
properties being determined by those of $Q$.  Unfortunately, however, the construction only works in the domain $\Re(s)>1$, 
since we do not know how to give the analytic continuation of the series for general $Q$.  
\enddemo

\demo{Example {\rm 5}}  If $s$ is an integer, then the theory of periodlike functions is somewhat different since
even in the domain $\CR$ the correspondence with periodic functions described in Proposition 2, Section~2, of Chapter~I
is no longer bijective.  If $s$ is a negative integer then there are sometimes polynomial solutions of (1.2) other than $\ps^-$,
the first examples being the functions 
  $$ z(z^2-1)^2(z^2-4)(4z^2-1)\in\text{\sans FE}_{-5}^+\,,\qquad\text{and}\qquad z^2(z^2-1)^3\in\text{\sans FE}_{-5}^-$$
for $s=-5$. Such solutions correspond to cusp forms 
of weight $2-2s$ on ${\rm SL}(2,\Z)$, as will be discussed in detail in Section~2 of Chapter~IV. For positive integral values of
$s$  one also has examples of rational periodlike functions,  discovered by M\. Knopp and studied by several subsequent authors 
(see [4] and the references there), e.g.,    $$ \frac{(z^2+1)(2z^4-z^2+2)}{z(z^2-z-1)^2(z^2+z-1)^2}$$     for $s=2$.
\enddemo
 
\section{Fundamental domains for periodlike functions}  

If a group $\G$ acts on a set $X$, then the functions on $X$ with a specified transformation behavior with respect to $\G$
are completely determined by their values in a fundamental domain, and these values can be chosen arbitrarily.  The fundamental
domain is not unique, but must merely contain one point from each orbit of $\G$ on $X$.  Similarly, for the three-term
functional equation there are ``fundamental domains" such that the solutions of the equation are completely determined by
their values in that domain, which can be prescribed arbitrarily.  The next proposition describes two such fundamental
domains for the even three-term functional equation on the positive reals.  

\nonumproclaim{Proposition}  {\rm a)} Any function on the half\/{\rm -}\/open interval $[1,2)$  is the restriction of a unique
element of $\;\FE^+(\R_+)$ for every $s\in\C${\rm .} 
\smallbreak
{\rm b)} Any function on the half\/{\rm -}\/open interval $(0,\h]$  is the restriction of a unique
element of $\;\FE^+(\R_+)$ for every $s\in\C$ such that $\bigl(\frac{3+\sqrt5}2\bigr)^s\ne1\,${\rm .} \endproclaim

\demo{Proof} a) Let $\psi_0$ be the given function on [1,2) and $\psi$ any extension of it to $\R_+$ satisfying (0.1).
It suffices to consider $\psi$ on $[1,\infty)$ since $\psi(1/x)=x^{2s}\psi(x)$.  We decompose $[1,\infty)$ as 
$\bigcup_{n\ge0}I_n$ with $I_n=[n+1,n+2)$.  For $x\in I_n$ with $n\ge1$ we have $x-1\in I_{n-1}$, $\dfrac x{x-1}\in I_0$
(unless $n=1$, $x=2$). The three-term equation then shows that the restriction $\psi_n:=\psi|I_n$ satisfies 
$\,\psi_n(x)=\psi_{n-1}(x-1)-(x-1)^{-2s}\psi_0\bigl(\dfrac x{x-1}\bigr)\,$, and hence by induction on $n$ that 
$\psi|[1,\infty)$ is given by $\psi(2)=\h\psi_0(1)$ and
  $$ \psi(x)=\psi_0\bigl(x-[x]+1\bigr)-\sum_{j=1}^{[x]-1}(x-j)^{-2s}\,\psi_0\bigl(1+\frac1{x-j}\bigr)\qquad(x\ne2)\,.\tag3.5$$
This proves the uniqueness of $\psi$, but at the same time its existence, since the function defined by this formula
has the desired properties.

\medbreak
b) The proof is similar in principle, but somewhat more complicated.  Instead of $[1,\infty)=\bigcup_{n\ge0}[n+1,n+2)$ we use the
decomposition $(0,1]=\{\a,1\}\cup\,\bigcup_{n\ge0}J_n$, where $\a=\h(\sqrt5-1)$ is the reciprocal of the golden ratio
and $J_n$ is the half-open interval with endpoints $F_n/F_{n+1}$ (excluded) and $F_{n+2}/F_{n+3}$ (included), $F_n$ being
the $n^{\rm th}$ Fibonacci number.  Note that the orientation of these intervals depends on the parity of $n$, the even indices
giving intervals $J_0=(0,\h]$, $J_2=(\h,\frac35]$, \dots to the left of $\a$ and the odd indices 
giving intervals $J_1^{\phantom{|}}=[\frac23,1)$, $J_3=[\frac58,\frac23)$, \dots to the right of $\a$.  The inductive procedure\break
now takes the form $\,\psi_n(x)=x^{-2s}\psi_{n-1}\bigl(\dfrac{1-x}x\bigr) -\psi_0(1-x)\,$ for the function\break $\psi_n:=\psi|J_n$,
leading to the closed formula
$$ \align \psi(x)&=\bigl|F_nx-F_{n-1}\bigr|^{-2s}\,\psi_0\biggl(\frac{-F_{n+1}x+F_n}{F_nx-F_{n-1}}\biggr)\tag3.6\\
 &\quad -\sum_{j=1}^n\,\bigl|F_{j-1}x-F_{j-2}\bigr|^{-2s}\,\psi_0\biggl(\frac{F_{j+1}x-F_j}{F_{j-1}x-F_{j-2}}\biggr)
  \qquad\text{for $x\in J_n$, $n\ge0$.}\endalign$$
This defines $\psi(x)$ for all $x\in(0,1]$ except $\a$ and $1$. At those two points the values must be given
by $\psi(\a)=\bigl(\a^{-2s}-1\bigr)^{-1}\,\psi_0(\a^2)$ and $\psi(1)=2^{1-2s}\psi_0(\h)$, as one sees by
taking $z=\a$ and $z=1$ in (0.1) and using the relation $\psi(1/x)=x^{2s}\psi(x)$.  (This is where we use the
condition $\a^{2s}\ne1$.)  One again checks that formula (3.6), completed in this way
at the two missing points 1 and $\a$, does indeed give a periodlike extension of $\psi_0$.  \enddemo

\demo{{R}emarks} 1.  In a) and b) we could have chosen the intervals $(\h,1]$ and $[2,\infty)$ instead, since the function 
we are looking for transforms in a known way under $x\to1/x$. The proposition also holds for the odd functional equation, 
except that the fundamental domain in (a) must be taken to be $(1,2]$ instead of $[1,2)$ because 
$\psi(1)$ is now automatically 0 and does not determine $\psi(2)$. For the uniform case (i.e., solutions of (0.2))
the corresponding fundamental domains are $(\h,2]$ and $(0,\h]\cup[2,\infty)$.
\medbreak
{2.} There are many other choices of ``fundamental domains." Two simple ones, again for the even functional equation, are
 $[1+\a,2+\a)$ and $[1,1+\a]\cup(2,2+\a)$, where $\a=\h(\sqrt5-1)$ as before. The proofs are similar.

\medbreak
{3.} Related to the proof of part (b) of the proposition is the following amusing alternative form of the even three-term 
functional equation when $\Re(s)>0$ and $\psi$ is continuous:
$$ \psi(z)=\sum_{n=1}^\infty\,\bigl(F_nz+F_{n+1}\bigr)^{-2s}\,\psi\biggl(\frac{F_{n-2}z+F_{n-1}}{F_nz+F_{n+1}}\biggr)\tag3.7$$
(and similarly for the odd case, but with the $n^{\rm th}$ term multiplied by $(-1)^n\,$).

\medbreak
{4.} Another natural question is whether there are also fundamental domains $D$ in the sense of this section for 
the three-term functional equation in the cut plane $\Cc$. We could show using the axiom of choice that such domains
exist, but did not have any explicit examples. However, Roelof Bruggeman pointed out a very simple one, generalizing~(a)
 of the proposition, namely the strip $D=\{z\in\C\mid1\le\Re(z)<2\}$.  (This is for the even case; otherwise take
the union of $D$ and its image under $\tau$.)  A sketch of the argument showing that this domain works is as follows:
for $\Re(z)\ge2$ use the three-term relation to express $\psi(z)$ in terms of the values of $\psi$ at $z/(z-1)$, which is
in $D$, and $z-1$, where $\psi$ is known by induction; then use the evenness to define $\psi$ in the reflected domain
$\Re(1/z)\ge1$; and finally, in the remaining domain $\,\max\{\Re(z),\Re(1/z)\}<1$ use the three-term relation to express
the value of $\psi$ at $z$ inductively in terms of its values at $1+z$ and $1+1/z$, which are both nearer to $D$ than $z$ is. One
checks fairly easily that this uniquely and consistently defines a periodlike function on all of $\Cc$.

\medbreak
{5.} In the proposition we considered simply functions $\psi:\R_+\to\C$ defined pointwise, with no requirements of
 continuity or other analytic properties. In case~(a) a necessary condition for continuity is that the  
given function $\psi_0$ on [1,2) extends continuously to $[1,2]$ and satisfies $\psi_0(2)=\h\psi_0(1)$,
and one can check easily from formula (3.5) that this condition in fact ensures the continuity of $\psi$ everywhere.
 (Of course it suffices to check the match-up at the endpoints $x=n$ of the intervals $I_n$.) The situation for smoothness
is similar: if~$\psi_0$ extends to a $C^\infty$ function on [1,2] and the derivatives of both sides of (0.1) agree to 
all orders at $z=1$, then the extension $\psi$ defined in the proof of part (a) is $C^\infty$ everywhere. 
In a sense, this holds for (real-) analytic also: if $\psi_0$ is analytic on [1,2) and extends analytically to
(a~neighborhood of) [1,2] with the same matching conditions on its derivatives, then $\psi$ is automatically also
analytic. But this is no longer useful as a construction because we have no way of ensuring that the function defined
by analytic continuation starting from its Taylor expansion at one endpoint of the fundamental domain will satisfy
the required matching conditions at the other endpoint. 

For case (b) the story is more complicated, even for continuity. The single condition
$$\,\lim_{x\to0}\bigl(\psi_0(x)-(1+x)^{-2s}\psi_0\bigl(\frac x{1+x}\bigr)\bigr)=2^{1-2s}\psi_0(\h)$$
on the function $\psi_0:(0,\h]\to\C$ is enough to ensure the matching of the functions $\psi_n$ and $\psi_{n+2}$ at the common
endpoint of their intervals of definition, but the resulting function on $(0,\a)\cup(\a,1]$ will in general be highly
discontinuous at the limit point~$\a$ of these intervals. Requiring continuity at $\a$ imposes far severer restrictions
on $\psi_0$, namely, that the identity (3.7) should hold (with $\psi$ replaced by $\psi_0$) for all $z\in(0,\h]$. This has
the remarkable consequence that a ``fundamental domain" for continuous functions is {\it smaller} than for arbitrary 
functions.  For instance, if $\psi$ is assumed to be continuous then its values on $(0,\sqrt2-1]$ already determine it,
since for $\sqrt2-1<z\le\h$ all the arguments on the right-hand side of (3.7) are less than $\sqrt2-1$; and similarly
$\psi$ is completely determined by its values on $[\a^2,\h]$ because if $0<x<\a^2$ then the $n=1$ term in (3.7) for
$z=x/(1-x)$ has argument $x$ and all the other terms have arguments strictly between $x$ and $\h$.

\section{Asymptotic behavior of smooth periodlike functions} 

 In this section will study smooth solutions of the three-term functional equation (0.2) on $\R_+$. (By ``smooth" we will always
mean $C^\infty$.)  We denote the space of such functions by $\FE(\R_+)_\infty$.  As usual, we use $\s$ to denote $\Re(s)$.

 \nonumproclaim{Proposition} Let $s\in\C\sm\{\h,\,0,\,-\h,\,\ldots\}${\rm .}
 Then any function $\psi\in\FE(\R_+)_\infty$ has asymptotic expansions of the form
$$\align
\noalign{\vskip8pt}
&\tag 3.8a\\
\noalign{\vskip-26pt} \qquad
\psi(x)\;&\sim\;x^{-2s}\,Q_0\bigl(\frac1x\bigr)\;+\,\sum_{m=-1}^\infty C_m^*\,x^m\;\,\qquad\qquad\text{as $x\to0$},\\
  \psi(x)\;&\sim\;Q_\infty(x)\;+\sum_{m=-1}^\infty (-1)^{m+1}\,C_m^*\,x^{-m-2s}\quad\,\text{as $x\to\infty$},\tag 3.8b
\endalign$$ where $Q_0,\;Q_\infty:\,\R\to\C$ are smooth periodic functions{\rm .}  The coefficients $C_m^*$ are given in terms of
the Taylor coefficients 
$C_n=\psi^{(n)}(1)/n!$ of $\psi(x)$ at $x=1$ by 
  $$ C_m^*=\frac1{m+2s}\,\sum_{k=0}^{m+1}(-1)^k\,B_k\,\binom{m+2s}k\,C_{m+1-k}\qquad(m\ge-1)\,,\tag3.9$$
where $B_k$ is the $k^{\rm th}$ Bernoulli number{\rm .}  \endproclaim

\nonumproclaim{{C}orollary} Let $s$ and $\psi$ be as above{\rm .}  Then\/{\rm :}
$$\psi(x)\,\text{ is $\;\O\bigl(x^{-\max(2\s,1)}\bigr)\,$ as $x\to0$, $\quad\O\bigl(x^{\max(0,1-2\s)}\bigr)\,$ as $x\to\infty$}.\tag3.10$$ 
\endproclaim

\demo{Proof}  We start with the expansion near 0.  Suppose first that $\s>\frac12$, and define $Q_0(x)$ for $x>0$ by
$$Q_0(x)=\frac1{x^{2s}}\,\psi\bigl(\frac1x\bigr)-\sum_{n=0}^\infty\frac1{(n+x)^{2s}}\,\psi\bigl(1+\frac1{n+x}\bigr)\,.\tag 3.11$$
Clearly this converges to a smooth function on $\R_+$. Moreover,
$$Q_0(x)-Q_0(x+1)=\frac1{x^{2s}}\,\psi\bigl(\frac1x\bigr)-\frac1{(x+1)^{2s}}\,\psi\bigl(\frac1{x+1}\bigr)
   -\frac1{x^{2s}}\,\psi\bigl(1+\frac1x\bigr)=0\,,$$
so $Q_0(x)$ is periodic. In the general case we replace (3.11) by
$$\align 
\noalign{\vskip8pt}
&\tag3.12\\
\noalign{\vskip-28pt} \qquad Q_0(x)\;=\;&x^{-2s}\,\psi\bigl(\frac1x\bigr)-\sum_{m=0}^M C_m\,\zeta(m+2s,x) \\ &\quad-\,\sum_{n=0}^\infty
  \frac1{(n+x)^{2s}}\,\biggl(\psi\bigl(1+\frac1{n+x}\bigr)-\sum_{m=0}^M\frac{C_m}{(n+x)^m}\biggr)\,,\endalign$$
where $M$ is any integer with $M+2\s>0$ and  $\zeta(a,x)$ is the Hurwitz zeta function, defined for $\Re(a)>1$ as 
$\sum\limits_{n=0}^\infty\dfrac1{(n+x)^a}$ and for other values of $a\ne1$ by analytic continuation.  
The infinite sum in (3.12) converges absolutely and one checks easily that the definition is independent of $M$, agrees with the previous 
definition if $\s>\frac12$, and is again periodic.  From the estimate $\zeta(a,x)=\text O(x^{1-a})$ we obtain that
$$\frac1{x^{2s}}\,\psi\bigl(\frac1x\bigr)=Q_0(x)+\sum_{m=0}^M C_m\,\zeta(2s+m,x)\;+\;\text O\bigl(x^{-2s-M}\bigr)\qquad(x\to\infty)$$
for any integer $M$ with $M+2\s>0$, and if we now use the full asymptotic expansion of the Hurwitz zeta function 
$$\zeta(a,x) \sim \frac1{a-1}\;\sum_{k\ge0}(-1)^k\,B_k\,\binom{k+a-2}k\,x^{-a-k+1}\qquad(x\to\infty)\,,$$  
which is an easy consequence of the Euler-Maclaurin summation formula, we obtain
the asymptotic expansion (3.8a) with coefficients $C_m^*$ defined by (3.9). 

The analysis at infinity is similar.  If $\s>\frac12$, then we define
  $$Q_\infty(x)=\psi(x)-\sum_{n=1}^\infty(n+x)^{-2s}\,\psi\bigl(1-\frac1{n+x}\bigr)\tag 3.13$$
and find from (0.2) as before that $Q_\infty(x+1)=Q_\infty(x)$.  The same trick as before permits us to define $Q_\infty$
also when $\s\le\frac12$, obtaining the formulas
  $$\psi(x)= Q_\infty(x)+\sum_{m=0}^M(-1)^m\,C_m\,\zeta(m+2s,x+1)+\text O\bigl(x^{-2s-M}\bigr)\qquad(x\to\infty)$$
for $M$ sufficiently large. Now using the asymptotic expansion of $\z(a,x+1)$, which is identical with that of $\z(a,x)$ 
but without the factor $(-1)^k$, we get (3.8b).  The corollary follows because any smooth periodic function on $\R$ is bounded. 
\enddemo

\demo{Remarks} 1. In principle it would have sufficed to treat just one of the expansions at 0 and $\infty$, since
if $\psi(x)$ satisfies the three-term functional equation, then so does $\psi^\tau(x)=x^{-2s}\,\psi(1/x)$,
and replacing $\psi$ by $\psi^\tau$ simply interchanges the two asymptotic formulas in (3.8), with the roles of $Q_0$ and $Q_\infty$ 
exchanged and $C_m^*$ multiplied by $(-1)^{m+1}$. However, it is not obvious (though it will become so in the next remark) that 
replacing $\psi$ by $\psi^\tau$ changes the Taylor coefficients $C_n$ in such a way as to multiply the right-hand side of (3.9) 
by $(-1)^{m+1}$, so it seemed easiest to do the expansions at 0 and $\infty$ separately.

\medbreak {2.} There is a more direct way to relate the coefficients of the asymptotic series in (3.8) to the Taylor coefficients
$C_n$. For convenience of notation we write the expansions (3.8) as
  $$ \psi(x)\,\underset{x\to0}\to\sim\,\frac{Q_0(1/x)}{x^{2s}}\,+\,P_0(x)\,,\qquad 
    \psi(x)\,\underset{x\to\infty}\to\sim\,Q_\infty(x)\,+\,\frac{P_\infty(1/x)}{x^{2s}}$$
where $P_0(t)$ and $P_\infty(t)$ are Laurent series in one variable. Similarly, we write
  $$ \psi(x)\,\underset{x\to1}\to\sim\,P_1(x-1)\,,\quad\text{where}\quad P_1(t)=\sum_{n=0}^\infty C_nt^n\in \C[[t]]\,,$$
for the asymptotic expansion of $\psi$ around 1. (Notice that saying that a function is $C^\infty$ at a point is equivalent to saying that
it has an asymptotic power series expansion at that point.) Now, setting $x=t$ or $x=-1-1/t$ in the functional equation (0.2) and letting
$t$ tend to 0 from the right or left, respectively, we find that the three formal Laurent series $P_0$, $P_\infty$ and $P_1$ are related by
  $$ P_0(t)-(1+t)^{-2s}P_0\bigl(\frac t{1+t}\bigr)\,=\,P_1(t)\,=\,(1+t)^{-2s}P_\infty\bigl(\frac{-t}{1+t}\bigr)-P_\infty(-t)\,.$$
Denoting the coefficients of $P_0$ by $C_m^*$ and expanding by the binomial theorem, we find that the first of these equations is
equivalent to the system of equations
  $$  C_n = \sum_{m=-1}^{n-1} (-1)^{n-m-1}\,\binom{n-1+2s}{n-m} \,C_m^*\qquad\quad(n\ge0)\,,\tag 3.14$$
while the second equation is the identical system but with $C_m^*$ replaced by the coefficient of $t^m$ in $-P_\infty(-t)$.
But the system (3.14) is clearly invertible, since the leading coefficients $n-1+2s$ are nonzero for all $n$.  It follows that
the $C_n^*$ are uniquely determined in terms of the $C_n$ and hence also that $P_\infty(t)=-P_0(-t)$; i.e., that the $m^{\rm
th}$ coefficient of $P_\infty$ is $(-1)^{m+1}C_m^*$ as asserted in (3.8b).  Finally, to see that the inversion of the system
(3.14)  is given explicitly by (3.9) we rewrite these equations in terms of generating functions as
$$ \sum_{n\ge0}^\infty\frac{C_n}{\G(n+2s)}\,w^n\;=\;\bigl(1-e^{-w}\bigr)\,\sum_{m=-1}^\infty\frac{C_m^*}{\G(m+2s)}\,w^m\,,\tag3.15$$
which can be inverted immediately by multiplying both sides by $(1-e^{-w})^{-1}=\sum B_kw^{k-1}/k!\,$.  The reader may
recognize formula (3.15) as being identical (in the case of the expansions associated to an even Maass form $u$) with the relation
between the functions $g(w)$ and $\phi(w)$ discussed in Chapter II, Section~4 (compare equations (2.29), (2.23) and (2.32)),
while equation (2.33) of that section identifies the coefficients $C_m^*$ in the Maass case as special values of the $L$-series
associated to $u$.

\medbreak{3.} Finally, we say a few words about the case excluded so far when $2s=1-h$ for some integer $h\ge0$.  Equation
(3.12) defining 
$Q_0$ (the case of $Q_\infty$ is similar and will not be mentioned again) no longer makes sense since the term $C_h\z(1,x)$ 
occurring in it is meaningless. However, the function $\;\lim_{a\to1}\bigl(\z(a,x)-\z(a)\bigr)$ exists and equals 
$-\G'(x)/\G(x)-\g$ for some constant $\g$ (namely, Euler's).  If we modify the definition (3.12) by replacing the
undefined term $C_h\z(1,x)$ by $-C_h\G'(x)/\G(x)\,$, then we find that $Q_0$ is again smooth and periodic but that the expansion
(3.8b) must be changed by replacing the term $C_{h-1}^*x^{h-1}$ by $(C_h\log x +C_{h-1}^*)\,x^{h-1}$.  The value of $C_{h-1}^*$ here
is arbitrary, since we can change it simply by adding a constant to the periodic function $Q_0$. (This corresponds to the
arbitrary choice of additive constant in our renormalization of $\z(1,x)$.)  Everything else goes through as before, and
(3.10) is also unaffected except in the case $s=\h$, when the logarithmic term becomes dominant.  We can summarize the growth
estimates of $\psi$ for arbitrary $s\in\C$ by the following table:

\centerline{\BoxedEPSF{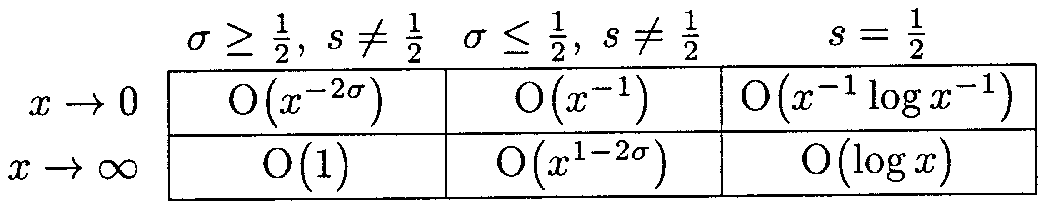 scaled 1000}}
 
We illustrate the proposition by describing the expansions (3.8) for each of the special periodlike functions given in Section~1.

\demo{Example {\rm 1}} Here $Q_\infty(x)\equiv1$,  $Q_0(x)\equiv-1$, $C_m^*=\delta_{m,0}$, and $C_n$ equals 0 
for $n=0$ and $-\binom{-2s}n$ for $n\ge1$, in agreement with (3.14).
\enddemo

\demo{Example {\rm 2}} Here $\psi=\psi^\tau_{\phantom{|}}$, so we must have $Q_0=Q_\infty$ and $C_m^*=\break (-1)^{m+1}C_m^*$, i.e.\
$C_m^*=0$ for
$m$ even.  A simple calculation shows that\break $Q_0(z)=\h\z(2s)$ (constant function) and 
 $\dsize C_m^*=B_{m+1}\binom{m+2s}{m+1}^{\phantom{\int}}\!\!\dfrac{\z(m+2s)}{m+2s}$ for odd $m\ge-1$.
\enddemo

\demo{Example {\rm 3}} Here (assuming that the function $Q$ is smooth) $Q_\infty(x)=Q(x)$, $Q_0(x)=-Q(-x)$, $\;C_{-1}^*=0$ 
and $C_m^*=Q^{(m)}(0)/m!$ for $m\ge0$.
\enddemo

\demo{Example {\rm 4}} Here we find $Q_0(x)=\h\,Q(-x),\quad Q_\infty(x)=\h\,Q(x)$.  The coefficients $C_m^*$ can be calculated 
as special values of Dirichlet series; e.g.
  $$ C_{-1}^*=\sum_{c=1}^\infty\,\sum\Sb \text{$a$ mod $c$}\\(a,c)=1\endSb Q\bigl(\frac ac\bigr)\;c^{-2s}\,,$$
where the series converges because $\Re(s)>1$ and $Q$ is bounded.
\enddemo

\demo{Example {\rm 5}} For the examples of polynomial periodlike functions $\psi(x)$ for negative integral values of $s$ (like the two 
given for $s=-5$), we can simply take $Q_0=Q_\infty=0$ and define $C_m^*$ as the $m^{\rm th}$ coefficient of the polynomial $\psi$, so that $C_m^*=0$
for $m=-1$ or $m>2-2s$. Then (3.8a) is obvious and (3.8b) follows from the property $\psi(x)=-x^{-2s}\psi(-1/x)$ which always hold for such
polynomial solutions.  Of course this is not quite unique, since we could also take $Q_0(x)=a$, $Q_\infty(x)=-a$ for any constant $a$ and
change the value of the coefficient $C^*_{-2s}$ by $a$.  In the case of rational periodlike functions for positive integral $s$
we again take $Q_0=Q_\infty=0$ and define the coefficients from the expansion of $\psi$ at either 0 or $\infty$ (e.g. $C_{-1}^*=2$, $C_1^*=13$,
$C_3^*=57,\,\ldots\,$ and $C_0^*=C_2^*=\ldots=0$ for the example given at the end of \S 2 with $s=2$), the results of the two
calculations agreeing because we again always have the invariance property $\psi(x)=-x^{-2s}\psi(-1/x)$.
\medbreak

Let us return to the general case and suppose again that $s\hskip-.65pt \notin\hskip-.65pt \{\h,0,-\h,\ldots\}$. Then there is no ambiguity in
the  decomposition of $\psi$ into a ``periodic" and an ``asymptotic" part at 0 or $\infty$, so that equations (3.8) give a well-defined map
  $$\align \FE(\R_+)_\infty\;&\longrightarrow\;C^\infty(\R/\Z)\,\oplus\,C^\infty(\R/\Z)\tag3.16\\ 
\noalign{\vskip6pt}
  \psi\quad &\;\,\mapsto\qquad\;(Q_0,\quad Q_\infty)\,.\endalign$$
A natural question is whether this map is surjective. The construction of Example~3 shows that all pairs $(Q_0,Q_\infty)$ of the form
$(Q(x),\,-Q(-x))$ are in the image. The construction of Example~4 gives the complementary space $\{(Q(x),Q(-x))\}$, and hence the full
surjectivity, if $\s>1$.  This construction also works in the analytic category and will be essentially all we can prove about surjectivity there 
(cf.~\S4). For $\s<1$ we do not know whether the surjectivity is true.  It would suffice to construct an even $C^\infty$ periodlike function
$\psi=\psi^\tau$ having a given periodic function $Q$ as its~$Q_0$. A simple construction using the second fundamental domain construction
(part (b) of the proposition) of Section~2 produces an infinity of even periodlike functions with given (smooth) $Q_0$ which are $C^\infty$
except at the two points $\a=\h(\sqrt5-1)$ and $\a^{-1}$, but we do not know whether any of them can be made smooth at these two points.

\section{``Bootstrapping"}

\nonumproclaim{Theorem 2} Let $s$ be a complex number with $\s>0$ and $\psi$ any real\/{\rm -}\/analytic solution  of the
functional equation $(0.2)$ on $\R_+$ such that\/{\rm :}\/
$$\psi(x)\,\text{ is $\;\o\bigl(x^{-\min(2\s,1)}\bigr)\,$ as $x\to0$, $\quad\o\bigl(x^{\min(0,1-2\s)}\bigr)\,$ as $x\to\infty$}.\tag3.17$$ 
Then $\psi$ extends holomorphically to all of $\Cc$ and satisfies the growth conditions
$$\psi(z)\ll\cases|z|^{-2\s}&\text{if $\Re(z)\ge0,\;|z|\ge1$},\\ \;1&\text{if $\Re(z)\ge0,\;|z|\le1\,,$}\\
     |\Im(z)|^{-\s}&\text{if $\Re(z)<0\,.$}   \endcases\tag3.18$$    \endproclaim 

\nonumproclaim{{C}orollary} Under these hypotheses{\rm ,} $\psi$ is the period function associated to a Maass wave form with eigenvalue $s(1-s)$.
In particular{\rm ,} if $\s=\h$ then any bounded and real\/{\rm -}\/analytic even or odd periodlike function comes from a Maass form{\rm .} \endproclaim

The corollary follows from Theorem 2 and from Theorem 1 of Chapter I, Section~1, since the estimates (3.18) are the same as (1.3) with
$A=\s$. Observe that the growth estimate (3.17) which suffices to imply that $\psi$ comes a Maass form differs from the estimate (3.10)
which  holds automatically for smooth periodlike functions only by the replacement of ``max" by ``min" and of ``O" by ``o".  In particular,
in the case of most interest when $\s=\h$, merely changing ``O" to ``o" suffices to reduce the uncountable-dimensional space 
$\FE(\R_+)_\omega$\break (cf\. \S1) to the finite- (and usually zero-) dimensional space of period functions of Maass wave forms, as
already  discussed in the introduction to the paper.

\demo{Proof} We first observe that, since both the hypothesis (3.17) and the conclusion (3.18) are invariant under $\psi\mapsto\psi^\tau$,
we can split $\psi$ into its even and odd parts $\h(\psi\pm\psi^\tau)$ and treat each one separately.  We therefore can (and will) assume
that $\psi$ is either invariant or anti-invariant under $\tau$.  This is convenient because it means that we can restrict our attention to
either $z$ with $|z|\le1$ or $|z|\ge1$ (usually the latter), rather than having to consider both cases.  In particular, we only need to use 
one of the two estimates (3.17), and only have to prove the first and the last of the inequalities (3.18).

We begin by proving the analytic continuation.  The key point is that the estimates (3.17) imply that the periodic function $Q_\infty$
and the coefficient $C_{-1}^*$ in (3.8b) (as well, of course, as the periodic function $Q_0$ in (3.8a)) vanish.  From (3.9) or (3.14) it
then follows that $\psi(1)$ ($=C_0$) also vanishes.  This in turn implies that equation (3.13) holds even if $\s$ is not bigger than $\h$ 
(recall that $\s>0$ by assumption), so that the vanishing of $Q_\infty$ gives the identity
   $$ \psi(z)\,=\,\sum_{n=1}^\infty(n+z)^{-2s}\,\psi\bigl(1-\frac1{n+z}\bigr)\tag 3.19$$
for all $z$ on (and hence for all $z$ in a sufficiently small neighborhood of) the positive real axis.  
This formula will play a crucial role in what follows.

We first show that $\psi$ extends to the wedge $W_\delta=\{z:|\arg(z)|<\delta\}$ for some $\delta>0$.
The function $\psi$ is already holomorphic in a neighborhood of $\R_+$ and in particular in the disk $|z-1|<\ve$ for some $\ve$.  Formula
(3.19) then defines $\psi$ (as a holomorphic function extending the original values) in the subset $\,\{z\in\Cal R:\,|z|>1/\ve\}$ of the
right half-plane $\,\Cal R=\{z\in\C:\,x>0\}$.  (We use the standard notations $x$ and $y$ for the real and imaginary parts of $z$.)   
The assumed invariance or anti-invariance of $\psi$ under $\tau$ gives us $\psi$ also in the half-disk $\,\{z\in\Cal R:\,|z|<\ve\}$.  
Since the interval $[\ve,\,1/\ve]$ is compact, this suffices to define $\psi$ in a wedge of the form stated.
\comment
We first show that $\psi$ extends to the wedge $W_\delta=\{z:\,|\arg(z)|<\delta\}$ for some $\delta>0$.
The function $\psi$ is holomorphic at $z=1$ and hence in the disk $|z-1|<\ve$ for some $\ve>0$.  Formula
(3.19) then defines $\psi$ (as a holomorphic function extending the original values) in the subset $\,\{z\in\Cal R:\,|z|>1/\ve\}$ of the
right half-plane $\,\Cal R=\{z\in\C,\,x>0\}$.  (We use the standard notations $x$ and $y$ for the real and imaginary parts of $z$.)   
The assumed invariance or anti-invariance of $\psi$ under $\tau$ gives us $\psi$ also in the half-disk $\,\{z\in\Cal R,\,|z|<\ve\}$.  
Since $\psi$ is already holomorphic in some neighborhood of the compact interval $[\ve,\,1/\ve]$ is compact, this suffices to define
$\psi$ in some wedge $W_\delta$.
\endcomment

We next show that any periodlike function in $W_\delta$ (even without any growth or continuity assumptions) extends uniquely to a periodlike 
function on all of $\Cc$. Denote by $\p$ the sub-semigroup of $\,{\rm SL}(2,\Z)\,$ generated by the matrices $T=\ma1101$ and $T'=\ma1011$. Note that any
element of $\p$ has nonnegative entries. (In fact, $\p$ consists of all matrices in $\,{\rm SL}(2,\Z)\,$ with nonnegative entries, but we will not
need this fact.) We claim that for any $z\in\Cc$ there are only finitely many elements $\g=\ma abcd$ of $\p$ for which $\g(z)=(az+b)/(cz+d)$
does {\it not} belong to $W_\delta$.  Indeed, this is obvious for the matrices $\ma1b01$ with $c=0$, and no element of $\p$ has $a=0$, so
we can assume that $ac\ne0$.  Then
$$ \arg\bigl(\g(z)\bigr)\,=\,\arg\biggl(\frac ca\cdot\frac{az+b}{cz+d}\biggr)=\arg\biggl(1-\frac1{a(cz+d)}\biggr)\,.$$
Choose $M$ so large that $W_\delta$ contains a $(1/M)$-neighborhood of 1.  There are only finitely many pairs of integers $(c,d)$ with
$|cz+d|<M$, and for each such pair only finitely many integers $a$ with $0<a<M/|cz+d|$, which proves the claim.  Now let $\p_n$ ($n\ge0$) be the
subset of $\p$ consisting of words in $T$ and $T'$ of length exactly $n$. By what we just showed, we know that for any $z\in \Cc$ and $n$
sufficiently large (depending on $z$) we have $\g(z)\in W_\delta$ for all $\g\in \p_n$.  We choose such an $n$ and define
 $$ \psi(z)=\sum_{\g\in \p_n}\,(\psi|\g)(z)\,,\tag 3.20$$
where $(\psi|\g)(z):=(cz+d)^{-2s}\psi(\g(z))$ for $\g=\ma abcd\in \p$. (The power $(cz+d)^{-2s}$ is well-defined because $cz+d\in\Cc$.) 
The right-hand side is well-defined because each argument $\g(z)$ is in the domain $W_\delta$ where $\psi$ is already defined; it
is independent of $n$ because $\psi$ satisfies the  three-term functional equation $\psi=\psi|T+\psi|T'$ in $W_\delta$ and 
$\p_{n+1}=\p_nT\bigsqcup \p_nT'$ (disjoint union); and it satisfies the three-term functional equation because $\p_{n+1}=T\p_n\bigsqcup T'\p_n$.  
It is also clear (since the sum in (3.20) is finite and one can choose the same $n$ for all $z'$ in a neighborhood of $z$) that this 
extension is holomorphic if $\psi|W_\delta$ is holomorphic. This completes the proof of the analytic continuation.

We now turn to the proof of the estimates (3.18). We again proceed in several steps. In principle the proof of the estimates mimics the
proof of the analytic continuation, with an inductive procedure to move outwards from the positive real axis to all of $\Cc$. However,
if we merely estimated $\psi$ in the wedge $W_\delta$ and used (3.20) directly, we would get a poorer bound than we need, so we have
to organize the induction in a more efficient manner.  

We start by estimating $\psi$ away from the cut. The first inequality in (3.18) follows immediately from formula (3.19) and the vanishing 
of $\psi$ at 1. In fact, this gives the full asymptotic expansion $\psi(z)\sim\sum\limits_{m=0}^\infty (-1)^{m+1}C_m^*z^{-m-2s}$ as 
$|z|\to\infty$ in $\Cal R$, with the same $C_m^*$ as in~(3.8b). (Note that $|z^{-2s}|\ll|z|^{-2\s}$ because $\,\arg(z)\,$ is bounded.) 
The same method works to estimate $\psi(z)$ for\break $z\in\Cal L$ (left half-plane) with $|y|>\h$, the bound obtained now being $|y|^{-2\s}$ 
instead of $|z|^{-2\s}$.  This proves the third estimate in (3.18), and in fact a somewhat sharper bound, in this region.

It remains to consider the region $X=\{z\in\Cal L:\,|z|\ge1,\;|y|\le\h\}$.  (Recall that we can restrict to $|z|\ge1$ because $\psi^\tau=\pm\psi$.)
To implement the induction procedure in this region we define a map $J:X\to\{z\in\Cc:\,|z|\ge1\}$ which moves any $z\in X$ 
further from the cut. This map sends $z$ to $\dfrac{z+N'}{z+N}$, where $-N$ is the nearest integer to $z$ and $-N'$ the second 
nearest, or more explicitly 
$$ J(z)=\cases\dfrac{z+N-1}{z+N}&\quad\text{if $\;-N\le x<-N+\h\,,\quad N\ge1\,,$}\\
  \dfrac{z+N+1}{z+N}&\quad\text{if $\;-N-\h\le x<-N,\quad N\ge1\,.$}\endcases\tag3.21$$
We claim that    
  $$ \psi(z)\,=\,\pm(z+N)^{-2s}\,\psi\bigl(J(z)\bigr)\,+\,\O(1)\,.\tag3.22$$
(Here and in the rest of the proof all O-estimates are uniform, depending only on~$\psi$.) This follows in the first case of (3.21) from the calculation
$$ \align \psi(z)-(z+N)^{-2s}\,\psi\bigl(J(z)\bigr) &=\psi(z+N)+\sum_{n=1}^{N-1}(z+n)^{-2s}\,\psi\bigl(1-\frac1{z+n}\bigr)\\
   & =\O(1)+\sum_{n=1}^{N-1}\O\bigl((N-n-\h)^{-2\s-1}\bigr)\endalign$$
 and in the second case from a similar calculation using
 $$ \psi(z)\mp(z+N)^{-2s}\,\psi\bigl(J(z)\bigr)=\psi(z+N+1)\pm\sum_{n=1}^{N-1}(z+n)^{-2s}\,\psi\bigl(1+\frac1{z+n}\bigr)\,.$$
Multiplying both sides of (3.22) by $y^s$ and noting that $\,\Im(J(z))=y/|z+N|^2$, we obtain the estimate
 $$ F(z)\;=\;F(J(z))\,+\,\O\bigl(|y|^\s\bigr)\tag3.23$$
for the function $F(z):=|y|^\s\,|\psi(z)|$. On the other hand, for $z\in X$ we have
  $$ \bigl|\Im\bigl(J(z)\bigr)\bigr| =\frac{|y|}{|z+N|^2}=\frac{|y|}{(x+N)^2+y^2}\ge\frac{|y|}{\frac14+\frac14}=2\,\bigl|\Im(z)\bigr|\,.\tag3.24$$  
The desired bound $F(z)=\O(1)$ follows easily from (3.23) and (3.24).  Indeed, (3.24) implies that some iterate $J^k$ of $J$ maps $z$
outside $X$. Write $z_0=z$ and $z_j=J(z_{j-1})$ for $1\le j\le k$, so that $z_0,\ldots,z_{k-1}\in X$ and $z_k\notin X$ but $|z_k|\ge1$. 
We have 
 $$  F(z)\,=\,\O\bigl(|y_0|^\s\bigr)+\O\bigl(|y_1|^\s\bigr)+\cdots+\O\bigl(|y_{k-1}|^\s\bigr)+F(z_k)\quad\;\bigl(y_j:=\Im(z_j)\bigr)\tag3.25$$
by (3.23). But $F(z_k)=\O(1)$ by the estimates away from the cut proven earlier, and $|y_j|\le 2^{-k+j}$ 
for $0\le j\le k-1$ by (3.24). The result follows. \enddemo

\demo{{R}emark} From the proof of Theorem 2 it is clear that we could weaken the hypothesis (3.17) somewhat and still obtain the same
conclusion: it would suffice to assume that 
$$\psi(x)=\o\bigl(x^{-2\s}\bigr)\,\text{ as $x\to0$, $\qquad\psi(x)=\o\bigl(1\bigr)\,$ as $x\to\infty$}$$ 
(which is weaker than (3.17) if $\s>\h$) and that $\psi(1)=0$, since these hypotheses would already imply that $Q_0\equiv Q_\infty\equiv0$,
$C_{-1}^*=0$ in (3.8) and this is all that was used in the proof.  The contents of Theorem 2 and its corollary can therefore be summarized by
saying that the sequence
  $$0\;\longrightarrow\;\M_s\;\overset\a\to\longrightarrow\;\FE(\R_+)_\omega\;\overset\beta\to\longrightarrow\;
     C^\omega(\R/\Z)\,\oplus\,C^\omega(\R/\Z)\oplus\C\;,$$
where $\a$ sends a Maass wave form to its period function and $\beta$ sends a real-analytic periodlike function $\psi$ to
$(Q_0,Q_\infty,C_{-1}^*)\,$, is exact.  As in Section~3, it is natural to ask whether $\beta$ is surjective. The same constructions
as used there for the $C^\infty$ case (Examples~3 and~4 of~\S1) 
show that the image of $\beta$ contains at least one copy $\,\{Q(x),\,-Q(-x)\}\,$ of $C^\omega(\R/\Z)$, and that it contains all of
$C^\omega(\R/\Z)\,\oplus\,C^\omega(\R/\Z)$ if $\s>1$.  But the latter case is not very interesting, since then $\M_s=\{0\}$,
and for $0<\s<1$ we do not know how to decide the surjectivity question.
 

\def\Z{\Bbb Z}    \def\R{\Bbb R}  \def\C{\Bbb C}  \def\z{\zeta} 
\def\H{\Cal H} \def\V{\bold V}  \def\h{\frac12} \def\D{\Delta} \def\L{\Cal L} 
\def\G{\Gamma} \def\g{\gamma}  \def\N{\Cal N}
\def\a{\alpha} \def\s{\sigma}    
\def\ma#1#2#3#4{\pmatrix #1&#2\\#3&#4\endpmatrix} \def\abcd{\ma abcd} 
\def\sma#1#2#3#4{\bigl(\smallmatrix #1&#2\\#3&#4\endmatrix\bigr)}

\def\={\,=\,} \def\O{\text O} \def\o{\text o} \def\Dk{\Bbb D} \def\bD{\overline{\Dk}} \def\Cm{\C\,'} \def\lr{\leftrightarrow}
\def\sums{\sideset\and^*\to\sum}    \def\O{\text{\rm O}} \def\o{\text{\rm o}}  \def\bs{\gtrless}
\def\p{\partial} \def\wf{\widetilde f} \def\|#1{\bigl|_{#1}\bigr.\,}  \def\wU{\widetilde U} 

\def\R{\Bbb R}  \def\C{\Bbb C}  \def\Cc{\C\,'} \def\Z{\Bbb Z}   
\def\H{\Cal H} \def\G{\Gamma}  \def\D{\Delta} \def\ve{\varepsilon} \def\h{\tfrac12}  \def\sm{\smallsetminus}
\def\e{\doteq} \def\lr{\leftrightarrow} \def\inv{^{-1}}
\vglue24pt \centerline{\bf Chapter IV. Complements}
\vglue16pt

In this chapter we describe the extension of the period theory to the noncuspidal case, its connection with
periods of holomorphic modular forms, and its relationship to Mayer's theorem expressing the Selberg zeta
function of $\G$ as a Fredholm determinant.  Other ``modular" aspects of the theory (such as the action of Hecke
operators and the Petersson scalar product) will be treated in Part II of the paper.

\advance\sectioncount by -4
\section{The period theory in the noncuspidal case}

So far we have been considering only cuspidal Maass forms.  If we define more generally a Maass form with spectral
parameter $s$ to be a $\G$-invariant solution of the equation $\D u=s(1-s)u$ in the upper half-plane which grows
less than exponentially as $y\to\infty$, then the space of such forms for a given value of $s$ contains the 
cusp forms as a subspace of codimension~1, the extra function being the nonholomorphic Eisenstein series $E_s(z)$.
This is the function defined for $\Re(s)>1$ by
  $$ E_s(z) \= \frac12\!\sum_{(m,n)\in\Z^2\sm\{(0,0)\}}\frac{y^s}{|mz+n|^{2s}} \tag 4.1 $$
and for arbitrary $s$ by the Fourier expansion
  $$ \align  \noalign{\vskip8pt}
&\tag 4.2\\
\noalign{\vskip-26pt}
\qquad E_s(z) &\= \z(2s)\,y^s\,+\,\frac{\pi^{1/2}\G(s-\h)}{\G(s)}\,\z(2s-1)\,y^{1-s} \\
   &\quad+\frac{4\pi^s}{\G(s)}\,\sqrt y\,\sum_{n=1}^\infty n^{\h-s}\s_{2s-1}(n)\,K_{s-\h}(2\pi ny)\,\cos(2\pi nx)\,,\endalign$$
where $\s_\nu(n)=\sum_{d|n}d^\nu$. 
If our theory of period functions of Maass forms is to extend to the noncuspidal case, there should therefore be a
solution $\psi$ of the (even) three-term functional equation associated to $E_s$.  We claim that this is the case,
with $\psi$ being the function $\psi_s^+$ introduced in Example 2 of Section~1, Chapter III.  There are three ways of seeing this:
\medbreak
1. Substitute $u(z)=E_s(z)$ into the integral (2.2).  From (4.2) we see that $E_s(iy)=\O\bigl(y^{\max(\s,1-\s)}\bigr)$
as $y\to\infty$, and the invariance under $y\mapsto1/y$ gives the corresponding statement as $y\to0$, so the integral
converges for all $s$ with $\s>0$.  For $\s>1$ we can substitute the convergent series (4.1) into (2.2) and integrate
term by term using the integral formula
 $$ 2z\int_0^\infty\frac{t^{2s}\,dt}{(m^2t^2+n^2)^s\,(z^2+t^2)^{s+1}}=\frac{\G(s+\h)\G(\h)}{\G(s+1)}\,\frac1{(mz+n)^{2s}}$$
(which is just (2.38) with the omitted constant re-inserted) to get $\psi\e\psi_s^+$.  This calculation was also done
by Chang and Mayer [2].

\medbreak
2. Start with $\psi=\psi_s^+$ and try to work out the corresponding Maass form by using the correspondences $\psi\lr f\lr u$
from Chapter I.  The first step is easy: if $f(z)$ is the function defined by equation (0.6) with $\psi=\psi_s^+$, then
  $$\align f(z)&\,\e\,\psi_s^+(z)+z^{-2s}\,\psi_s^+(-1/z)\,=\,\sums_{m\ge0}\sum_{n\in\Z}\frac1{(mz+n)^{2s}}\tag4.3\\
  &\=\h\bigl(1+e^{-2\pi is}\bigr)\z(2s)+\frac{(-2\pi i)^{2s}}{\G(2s)}\,\sum_{n=1}^\infty\s_{2s-1}(n)\,e^{2\pi inz}\endalign$$
for $z\in\H$ (the second equality follows from the Lipschitz formula for\break $\sum(z+n)^{-\nu}$), and of course $f(-z)=-f(z)$
since we are in the even case. Comparing this with the Fourier expansion (4.2), we see that indeed the two functions $f$
and $u=E_s$ are related up to a constant factor by the recipe given in Chapter I (equations (1.9) and (1.11)), except that there we had no
constant terms and we did not specify how to get the coefficients of $y^s$ and $y^{1-s}$ for a noncuspidal $u$ from $f$.
Comparing (4.3) and (4.2) suggests that the former coefficient should come from the constant term of the Fourier expansion
of $f$, but it is not yet clear where the coefficient of $y^{1-s}$ comes from. (We will see the answer in a moment.)
\medbreak
3. According to the proposition at the end of Section~4 of Chapter II, the period function of an even Maass cusp form $u$ is $C^\infty$
from the right at 0 and its $m^{\rm th}$ Taylor coefficient at $0$ vanishes for $m$ even and is a simple multiple of $L_0(m+s+\h)$ for
$m$ odd, where $L_0(\rho)$ is the $L$-series of $u$. On the other hand, the function $\psi_s^+(x)$ has an expansion at $0$ given by
  $$ \psi_s^+(x)\sim\,\frac{\z(2s)}2\,x^{-2s}+\frac{\z(2s-1)}{2s-1}x^{-1}
    +\sum\Sb m\ge1\\m\text{ odd}\endSb \binom{m+2s-1}m\,\frac{B_{m+1}}{m+1}\,\z(m+2s)\,x^m\,, $$
where $B_n$ denotes the $n^{\rm th}$ Bernoulli number (this easily proved result was already mentioned under ``Example 2" at the
end of \S3 of Chapter III); and since by (4.2) the $L$-series of $E_s$ is a multiple of $\z(\rho-s+\h)\z(\rho+s-\h)$, we see
that the same relationship holds for the pair $u=E_s$, $\psi=\psi_s^+$. Moreover, from this point of view we can also see where
the two first coefficients in (4.2) come from: they are (up to simple multiples) the coefficients of $x^{-2s}$ and~$x\inv$ in
the asymptotic expansion of $\psi(x)$ at $x=0$.  This suggests the following theorem.

\nonumproclaim{Theorem} Let $s$ be a complex number with $\Re(s)>0$, $s\notin\Z${\rm . } 
\smallbreak
{\rm a)} If $u$ is a $\G$-invariant function in $\H$ with Fourier expansion
$$ u(z)=c_0y^s+c_1y^{1-s}+2\sqrt y\,\sum_{n=1}^\infty A_n\,K_{s-\h}(2\pi|n|y)\,\cos(2\pi nx)\tag 4.4$$
and we define a periodic holomorphic function $f:\C\sm\R\to\C$ by
$$ \pm f(z)\= \frac{\pi^{\h-s}}{\G(\h-s)}\,c_0 + \sum_{n=1}^\infty n^{s-\h}A_n\,e^{\pm2\pi inz} \qquad\bigl(\Im(z)\bs0\bigr)\,,\tag 4.5$$
then the solution $\psi$ of the three\/{\rm -}\/term functional equation $(0.1)$ defined by $(0.5)$ extends holomorphically to $\Cc$ and satisfies
  $$ \psi(x) \=\frac{\pi^\h\G(s+\h)}{\G(s)}\,\frac{c_0}{x^{2s}}\,+\,\frac{c_1}x\,+\,\O(1) \qquad(x\to0)\,.\tag4.6$$
\medbreak
{\rm b)} Conversely{\rm ,} if $\psi$ is a real\/{\rm -}\/analytic solution of $(0.1)$ on $\R_+$ with asymptotics of the form $(4.6)${\rm ,}
then $\psi$ extends holomorphically to $\Cc${\rm ,} the function $f$ defined by~$(0.6)$ has a Fourier expansion of the form $(4.5)${\rm ,}
and the function $u$ defined by $(4.4)$ is $\G$\/{\rm -}\/invariant.{\rm }\endproclaim
 
\demo{Proof} a) The function $u$ is the sum of a cusp form and a multiple of $E_s$.  The assertion is true for cusp forms by
the results of Chapter~I (with $c_0=c_1=0$) and for the Eisenstein series by the discussion above.
\medbreak
b) By the evenness of $\psi$, the assumed asymptotic behavior is equivalent to the assertion that $\psi(x)\sim cx^{1-2s}+c'+\O(x^{-2s})$
as $x\to\infty$ for some $c,\,c'\in\C$. If $c=c'=0$ then the asserted facts are the contents of Theorems~1 and~2. We indicate
how to modify the proofs of these theorems to apply to the new situation.  

To prove that $\psi$ extends holomorphically to $\Cm$ we follow the ``bootstrapping" proof of Chapter~III, the only change
being that equation (3.19) is replaced by
   $$ \psi(z)\,=\,c'\,+\,\psi(1)\,\z(2s,z+1)+\sum_{n=1}^\infty(n+z)^{-2s}\,\biggl(\psi\bigl(1-\frac1{n+z}\bigr)-\psi(1)\biggr)\,.\tag4.7$$
({\it Proof}\/: Note first that $\psi(1)=(2s-1)c$ by the three-term functional equation. Using the functional equation
we deduce that difference of the two sides of (4.7) is a periodic function, and since it is also $\,\o(1)\,$ at infinity it must vanish.
This argument is essentially the same as the proof of the proposition in \S3 of Chapter III.) At the same time we find that the
assumed asymptotics of $\psi(x)$ at 0 and $\infty$ remain true for $\psi(z)$ in the entire right half-plane, and that $\psi(z)$ is
bounded by a negative power of $|y|$ near the cut in the left half-plane. (To prove the latter statement we note that replacing
(3.19) by (4.7) gives (3.23) and (3.25) with the exponent $\s$ replaced by $-C$ for some $C>0$, and using $|y_{i+1}|\ge2|y_i|$
we obtain $F(z)=\O\bigl(|y|^{-C}\bigr)$.)

To  get (4.5), observe that the asymptotic expansions of $\psi(iy)$ at 0\break and $\infty$ imply that the function
$c^\star(s)f(iy)=\psi(iy)+(iy)^{-2s}\psi(i/y)$ equals $c'(1+e^{-2\pi is})+\o(1)$ as $y\to\infty$. Also, $f(z)$ is periodic by Proposition 2 of
Chapter I, Section~2. It follows that $f(z)$ has a Fourier expansion $\sum_{n\ge0}a_ne^{2\pi inz}$ in $\H$ with
$c^\star(s)a_0=c'(1+e^{-2\pi is})$. Now use the relation between $c'$ and $c_0$ in (4.6) together with formula (1.12).

\def\wf{\widetilde f}  
We now have the coefficients $A_n$ and can define $u$ by (4.4). This function is automatically an eigenfunction of $\D$ and
periodic, so we only need to show the invariance of $u(iy)$ under $y\mapsto1/y$. We will follow the $L$-series proof of Chapter~I
with suitable modifications. Note first that, by virtue of the estimate of $\psi$ near the cut given above, the coefficients
$A_n$ have at most polynomial growth, so that the $L$-series $L_0(\rho)=2\sum A_nn^{-\rho}$ converges in some half-plane. If we now define 
  $$u_0(y)=\frac1{\sqrt y}\,\bigl(u(iy)-c_0y^s-c_1y^{1-s}\bigr)=2\,\sum_{n=1}^\infty A_n\,K_{s-\h}(2\pi|n|y)$$
(cf\. (1.6)), then the Mellin transform of $u_0(y)$ is $L_0^*(\rho):=\g_s(\rho)L_0(\rho)$ for $\Re(\rho)$ sufficiently large.
Similarly the Mellin transform $\wf(\rho)$ of $f(iy)-a_0$ is a multiple of $L_0(\rho-s+\h)$, the formulas being the same
as in (1.16) (with $\wf_\pm=\pm\wf$ and $L_1=L_1^*=0$).  Now the same arguments as in Section~4 of Chapter I let us deduce from
the analytic continuability of $\psi$ across $\R_+$ that $L_0^*(\rho)=L_0^*(1-\rho)$ and from this that $u(iy)=u(i/y)$.
The effect of subtracting the powers of $y$ from $u(iy)$ and the constant term from $f(iy)$ in order to get convergence in a half-plane
is that the function $L_0^*(\rho)$ is now no longer entire, but acquires four simple poles, at $\h-s$, $\frac32-s$,
$s-\h$ and $s+\h$, with residues $-c_0$, $c_1$, $-c_1$ and $c_0$, respectively. Also, the Mellin transforms of $\psi(x)$ 
and $\psi(\pm iy)$ do not necessarily converge in any strip and we must subtract off a finite number of elementary functions 
from $\psi$ in order to define and compute these transforms. The details of the argument are left to the reader. \enddemo

\demo{Remarks} 1.  Part (b) of the theorem and the fact that cusp forms have codimension~1 in the space of all Maass forms imply that if
$\psi(x)$ is any analytic even periodlike function satisfying $\psi(x)\sim cx^{1-2s}+c'+\O(x^{-2s})$ as $x\to\infty$,
then $c=\lambda\,\dfrac{\z(2s-1)}{2s-1}$, $\;c'=\lambda\,\dfrac{\z(2s)}2$ for some $\lambda\in\C$. This fact will be used in Section~3.

\medbreak {2.} We stated the theorem only in the even case. The odd case is uninteresting.  On the one hand, since the function $E_s$ is
even, any odd Maass form is cuspidal.  In the other direction, if an odd periodlike function has an asymptotic expansion of the
form (4.6), then $c_1=\psi(1)/(2s-1)=0$ by the three-term equation, while $c_0$ can be eliminated by subtracting from $\psi$ a
multiple of the trivial periodlike function $\psi_s^-(z)=1-z^{-2s}$ of (3.1).
 
\section{Integral values of $s$ and connections\\ with holomorphic modular forms} 

In this section, we will first review the classical Eichler-Shimura-Manin theory
of period polynomials of holomorphic modular forms and describe the analogies 
between the properties of these polynomials and of the holomorphic function 
$\psi(z)$ associated to a Maass wave form, justifying the title of the paper.  
Second, and more interesting, we will show that the two theories are not only 
analogous, but are in fact related to one another in the special case when the 
spectral parameter $s$ is an integer.  In that case, there are no Maass cusp
forms, but there are ``nearly automorphic" eigenfunctions $u$ of the Laplace
operator whose associated periodic and holomorphic function $f(z)$  is a 
holomorphic cusp form of weight $2k$ in the case $s=k\in\Z_{>0}$ and the Eichler integral
of such a form if $s=1-k\in\Z_{\le0}$.

Let $f(z)$ be a holomorphic cusp form of weight $2k$ on $\G$; to fix notation we recall that this means that
$f$ satisfies $f\|{2k}\g=f$ for all $\g\in\G$, where the weight $2k$ action of $G$ on functions is defined by
$\bigl(F\|{2k}\sma abcd\bigr)(z)=\break (cz+d)^{-2k}F\bigl(\dfrac{az+b}{cz+d}\bigr)$, and that $f$ has a Fourier expansion of the form
   $$ f(z)=\sum_{n=1}^\infty a_n\,q^n\qquad(z\in\H)\,,\tag 4.8$$
where we have used the standard convention $q=e^{2\pi iz}$.  Associated to $f$ is a polynomial $r_f$ of degree $2k-2$,
the {\it period polynomial} of $f$.  It can be defined in three ways (we now use the symbol $\e$ to denote 
equality up to a constant depending only on~$k$):
\smallbreak
(i) by the identity 
    $$ r_f(z) \,\e\, \wf(z) - z^{2k-2}\wf(-1/z)\qquad(z\in\H)\,,\tag 4.9$$
where $\wf$ is the {\it Eichler integral} of $f$, defined by the Fourier expansion 
  $$ \wf(z)=\sum_{n=1}^\infty \frac{a_n}{n^{2k-1}}\,q^n\qquad(z\in\H)\,;\tag 4.10$$

\smallbreak
(ii) by the integral representation
  $$ r_f(X) \,\e\, \int_0^\infty f(\tau)\,(\tau-X)^{2k-2}\,d\tau\,,\tag 4.11$$
where the integral is taken over the positive imaginary axis;

\smallbreak
(iii) by the closed formula
  $$ r_f(X) \,\e\, \sum_{r=0}^{2k-2}\frac{(-2\pi i)^{-r}}{(2k-2-r)!}\,L_f(r+1)\,X^r\,, \tag 4.12$$
where $L_f(\rho)=\sum_{n=1}^\infty a_nn^{-\rho}$ (or its analytic continuation) is the Hecke $L$-series associated to $f$.

The proofs that these definitions agree are simple. Denote by $D$ the differential operator 
  $$  D \= \frac1{2\pi i}\,\frac d{dz} \= q\,\frac d{dq}\,, $$
so that the relationship between the functions with the Fourier expansions (4.8) and (4.10) is given by
  $$   D^{2k-1}(\wf) \= f \,. \tag 4.13 $$
The key point is that the $(2k-1)$st power of the operator $D$ {\it intertwines} the actions of $G={\rm PSL}(2,\R)$
in weights $2-2k$ and $2k$; i.e.
  $$ D^{2k-1}\bigl(F\|{2-2k}g\bigr) \= \bigl(D^{2k-1}F\bigr)\|{2k}g  \tag 4.14$$
for any $g\in G$ and any differentiable function $F$.  (The identity (4.14) is known as Bol's identity; we will see
why it is true a little later.) It follows that the function $r_f$ defined by (4.9) satisfies
  $$ \align D^{2k-1}\bigl(r_f\bigr) &\=  D^{2k-1}\bigl(\wf-\wf\|{2-2k}S\bigr) \\
\noalign{\vskip5pt}
    &\=  D^{2k-1}\wf-\bigl(D^{2k-1}\wf\bigr)\|{2k}S \= f-f\|{2k}S \= 0 \endalign $$
(where $S=\sma0{-1}10$ as usual) and hence is indeed a polynomial of degree $2k-2$. To show that it is proportional to the polynomial
defined by (4.11), we observe that the $(2k-1)$\-fold primitive $\wf$ of $f$ can also be represented by the integral
$\wf(z)\e\int_z^\infty (z-\tau)^{2k-2}\,f(\tau)\,d\tau\;$ ({\it Proof}: the right-hand side is exponentially small at infinity
and its $(2k-1)$st derivative is a multiple of $f$), and from this and the modularity of $f$ we get 
$z^{-2k}\wf(-1/z)\e\int_z^0(z-\tau)^{2k-2}\,f(\tau)\,d\tau$, from which the asserted equality follows.
Finally, the equality of the right-hand sides of (4.11) and (4.12) follows from the representation of $L_f(r+1)$ as a multiple 
of $\int_0^\infty \tau^r f(\tau)\,d\tau$.

It is now clear why throughout this paper we have been referring to the function $\psi(z)$ associated to a Maass
form $u$ as its ``period function," for each of the defining properties (i)--(iii) has its exact analogue in the
theory we have been building up.  Formula (4.9) is the analogue of formula (0.5) expressing $\psi$ as $f|(1-S)$
where $f$ is the periodic holomorphic function attached to $u$; formula (4.11) corresponds to the expression (2.8) given
in \S2 of Chapter II for $\psi$ as an integral of a certain closed form attached to $u$; and formula (4.12) is the
analogue of the result given in Chapter II (eq.~(2.33)) for the Taylor coefficients of $\psi$ as multiples of
the values of the $L$-series of $u$ at (shifted) integer arguments.

We can in fact make the analogy even more precise. Denote by $P_{2k-2}$ the space of polynomials of degree $\le2k-2$,
with the action $\|{2-2k}$ (we will drop the subscript from now on) of $G$.  (Note that $P_{2k-2}$ with this action 
is a sub-representation of the space $\text{\script V}_{1-k}$ defined at the end of Section~5 of Chapter II.) It is easily shown using
the above definitions that the period polynomial $r_f$ of a cusp form belongs to the space
  $$ W_{2k-2}\,:=\,\bigl\{F\in P_{2k-2}\,:\, F|(1+S)=F|(1+U+U^2)=0\bigr\}\,. $$
Here $S$ and $U=TS$ are the standard generators of $\G$ with $S^2=U^3=1$ and we have extended the action of $\G$ 
on $P_{2k-2}$ to an action of the group ring $\Z[\G]$ by linearity, so that e.g\. $F|(1+U+U^2)$ means $F+F|U+F|U^2$.
In fact it is known that the functions $r_f(X)$ and $\overline{r_f(\bar X)}$ as $f$ ranges over all cusp forms of weight $2k$
span a codimension~1 subspace of $W_{2k-2}$.  (The missing one-dimensional space comes from an Eisenstein series; see
below.)  The relations defining $W_{2k-2}$ express the fact that there is a 1-cocycle $\G\to P_{2k-2}$ sending the
generators $T$ and $S$ of $\G$ to 0 and $F$, respectively, this cocycle in the case $F=r_f$ being the map defined by
$\g\mapsto\wf|(1-\g)$. We now have:

\vglue8pt

\nonumproclaim{Proposition} Every element of $W_{2k-2}$ is a solution of the three\/{\rm -}\/term functional
 equation $(0.2)$ with parameter $s=1-k${\rm ,} 
and conversely if $k>1$ every periodlike function of parameter $1-k$ which is a polynomial is an element of $W_{2k-2}${\rm .} \endproclaim

\vglue8pt

\demo{Proof} The elements $T=\sma1101$ and $T'=\sma1011$ are represented in terms of the generators $S$ and $U$ by $T=US$, $T'=U^2S$,
respectively, so for $F\in W_{2k-2}$ we find $F|(1-T-T')= F|(1+S)-F|(1+U+U^2)|S=0$, which is precisely the three-term functional
equation. For the converse direction, we reverse\break the calculation to find that a periodlike function $F$ satisfies
$F|(1+S)=\break
F|(1+U+U^2)$. Writing $H$ for the common value of $F|(1+S)$ and $F|(1+U+U^2)$, we find that $H$ is invariant under both $S$ and
$U$ and hence under all of $\G$. 
It follows that $H=0$ (the only $T$-invariant polynomials are constants, and these are not $S$-invariant for $k>1$). Hence $F$ 
satisfies the equations defining $W_{2k-2}$, while from $F|S=-F$ it follows that $\,\text{deg}(f)\le 2k-2$, so $F\in P_{2k-2}$.
 \enddemo

\vglue8pt

The proposition and the preceding discussion show that for the parameter $s=1-k$ the period polynomials of holomorphic cusp
forms of weight $2k$ produce holomorphic solutions of the three-term functional equation with reasonable growth properties 
at infinity. How does this fit into our Maass picture?  The answer is very simple.  For each integer $h$ the differential operator 
$\p_h := D-ih/(2\pi y)$ (where $z=x+iy$ as usual) intertwines the actions of $G$ in weights $2h$ and $2h+2$, i.e.
  $$ \p_h\bigl(F\|{2h}g\bigr) \= \p_h(F)\|{2h+2}g \,. $$
In particular, if $F$ is modular (or nearly modular) of weight $2h$, then $\p_hF$ is modular (or nearly modular) of weight
$2h+2$.  Iterating, we find that the composition $\p_h^n:=\p_{h+n-1}\circ\dots\circ\p_h$ intertwines the actions
of $G$ in weights $2h$ and $2h+2n$ and in particular sends modular or nearly modular forms of weight $2h$ to modular or 
nearly modular forms of weight $2h+2n$.  On the other hand, by induction on $n$ one proves the formula
  $$ \p_h^n \= \sum_{m=0}^n\frac{n!}{(n-m)!}\,\binom{n+2h-1}m\,\biggl(\frac{-1}{4\pi y}\biggr)^m\,D^{n-m}\,.\tag 4.15$$
In the special case when $h=1-k$ and $n=2k-1$ this reduces simply to $\p_h^n=D^n$, which explains Bol's identity (4.14):
the operator $\p_h^n$ always preserves modularity (shifting the weight by $2n$) but in general destroys holomorphy, while $D^n$ 
preserves holomorphy but in general destroys modularity, so that in the case  $n+2h=1$ when they agree both properties are preserved. 

Now, dropping the lower indices on the $\p$'s for convenience, we can factor the identity $D^{2k-1}=\p^{2k-1}$ as
$D^{2k-1}=\p^k\circ\p^{k-1}$, leading to the picture indicated by the following diagram:
\medbreak
 \centerline{\BoxedEPSF{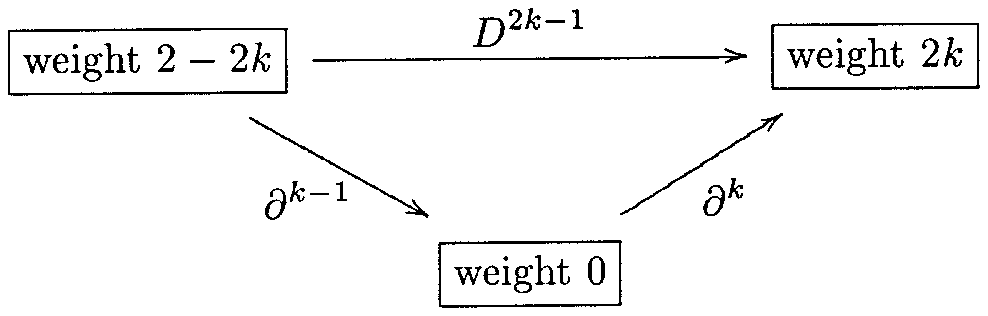 scaled 1000}}
\medbreak
\noindent 
Hence we can factor the relation (4.13) into steps, writing
 $$ f\=\p_0^k(u)\qquad\text{where}\qquad u:=\p_{1-k}^{k-1}(\wf)\,. \tag 4.16$$
From (4.15) and (4.10) we immediately obtain the Fourier expansion of $u$, and comparing the result with the well-known formulas for
Bessel functions of half-integer index in terms of elementary functions, we find to our delight the familiar-looking expression
  $$ u(z) \,\e\, \sqrt y\,\sum_{n=1}^\infty A_n\,K_{k-1/2}(2\pi ny)\,e^{2\pi i nx}\,, \tag 4.17$$
with Fourier coefficients $A_n$ given by
  $$ A_n\= \cases n^{-k+1/2}\,a_n&\text{for $n>0\,,$}\\0&\text{for $n<0\,.$}\endcases\tag 4.18$$
In particular, the function $u$ defined in (4.16) is an eigenfunction of the Laplace operator with eigenvalue $k(1-k)$
and is $T$-invariant and small at infinity.  It is not quite $\G$-invariant, but its behavior under the action of the
second generator $S$ of $\G$ is easily determined:
  $$ \align u(z)\,-\,u(-1/z) &\= u\|0(1-S) \\ \noalign{\vskip4pt}  &\e \bigl(\p^{k-1}_{1-k}\wf\bigr)\|0(1-S) \\
\noalign{\vskip4pt}
   &\= \p^{k-1}_{1-k}\bigl(\wf\|{2-2k}(1-S)\bigr) \\
\noalign{\vskip4pt} &\=  \p^{k-1}_{1-k}\bigl(r_f\bigr) \endalign$$
by (4.9) and the intertwining property of $\p_h^n$.   Thus $u(z)-u(-1/z)$, and hence also $u(z)-u(\g(z))$ for every $\g\in\G$,
belongs to the to the space $P_k^0$ of polynomials in $x$, $y$ and $1/y$ which are annihilated by $\D-k(1-k)$.
In fact, the map $\p^{k-1}_{1-k}$ is an isomorphism from $P_{2k-2}$ to $P_k^0$ and transfers the original
cocycle $\g\mapsto \wf\|{2-2k}(1-\g)$ with coefficients in $P_{2k-2}$ to the cocycle $\g\mapsto u|_0(1-\g)$ with
coefficients in $P_k^0$. 

Let us compare this new correspondence with our usual correspondence between Maass cusp forms and their period functions.
There are several points of difference as well as of similarity.  The obvious one is that in our new situation the
eigenfunction $u$ is no longer a Maass form, but only ``Maasslike" in the sense just explained.  But there are other
differences.  First we point out a property of the usual $u\lr\psi$ correspondence  which we have not previously
emphasized: associated to a Maass form $u$ there is not just one, but {\it two} period functions.  Namely, if $u$ is an
eigenfunction of $\D$ with eigenvalue $\lambda$ and
we write its Fourier expansion in the form (1.9), we can freely choose between the parameter $s$ and the
parameter $1-s$, since the $K$-Bessel function $K_\nu(t)$ is an even function of its index $\nu$.  But when we
write down the associated periodic function $f$ in $\C\sm\R$ by (1.11), we have broken the $s\lr1-s$ symmetry and chosen
one of the two roots of $s(1-s)=\lambda$, so that there is actually a {\it second} holomorphic periodic function $\wf$,
 defined by the same formula (1.11) but with the exponent $s-\h$ replaced by $\h-s$ (compare eq.~(2.44), where there
were two boundary forms $U(t)$ and $\wU(t)$ associated to $u$), and similarly a second period function 
$\widetilde\psi=\wf\|{2-2s}(1-S)$. The choice was not important in the case of a Maass
cusp form, since then $\Re(s)=\h$ anyway, so that $\psi(z)\mapsto\overline{\psi(\bar z)}$ gives a correspondence between
the two possible choices for the period function.  But in our new situation the picture is different:
if $u$ is the Maasslike function defined by (4.17) and (4.18), then with the spectral parameter $s=k$ we see
that the ``$f$" defined by (1.11) is our original cusp form $f$ and the ``$\wf$" defined by (1.11) with $s$
replaced by $1-s$ is the associated Eichler integral (4.10). The reason for this dichotomy is that the 
correspondence between periodlike and periodic functions described in Proposition~2, Section~2, of Chapter~I 
breaks down when the parameter $s$ is an integer: the maps $\psi\mapsto\psi\|{2s}(1+S)$ and $f\mapsto f|_{2s}(1-S)$
still send periodlike to periodic functions and vice versa, but are no longer isomorphisms, since their
composition is~0 in both directions. (This is because $h\mapsto h\|{2s}g$ is a $G$-action for $s\in\Z$,
and $S^2=1$ in $G$.)  The two functions $f$ and $\wf$ associated to an eigenfunction $u$ with eigenvalue
$k(1-k)$ behave very differently under these correspondences: $f$ itself (if we make the choice $s=k$ and
if $u$ is the Maasslike form associated to a holomorphic cusp form) is $\G$-invariant in weight $2k$ and
hence is annihilated by $1-S$, while $\wf$ is mapped by $1-S$ to the period polynomial $r_f$ which in turn, unlike
the period functions of true Maass forms, is then in the kernel of $1+S$. The situation is summarized by
the following diagram:

 \vglue-12pt
\centerline{\BoxedEPSF{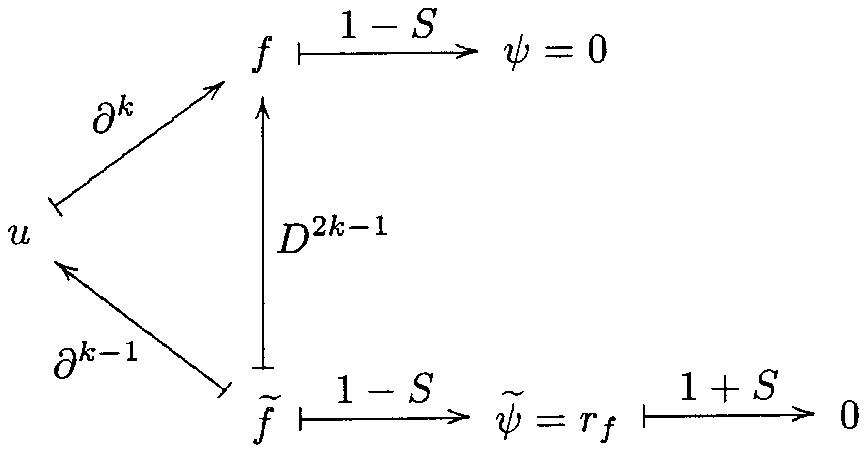 scaled 900}}

We mention one other point of difference: in the usual case a Maass form was typically even or odd (it is believed---and
certainly true for all eigenfunctions computed up to now---that the space of  Maass forms for a given eigenvalue
is at most one-dimensional, in which case every Maass form is necessarily an even or an odd one), while the
expansion (4.17) does not have either of the symmetry properties $A_{-n}=\pm A_n$.  If we split (4.17) into its
even and odd parts by writing $e^{2\pi inx}$ as the sum of $\cos(2\pi nx)$ and $i\,\sin(2\pi nx)$, we find two different 
cocycles in $H^1(\G,P_k^0)$, given by the even and odd parts of the period polynomial $r_f$.  This corresponds to the fact
mentioned above that the space $W_{2k-2}$ contains two isomorphic images of the space of holomorphic cusp forms of weight $2k$.


To complete the picture we should say something about the Eisenstein case. The space of all modular forms of weight $2k$
on $\G$ is spanned (for $k>1$) by the cusp forms and the classical Eisenstein series $G_{2k}$.  For $f=G_{2k}$
one can still define an associated function $r_f$ by appropriate modifications of each of
the three definitions (i)--(iii), as was shown in [15]. The result (proposition on p.~453 of [15]) is 
  $$ r_{G_{2k}}(X) \;\e\;\sum_{n=0}^k\frac{B_{2n}}{(2n)!}\,\frac{B_{2k-2n}}{(2k-2n)!}\,X^{2n-1}
   + \frac{\z(2k-1)}{(2\pi i)^{2k-1}}\bigl(X^{2k-2}-1\bigr)\,,$$
which is not quite an element of $P_{2k-2}$ but of the larger space (on which $G$ does not act) spanned by $\{X^n\mid-1\le n\le2k-1\}$.
The even and odd parts of $r_{G_{2k}}$ are (multiples of) the functions $\psi_{1-k}^+$ and $\psi_{1-k}^-$ discussed in Examples~1 and~2
of Chapter III, Section~1. This is in accordance with the correspondence for cusp forms described above, since from (4.2) and the fact that 
the $n^{\rm th}$ Fourier coefficient of $G_{2k}$ is $\sigma_{2k-1}(n)$ for $n>0$ we see that the ``$u$" associated to $f=G_{2k}$
should be precisely the nonholomorphic Eisenstein series $E_k$, whose period function for the choice of spectral parameter
$s=1-k$ is indeed the function $\psi_{1-k}^+$, as we saw in Section~1. This function belongs to $P_{2k-2}$ and spans the ``missing 
one-dimensional space" mentioned before the proposition above. The other period function $\psi_k^+$ of $E_k$, associated
to the choice of spectral parameter $s=k$, is not a polynomial but instead is the periodlike function whose image under
$1+S$ is the element $f=G_{2k}$ of $\,\ker(1-S)$.

 We make one final observation. All three definitions (i)--(iii) of classical period polynomials have analogies for the
periods of Maass forms, as discussed at the beginning of this section, but there is one aspect of the classical theory which does not
immediately generalize, namely the interpretation of the map $T\mapsto 0$, $S\mapsto r_f$ as a cocycle. The reason is
that for $s$ nonintegral the ``action" $\,\|{2s}\,$ is not in fact an action of the group $G$, because the automorphy
factors $(cz+d)^{-2s}$ have an ambiguity given by powers of $e^{2\pi is}$. This means that the relation $\psi\|{2s}(1-T-T')=0$
(three-term functional equation) does {\it not} imply that there is a cocycle on $\G$ (say, with coefficients in the space
of holomorphic functions on $\Cc$) sending $T$ to 0 and $S$ to~$\psi$.  Nevertheless, there are ways to interpret
$\psi$ as part of a 1-cocycle, and these in fact work for any discrete subgroup of $G$, permitting us to extend the theory
developed in this paper to groups other than ${\rm PSL}(2,\Z)$.  This will be the main subject of Part~II of this paper.

\section{Relation to the Selberg zeta function and Mayer's theorem}

The theme of this paper has been the correspondence between the spectral parameters of the group $\G={\rm PSL}(2,\Z)$ 
and the holomorphic solutions of the three-term functional equation (0.2).  On the other hand, there is a
famous relation between the same spectral parameters and the  set of lengths of the closed geodesics on the
Riemann surface $X=\H/\G$, namely the Selberg trace formula, which expresses the sum of the values of suitable
test functions evaluated at the spectral parameters as the sum of a transformed function evaluated on this
length spectrum.  The triangle is completed by a beautiful result of Mayer, which relates the length spectrum of~$X$,
as encoded by the associated Selberg zeta function $Z(s)$, to the eigenvalues of a certain linear operator 
$\L_s$ which  is closely connected with the three-term functional equation.
Combining this theorem, the underlying idea of whose proof is essentially elementary, with the theory 
developed in this paper yields a direct connection between Maass wave forms and the length spectrum of $X$,
and hence a new insight into the Selberg trace formula.

We begin by recalling the definitions of $Z(s)$ and $\L_s$ and the proof of Mayer's theorem, referring 
to [13] and [11] for more details. The function $Z(s)$ is defined for $\Re(s)>1$ by the product expansion 
$$Z(s)=\prod\Sb\{\g\}\text{ in }\G\\ \g\text{ primitive}\endSb \prod_{m=0}^\infty\bigl(1-\N(\g)^{-s-m}\bigr)\,.\tag 4.19$$
Here the first product is taken over $\G$-conjugacy classes of primitive hyperbolic elements of $\G$ (``hyperbolic"
means that the absolute value of the trace is bigger than 2, and ``primitive" that $\g$ is not a power of any
matrix of smaller trace), and the norm $\N(\g)$ is defined as $\frac14\bigl(|\text{tr}(\g)|+\sqrt{\text{tr}(\g)^2-4}\,\bigr)^2$
or equivalently as $\ve^2$ where $\g$ is conjugate in ${\rm PSL}(2,\R)$ to $\pm\sma\ve00{1/\ve\!}$.  The function $Z(s)$, or
rather its logarithmic derivative, arises by applying the Selberg trace formula to a particular family of test
functions parametrized by the complex number~$s$.  The trace formula then implies that $Z(s)$ extends
meromorphically to all $s$, with poles at negative half-integers and with zeros at the spectral parameters of $\G$,
together with the value $s=1$ and the zeros of $\z(2s)$.

The operator $\L_s$ is an endomorphism of the vector space~$\V$ of functions which are holomorphic in the disk 
$\Dk=\{z\in\C\mid |z-1|<\frac32\}$ and continuous in $\bD$. It is defined for $\Re(s)>\h$ by 
  $$(\L_sh)(z)=\sum_{n=1}^\infty\frac1{(z+n)^{2s}}\,h\bigl(\frac1{z+n}\bigr)\qquad(h\in\V)\,,\tag 4.20$$
where the holomorphy of $h$ at 0 implies that the sum converges absolutely and again belongs to $\V$.
We continue this meromorphically to all complex values of $s$ by setting
  $$(\L_sh)(z)= \sum_{m=0}^{M-1}c_m\,\z(2s+m,z+1)+(\Cal L_sh_0)(z)\,,\tag 4.21$$
where $\z(s,z)$ denotes the Hurwitz zeta function, $M$ is any integer greater than $1-2\Re(s)$, the $c_m$
($0\le m\le M-1$) are the first $M$ Taylor coefficients of $h(z)$ at 0, and 
$h_0(z)=h(z)-\sum_{m=0}^{M-1}c_mz^m$. This is clearly independent of $M$ and holomorphic except for simple poles
at $2s=1,\,0,\,-1,\ldots\;$. Mayer proves that the operator $\L_s$ is of trace class (and in fact
nuclear of order~0), from which it follows that the operators $1\pm\L_s$ have determinants in the Fredholm sense.   

\nonumproclaim{Theorem {\rm (Mayer [13], [14])}} \hskip-4pt The Selberg zeta function of $\H/\G$ is given~by  
$$ Z(s)=\det\bigl(1-\L_s\bigr)\det\bigl(1+\L_s\bigr)\,.\tag 4.22$$ \endproclaim

A simplified version of the proof is given in [11].  Roughly, the idea is as follows.  We may assume
$\Re(s)>1$. After an elementary manipulation, (4.19) can be rewritten
$$  \log Z(s)\=-\sum_{\{\g\},\,k}\frac1k\,\chi_s\bigl(\g^k\bigr)\, , $$
where the sum over $\g$ is the same as before and $k$ runs over all integers $\ge1$, and where
$\chi_s(\g)=\N(\g)^{-s}/\bigl(1-\N(\g)^{-1}\bigr)\,$.  By the reduction theory of quadratic forms, every
conjugacy class $\{\g^k\}$ of hyperbolic matrices in $\G$ has a finite number of ``reduced" representatives
(a matrix $\sma abcd\in\G$ is called reduced if $0\le a\le b,\,c\le d$); these have the form
$\rho_{n_1}\cdots\rho_{n_{2l}}$ where $n_1,\ldots,n_{2l}\in\Bbb N$ are the convergents in a periodic continued
fraction expansion associated to $\{\g^k\}$ and $\rho_n=\sma011n$, and the $l/k$ reduced representatives of the 
conjugacy class are all obtained by permuting $(n_1,\ldots,n_{2l})$ cyclically.  On the other hand, any reduced
matrix $\g=\sma abcd$ in $\G$ acts on $\V$ by $(\pi_s(\g)f)(z)=(cz+d)^{-2s}f\bigl(\frac{az+b}{cz+d}\bigr)$
(this makes sense because $\g(\Dk)\subset\Dk$ for $\g$ reduced), and this
operator is of trace class with Fredholm trace $\text{Tr}(\pi_s(\g))=\chi_s(\g)$.  Putting all of
this together and observing that $\L_s=\sum_{n=1}^\infty\pi_s(\rho_n)$, we find
$$ \log Z(s)\=-\sum_{l=1}^\infty\frac1l\,\text{Tr}\bigl(\L_s^{2l}\bigr) \= \log\det\bigl(1-\L_s^2\bigr)$$
as claimed. Actually, the formula (4.22) can be made a little more precise, as discussed in [5] and [11]: the function
$Z(s)$ has a natural splitting $Z_+(s)Z_-(s)$ where $Z_+(s)$ is the Selberg zeta function of $\G^+={\rm PGL}(2,\Z)$
(with the elements of $\G^+$ of determinant $-1$ acting on $\H$ by $z\mapsto(a\bar z+b)/(c\bar z+d)\,$) and 
where $Z_+(s)$ and $Z_-(s)$ have zeros at the spectral parameters of $\G$ corresponding to even and odd
Maass forms, respectively, and one in fact has $Z_\pm(s)=\det\bigl(1\mp\L_s)\,$.

We now come to the connection with period functions. Fredholm theory implies that the trace class
operators $\L_s$ share with operators of finite rank the property that $\det(1\mp\L_s)=0$ if and only 
if $\L_s$ has an eigenvector with eigenvalue $\pm1$.  Hence combining Mayer's theorem (in the sharpened
form just mentioned) with the known position of the zeros of $Z(s)$ implied by the Selberg trace formula, we obtain
\nonumproclaim{{C}orollary} Let $s\ne\h$ be a complex number with $\Re(s)>0${\rm .}  Then\/{\rm :}\/
\smallbreak
{\rm a)} there exists a nonzero function $h\in\V$ with $\L_sh=-h$ if and only if $s$ is the spectral parameter
corresponding to an odd Maass wave form on $\G$\/{\rm ;}\/

\smallbreak
{\rm b)} there exists a nonzero function $h\in\V$ with $\L_sh=h$ if and only if $s$ is either the spectral parameter 
corresponding to an even Maass form{\rm ,} or $2s$ is a zero of the Riemann zeta function{\rm ,} or $s=1${\rm .} \endproclaim

We now show how the main theorems of this paper give a constructive proof of this corollary, independent of Mayer's 
theorem and the Selberg trace formula. To do this, we use the following bijection between solutions of $\L_sh=\pm h$ and holomorphic 
periodlike functions.

\nonumproclaim{Proposition} Suppose that $\Re(s)>0$, $s\ne\h${\rm .} Then a function $h\in\V$
is a solution of $\L_sh=\pm h$ if and only if $h(z)$ is the restriction to $\Dk$ of $\psi(z+1)$ where
$\psi$ is a holomorphic solution in $\Cm$ of the even\/{\rm /}\/odd three\/{\rm -}\/term functional equation $(1.13)$ having the asymptotic
behavior $\psi(x)=cx^{1-2s}+\O\bigl(x^{-2s}\bigr)$ as $x\to\infty${\rm .} \endproclaim

\demo{Proof} First assume that $h\in\V$ is given with $\L_sh=\pm h$. Then the same bootstrapping
arguments as given in Section~4 of Chapter III (compare for instance equation (3.19))
 shows that the function $\psi(z):=h(z-1)$ extends analytically to all of $\Cm$. (The function $\L_sh$
for any $h\in\V$ is holomorphic in a much larger domain than $\Dk$, including in particular the
right half-plane $\Re(z)>-\frac35$, so $h=\pm\L_sh$ is also defined in this domain, and iterating
this argument we extend $h$ to a larger and larger region finally filling up the cut plane $\C\sm(-\infty,-1]$.)
The even or odd three-term functional equation for $\psi$ is obvious from the identity
$\L_sh(z-1)=\L_sh(z)+z^{-2s}h(1/z)$, which is valid for any $h\in\V$. The asymptotic expansion of $\psi$ at
infinity, with $c=h(0)/(2s-1)$, follows easily from equation (4.21).
 
Conversely, assume we are given a holomorphic function  $\psi$ in $\Cm$ which satisfies the even or
odd three-term functional equation and the given asymptotic formula. We must show that the function $h(z):=\psi(z+1)$,
which obviously belongs to $\V$, is a fixed-point of $\pm\L_s$.  The first observation is that, as pointed out in the
first part of the proof, the function $\L_sh$ is defined in the right half-plane and satisfies
$\L_sh(z-1)=\L_sh(z)+z^{-2s}h(1/z)$. From this and the functional equation of $\psi$ it follows that the function
$h_1(z):=\L_sh(z)\mp h(z)$ is periodic. But by letting $x\to\infty$ in (0.2) we find that the constant $c$ in the assumed
asymptotic formula for $\psi$ is given by $c=\psi(1)/(2s-1)$ (and in particular vanishes in the odd case); and from
this and equation (4.21) we find (since the leading terms cancel) that $h_1(x)$ is $\O(x^{-2s})$, and hence $\o(1)$, 
as $x\to\infty$, which together with the periodicity implies that $h_1\equiv0$. (Notice that this proof duplicates part 
of the proof of the theorem of Section~1: the assumptions on $\psi$ imply the hypotheses in part (b) of the theorem
with $c_0=0$, and the assertion $\L_sh=h$ is equivalent to equation (4.7) with $c'=0$.) \enddemo 

\vglue10pt

We can now write down the following explicit functions $h$ satisfying the conditions of the corollary: \medbreak
\item{(a)} If $s$ is the spectral parameter of an even or odd Maass form $u$ on $\G$, then the period function
associated to $u$ satisfies the hypotheses of the proposition with $c=0$ and hence gives a solution of $\L_sh=\pm h$ with $h(0)=0$. 
\medbreak
\item{(b)} If $\z(2s)=0$, then the function $\psi_s^+(z)$ studied in Section~1 of this chapter satisfies the conditions 
of the proposition, so $\psi_s^+(z+1)$ is a solution of $\L_sh=h$.
\medbreak
\item{(c)} If $s=1$, then the function $\psi(z)=1/z$ satisifes the hypotheses of the proposition, so the function $h(z)=1/(z+1)$
is a solution of $\L_1h=h$. (This can of course also be checked directly.) \medbreak\noindent
Conversely, any solution of $\L_sh=\pm h$ must be one of the functions on this list.  Indeed, in the odd case the
coefficient $c$ in the proposition is always 0, so the function $\psi(z)=h(z-1)$ is the period functions of an odd Maass cusp form
by the corollary to Theorem 2 (Chapter III).  The same applies in the even case if the coefficient $c$ vanishes.  If it doesn't,
and if $s$ is not an integer, then Remark~1 following the theorem of Section~1 of this chapter shows that $\z(2s)$ must be 0;
then $\psi$ must be a nonzero multiple of $\psi_s^+(z)$ because $s$ has real part less than 1/2 and hence cannot 
be the spectral parameter of a Maass cusp form. The case $s=1$ works like the cases when $\z(2s)=0$ (and in fact
can be absorbed into it if we notice that the family of
Maass forms $(s-1)E_s(z)$ is continuous at $s=1$ with limiting value a constant function and that the family of period functions
$(s-1)\psi_s^+(z)$ is continuous there with limiting value a multiple of $1/z$). Again this is the only possible solution for
this value of $s$ since by subtracting a multiple of $1/z$ from any solution we would get a period function corresponding
to a Maass cusp form, and they do not exist for this eigenvalue.  The case when $s$ is an integer greater than~1 was excluded
in the theorem in Section~1 because the $f\lr\psi$ bijection breaks down, but can be treated fairly easily by hand and turns out to be
uninteresting: there are no solutions of the three-term functional equation of the form demanded by the proposition, in
accordance with the corollary to Mayer's theorem. Finally,
we remark that our analysis could be extended to $\s<0$, but we omitted this to avoid further case distinctions
and because it turns out that the only solutions of the three-term functional equation with growth
conditions of the required sort are the ones for $s=1-k$ coming from holomorphic modular forms which were discussed
in the last section.  A complete analysis of the Mayer operator $\L_s$ in this case is given in [3].

\references
[1]
\name{R.\ Bruggeman}, Automorphic forms, hyperfunction cohomology, and
period functions, {\it J.\ Reine Angew.\ Math\/}.\ {\bf 492} (1997),
1--39.

[2] 
\name{C.\ Chang} and \name{D.\ Mayer}, The period function of the nonholomorphic
Eisenstein series for $\text{PSL}(2,{\bold Z})$, {\it Math.\ Phys.\
Electron.\ J\/}.\ {\bf 4} (1998), Paper 6, 8pp.

[3] 
\bibline, The transfer operator approach to Selberg's
zeta function and modular and Maass wave forms for $\text{PSL}(2,{\bold
Z})$,
in {\it Emerging Applications of Number Theory\/}, 73--141, 
{\it IMA Vol.\ Math.\
Appl\/}.\ {\bf 109}, Springer-Verlag, New York, 1999.

[4]
Y\name{.-J.\ Choie} and \name{D.\ Zagier}, Rational period functions for
$\text{PSL}(2,{\bold Z})$, in {\it A Tribute to Emil Grosswald\/}:
{\it Number Theory and Related Analysis\/}, {\it Contemp.\ Math\/}.\
{\bf 143}, 89--108, A.M.S., Providence, RI, 1993.

[5]
\name{I.\ Efrat}, Dynamics of the continued fraction map and the spectral
theory of $\text{SL}(2,{\bold Z})$, {\it Invent.\ Math\/}.\ {\bf 114}
(1993), 207--218.

[6]
A.\ Erd\'elyi\name{, et.\ al.}, {\it Tables of Integral Transforms\/}
(Bateman manuscript project), Vols.\ I, II, McGraw-Hill, New York,
1954.

[7]
\name{W.\ Ferrar},
Summation formulae and their relation to Dirichlet's series (II),
{\it Compositio Math\/}.\ {\bf 4} (1937), 394--405.

[8]
\name{S.\ Helgason}, {\it Topics in Harmonic Analysis on Homogeneous Spaces\/},
{\it Progr.\ in Math\/}.\ {\bf 13}, Birkh\"auser, Boston,  1981.

[9]
\name{J.\ B.\ Lewis}, Eigenfunctions on symmetric spaces with
distribution-valued
boundary forms, {\it J.\ Funct.\ Anal\/}.\ {\bf 29} (1978), 287--307.

[10]
\bibline, Spaces of holomorphic functions equivalent to
the even Maass cusp forms, {\it Invent.\ Math\/}.\ {\bf 127} (1997),
271--306.

[11]
\name{J.\ B.\ Lewis} and \name{D.\ Zagier}, Period functions and the Selberg
zeta function for the modular group, in {\it The Mathematical Beauty
of Physics\/}, 83--97, {\it Adv.\ Series in Math.\ Physics\/} {\bf 24},
World Sci.\ Publ., River Edge, NJ, 1997.

[12]
\name{H.\ Maa\ss}, \"Uber eine neue Art von nichtanalytischen automorphen
Funktionen und die Bestimmung Dirichletscher Reihen durch
Funktionalgleichungen, {\it Math.\ Ann\/}.\ {\bf 121} (1949),
141--183.

[13]
\name{D.\ Mayer}, The thermodynamic formalism approach to Selberg's
zeta function for\break $\text{PSL}(2,{\bold Z})$, {\it Bull.\ A.M.S\/}.\
{\bf 25} (1991), 55--60.

[14]
\bibline, Continued fractions and related transformations,
in {\it Ergodic Theory, Symbolic Dynamics and Hyperbolic Spaces\/}
(T.\ Bedford, M.\ Keane, C. Series, eds.), Oxford Univ.\ Press,
New York, 1991, 175--222.

[15]
\name{D.\ Zagier}, Periods of modular forms and Jacobi theta functions,
{\it Invent.\ Math\/}.\ {\bf 104} (1991), 449--465.

\endreferences
\bye

[1]  \name{R. Bruggeman}, 
Automorphic forms, hyperfunction cohomology, and period functions, {\it J. Reine Angew.\ Math.\/} {\bf 492} (1997), 1--39. 

[2]  \name{C. Chang} and \name{D. Mayer}, 
The period function of the nonholomorphic Eisenstein series for ${\rm PSL}(2,\bold Z)$, {\it Math.\ Phys.\ Electron.\ J}.
 {\bf 4} (1998), Paper 6, 8 pp.

[3]  \bibline
The transfer operator approach to Selberg's zeta function and modular and Maass wave forms for ${\rm PSL}(2,\bold Z)$,
in {\it  Emerging Applications of Number Theory}, 73--141, {\it IMA Vol\. Math\. Appl.\/} {\bf 109},  Springer, New York, 1999.

[4]  \name{Y.-J. Choie} and \name{D. Zagier},
Rational period functions for ${\rm PSL}(2,{\bf Z})$,
in {\it A Tribute to Emil Grosswald\/{\rm :} Number Theory and Related Analysis}, 
 {\it Contemp.\ Math.\ {\bf 143} (1993), 89--108.

[5]  \name{I. Efrat},
Dynamics of the continued fraction map and the spectral theory of ${\rm SL}(2,\bold Z)$, {\it Invent. Math. {\bf 114} (1993 \pages 207--218) 

[6]  A.~Erd\'elyi et.~al., Tables of Integral Transforms (Bateman manuscript project) Vols.~I, II
\publ McGraw-Hill \publaddr New York} (1954) 

[7]  W. Ferrar
Summation formulae and their relation to Dirichlet's series {\rm (II)}, {\it Comp.~Math. {\bf 4} (1937 \pages 394--405)

[8]  S. Helgason, Topics in Harmonic Analysis on Homogeneous Spaces
\publ Progress in Math.~{\bf13}, Birkh\"auser \publaddr Boston} (1981)

[9]  J.B. Lewis
Eigenfunctions on symmetric spaces with distribution-valued boundary forms, {\it J. Funct. Analysis{\bf 29\yr1978\pages287--307)

[10]  J.B. Lewis
Spaces of holomorphic functions equivalent to the even Maass cusp forms
\pages \jour Invent. Math.{\bf 127} (1997 \pages 271--306)

[11]  J.B. Lewis and D. Zagier
 Period functions and the Selberg zeta function for the modular group \pages 83-97
\inbook The Mathematical Beauty of Physics \publ Adv. Series in Math. Physics {\bf24}, World Scientific\publaddr Singapore} (1997)

[12]  H. Maa\ss \"Uber eine neue Art von nichtanalytischen automorphen Funktionen und
die Be\-stimmung Dirichletscher Reihen durch Funktionalgleichungen\jour Math\. Annalen {\bf 121} (1949\pages 141--183)

[13]  D. Mayer
The thermodynamic formalism approach to Selberg's zeta function for ${\rm PSL}(2,\bold Z)$, {\it Bull. AMS {\bf 25} (1991 \pages 55--60)

[14]  D. Mayer
Continued fractions and related transformations
\inbook Ergodic Theory, Symbolic Dynamics and Hyperbolic Spaces, T. Bedford, M. Keane, C. Series (Eds.)
 \publ Oxford University Press  \publaddr New York 1991 \pages 175--222)



\ref\key 15 ]  D. Zagier
Periods of modular forms and Jacobi theta functions, {\it Invent. Math. {\bf 104} (1991 \pages 449--465)

\endRefs 
\end